\documentclass[12pt,reqno]{amsart}

\usepackage{amscd}

\newcommand{\N}{{\mathbb N}}
\newcommand{\Z}{{\mathbb Z}}
\newcommand{\Q}{{\mathbb Q}}
\newcommand{\C}{{\mathbb C}}
\newcommand{\R}{{\mathbb R}}
\renewcommand{\P}{{\mathbb P}}
\renewcommand{\H}{{\mathbb H}}

\newcommand{\AAA}{{\mathcal A}}

\newcommand{\DD}{{\mathcal D}}

\newcommand{\HH}{{\mathcal H}}

\newcommand{\Nu}{{\mathcal V}}
\newcommand{\OO}{{\mathcal O}}

\newcommand{\QQ}{{\mathcal Q}}
\newcommand{\RR}{{\mathcal R}}
\newcommand{\TT}{{\mathcal T}}
\newcommand{\UU}{{\mathcal U}}
\newcommand{\XX}{{\mathcal X}}
\newcommand{\aaa}{{\bf a}}

\newcommand{\ddd}{{\rm d}}
\newcommand{\www}{\widetilde}
\newcommand{\oooo}{\overline}

\renewcommand{\Im}{{\rm Im}}
\renewcommand{\Re}{{\rm Re}}
\newcommand{\om}{{{\mathcal O}_M}}
\newcommand{\tm}{{\mathcal T}_M}
\newcommand{\tmne}{{\mathcal T}_M^{0,1}}
\newcommand{\tmen}{{\mathcal T}_M^{1,0}}
\newcommand{\am}{{\mathcal A}_M}
\newcommand{\ame}{{\mathcal A}_M^1}
\newcommand{\amz}{{\mathcal A}_M^2}
\newcommand{\amen}{{\mathcal A}_M^{1,0}}
\newcommand{\amne}{{\mathcal A}_M^{0,1}}
\newcommand{\amzn}{{\mathcal A}_M^{2,0}}
\newcommand{\amee}{{\mathcal A}_M^{1,1}}
\newcommand{\amnz}{{\mathcal A}_M^{0,2}}
\newcommand{\ampq}{{\mathcal A}_M^{p,q}}
\newcommand{\ciii}{C^\infty}
\newcommand{\ciik}{C^\infty(K)}
\newcommand{\chiii}{C^{h\infty}}
\newcommand{\cooo}{C^\omega}
\newcommand{\cook}{C^\omega(K)}
\newcommand{\chooo}{C^{h\omega}}
\newcommand{\eezz}{\frac{1}{z}}
\newcommand{\nmmm}{{\{0\}\times M}}
\newcommand{\zmmm}{{\{z\}\times M}}
\newcommand{\immm}{{\{\infty\}\times M}}
\newcommand{\cmmm}{{\C\times M}}
\newcommand{\csmmm}{{\C^*\times M}}
\newcommand{\pnmmm}{{(\P^1-\{0\})\times M}}
\newcommand{\pmmm}{{\P^1\times M}}
\newcommand{\dmmm}{{\Delta\times M}}

\newcommand{\oomm}{\Omega}
\newcommand{\paa}{\partial}
\newcommand{\dz}{{\partial_z}}
\newcommand{\zdz}{{z\partial_z}}
\newcommand{\dzz}{\frac{{\rm d}z}{z}}
\newcommand{\ndz}{\nabla_{\partial_z}}
\newcommand{\nzdz}{\nabla_{z\partial_z}}
\newcommand{\nnn}{\nabla}
\newcommand{\trrr}{\tau_{real}}

\newcommand{\Cz}{\C\{ z\}}
\newcommand{\vgii}{V^{>-\infty}}

\newcommand{\Exp}{{\rm Exp}}
\newcommand{\hiii}{H^\infty}
\newcommand{\hiie}{H^\infty_1}
\newcommand{\hiin}{H^\infty_{\neq 1}}
\newcommand{\hiil}{H^\infty_\lambda}
\newcommand{\hiia}{H^\infty_{e^{-2\pi i\alpha}}}
\newcommand{\hiib}{H^\infty_{e^{-2\pi i\beta}}}

\DeclareMathOperator{\Gl}{Gl}
\DeclareMathOperator{\Gr}{Gr}
\DeclareMathOperator{\Hom}{Hom}
\DeclareMathOperator{\id}{id}

\DeclareMathOperator{\Lie}{Lie}

\DeclareMathOperator{\rk}{rk}

\theoremstyle{plain}
\newtheorem{lemma}{Lemma}[section]
\newtheorem{theorem}[lemma]{Theorem}
\newtheorem{proposition}[lemma]{Proposition}
\newtheorem{corollary}[lemma]{Corollary}
\newtheorem{conjecture}[lemma]{Conjecture}

\theoremstyle{definition}
\newtheorem{definition}[lemma]{Definition}

\newtheorem{remark}[lemma]{Remark}
\newtheorem{remarks}[lemma]{Remarks}

\begin{document}

\title[$tt^{*}$ geometry, Frobenius manifolds, and singularities]
{$tt^{*}$ geometry, Frobenius manifolds, their connections, and
the construction for singularities} 

\author{
Claus Hertling}

%
%

\maketitle


\begin{abstract}
The base space of a semiuniversal unfolding of a hypersurface singularity
carries a rich geometry. By work of K. Saito and M. Saito it can be equipped
with the structure of a Frobenius manifold. By work of Cecotti and Vafa
it can be equipped with $tt^*$ geometry if the singularity is
quasihomogeneous. $tt^*$ geometry generalizes the notion of variation
of Hodge structures.

In the second part of this paper (chapters \ref{s6}--\ref{s8})
Frobenius manifolds and $tt^*$ geometry are constructed for any 
hypersurface singularity, using essentially oscillating integrals;
and the intimate relationship between polarized mixed Hodge structures
and this $tt^*$ geometry is worked out.

It builds on the first part (chapters \ref{s2}--\ref{s5}). There
$tt^*$ geometry and Frobenius manifolds and their relations are studied
in general. To both of them flat connections with poles are associated,
with distinctive common and different properties. A frame for a 
simultaneous construction is given.
\end{abstract}


\tableofcontents

\clearpage

\section{Introduction}\label{s1}
\setcounter{equation}{0}

\noindent
10 years ago Cecotti and Vafa \cite{CV1}\cite{CV2} considered moduli spaces
of $N=2$ supersymmetric quantum field theories and introduced
a geometry on them which is governed by the $tt^*$-equations.
This $tt^*$ geometry generalizes variations of Hodge structures.
Its most important ingredients are a real structure, given by a fiberwise
$\C$-antilinear involution $\kappa$, a hermitian metric $h$, and 
a Higgs field $C$.

A distinguished class of these moduli spaces is associated to 
quasihomogeneous singularities. They are the base spaces of semiuniversal
unfoldings of such singularities.

Already 10 years earlier K. Saito \cite{SK3}\cite{SK4} and M. Saito
\cite{SM2} introduced another geometry on these moduli spaces of 
singularities, the structure of Frobenius manifolds. This is a purely
holomorphic structure, but it shares many properties with the
$tt^*$ geometry. One can combine them to a structure which is called
CDV-structure in this paper because of the papers 
\cite{CV1}\cite{Du2}\cite{CV2}. It is defined below in definition 
\ref{t1.2}.

One purpose of this paper is to construct without physics CDV-structures
on the base spaces of semiuniversal unfoldings of all hypersurface
singularities, not only quasihomogeneous ones, and to study their
properties. This is done in chapters 6--8. It requires a careful
analysis of the common and different properties of $tt^*$ geometry and
Frobenius manifolds and of the connections with poles which belong to them.
This is carried out in chapters 2--5.

Both structures are studied first on abstract holomorphic vector bundles
$K\to M$, and only later on the holomorphic tangent bundle $TM$.
The $tt^*$ geometry on an abstract bundle is called CV-structure,
the part of a Frobenius manifold structure which makes sense is called
a Frobenius type structure. CV-structures make the contact to the work
of Simpson \cite{Si1}--\cite{Si4}. They are special cases of Simpson's
harmonic bundles. In a harmonic bundle one has a hermitian metric,
a Higgs field and the $tt^*$-equations, but not the real structure.
Simpson defined them as a generalization of complex variations of
Hodge structures and used them for his nonabelian Hodge theorem
\cite{Si3}. So his manifolds are not germs or small representatives
of germs, but compact K\"ahler manifolds.

A common feature of CV-structures and Frobenius type structures is that
there are two ways to describe them: in terms of geometry on a 
bundle $K\to M$ or in terms of a flat connection with poles
and a pairing on the lifted tangent bundle $\pi^* K$, where 
$\pi:\P^1\times M \to M$ is the projection.

The second point of view is most important. One can say that these
connections with poles present the true geometry in CV-structures
and Frobenius type structures. Several correspondences between 
connections with poles on $\pi^*K$ and geometry on $K$ are given
in chapters 2--5. They are quite intricate because the structures
on both sides are so rich. The correspondence for CV-structures
is richer than, but closely related to Simpson's correspondence
between harmonic bundles and variations of twistor structures
\cite{Si5}\cite{Sab6}.

The similarity of CV-structures and Frobenius type structures is visible
in the restriction of the connections to 
$\pi^*K|_\cmmm$, the difference in the extension to $\immm$.
On $\pi^*K|_\cmmm$ one has in both cases a flat connection with a pole
of Poincar\'e rank 1 along $\nmmm$ and a flat pairing $P$ on 
$\pi^*K|_\csmmm$ with certain properties. In the case of a CV-structure
the bundle $\pi^*K|_\csmmm$ has a real flat subbundle.
This structure is called a $(TERP)$-structure, here `T' stands for Twistor,
`E' for Extension of the connection in $z$-direction, so `TE' for a pole
of Poincar\'e rank 1, `R' for Real, `P' for Pairing.

Along $\immm$ the meromorphic connection of a Frobenius type structure has
a logarithmic pole, the connection for a CV-structure has
an `antiholomorphic twin' of a pole of Poincar\'e rank 1
(see definition \ref{t2.6} and theorem \ref{t2.19}).
There is freedom in the choice of the logarithmic pole, but the
antiholomorphic twin at $\immm$ of a pole of Poincar\'e rank 1 is obtained
uniquely with the real structure from the pole at $\nmmm$.

So $(TERP)$-structures are the common ferry to  CV-structures and to
Frobenius type structures. In singularity theory the Frobenius manifolds
arise from $(TERP)$-structures which come essentially from oscillating
integrals. The unique antiholomorphic extension of these 
$(TERP)$-structures thus induces a unique $tt^*$ geometry on these
Frobenius manifolds.

One may hope that $(TERP)$-structures and thus $tt^*$ geometry also exist
in quantum cohomology and in the Barannikov--Kontsevich construction
of Frobenius manifolds.

In definition \ref{t1.2} CDV-structures are defined using the 
connection with poles. First, the definition of Frobenius manifolds
is recalled.

\begin{definition}\label{t1.1}
A {\it Frobenius manifold} is a complex manifold $M$ with a commutative
and associative multiplication $\circ$ on the holomorphic tangent bundle $TM$,
with a unit field $e$, an {\it Euler field} $E$, and a holomorphic
metric $g$ such that the following holds. 
Let $C:\tm\to \oomm^1_M\otimes \tm$ be the Higgs field with
$C_XY:=-X\circ Y$. Let $\nnn^g$ be the Levi-Civita connection
of the metric $g$.

Then $\nnn^g$ is flat; the potentiality condition $\nnn^g(C)=0$ holds;
$e$ is flat; the Euler field satisfies $\Lie_E(\circ)=\circ$ and
$\Lie_E(g)=(2-d)\cdot g$ for some $d\in \C$; and
\begin{eqnarray} \label{1.1}
g(C_XY,Z)=g(Y,C_XZ) \mbox{ \ \ for }X,Y,Z\in \tm \ ,
\end{eqnarray} 
that means, $g$ is multiplication invariant.
\end{definition}

A meromorphic connection is associated to a Frobenius manifold as follows.
Define the endomorphisms $\UU$ and $\Nu$ of the holomorphic tangent bundle
$TM$,
\begin{eqnarray} \label{1.2}
\UU&:=&E\circ \ ,\\
\Nu &:=& \nnn^g_\bullet E -\frac{2-d}{2}\id = 
\nnn^g_E - \Lie_E -\frac{2-d}{2}\id\ , \label{1.3}
\end{eqnarray} 
and lift $\nnn^g, C,\UU,\Nu$ canonically to $\pi^*TM$, where $\pi:\pmmm\to M$.
Then it is well known (e.g. lemma \ref{t5.11}) that
\begin{eqnarray} \label{1.4}
\nnn:= \nnn^g +\eezz C+(\eezz\UU-\Nu)\dzz
\end{eqnarray} 
is a flat connection on $\pi^*TM|_\csmmm$. It is called (first) structure
connection of the Frobenius manifold. The flatness of $\nnn^g$ encodes
many properties in the definition of a Frobenius manifold.

\begin{definition}\label{t1.2}
A CDV-structure is a Frobenius manifold $(M,\circ,e,E,g)$ together with a
real structure on $TM$ given by a fiberwise $\C$-antilinear involution
$\kappa:T_tM\to T_tM$ such that the following holds.

$\kappa$ is an automorphism of $TM$ as a real analytic real vector bundle.
The form $h:=g(\cdot ,\kappa \cdot )$ is hermitian (but not necessarily
positive definitive) and satisfies
\begin{eqnarray} \label{1.5}
h(C_XY,Z)&=&h(Y,\kappa C_X\kappa Z) \mbox{ \ \ for }X,Y,Z\in \tm\ ,\\
\Lie_e(h)&=&0\ ,\label{1.6}\\
\Lie_{E-\oooo E}(h)&=&0 \label{1.7}\ .
\end{eqnarray} 
The metric connection $D$ on $TM$ for $h$ respects $\kappa$, i.e.
$D(h)=0$ and $D(\kappa)=0$. 
The number $d$ with $\Lie_E(g)=(2-d)\cdot g$ is real.
Define an endomorphism $\QQ$ of $TM$ as
real analytic complex vector bundle by
\begin{eqnarray} \label{1.8}
\QQ := D_E-\Lie_E-\frac{2-d}{2}\id \ .
\end{eqnarray} 
Lift $D$ and $\QQ$ canonically to $\pi^*TM$. Then 
\begin{eqnarray} \label{1.9}
\nnn^{CV} := D+\eezz C+z\kappa C\kappa + 
(\eezz \UU - \QQ - z\kappa \UU\kappa)\dzz 
\end{eqnarray} 
is a flat non holomorphic connection on $\pi^*TM|_\csmmm$.
\end{definition}

In \cite{He2} Frobenius manifold structures on the base spaces $M$ of
semiuniversal unfoldings $F$ of hypersurface singularities 
$f:(\C^{n+1},0)\to (\C,0)$ were constructed with the 
Gau{\ss}--Manin connection, following K. Saito and M. Saito.
In section \ref{s8.1} they are reconstructed using 
oscillating integrals and $(TERP)$-structures. A result of chapter 8 
is the following.

\begin{theorem}\label{t1.3}
Each of these Frobenius manifolds carries a canonical CDV-structure
outside of a real analytic subvariety $R\subset M$. The set $R$ is 
invariant under the flow of $E-\oooo E$. 
\end{theorem}

A CDV-structure with positive definite hermitian metric gives on 
the subvariety $\{t\in M\ |\ \UU|_t=0\}$ a variation of polarized
Hodge structures (chapter \ref{s3}). In the singularity case
this subvariety parametrizes the quasihomogeneous singularities.
The $\mu$-constant stratum is the larger set
$\{t\in M\ |\ \UU|_t \mbox{ is nilpotent }\}$. Here polarized
mixed Hodge structures turn up. In chapter \ref{s7} the relations between
CV-structures and polarized mixed Hodges structures are studied in 
general. In chapter \ref{s8} this is used to prove a part of the 
following conjecture.

\begin{conjecture}\label{t1.4}
Let $M$ be the base space of an unfolding $F$ of a singularity with 
a CDV-structure from chapter \ref{s8}. Let $t\in M$ be any point of $M$.

If one starts at $t$ and goes sufficiently far along the flow of the real
vector field $E+\oooo E$ then one does not meet $R$ anymore,
the metric $h$ is positive definite and the eigenvalues of $\QQ$ tend
to the set $\Exp (F_t)-\frac{n+1}{2}$.
\end{conjecture}

This conjecture is motivated by work in \cite{CV1}\cite{CV2}. There the 
behaviour of a CV-structure with respect to the renormalization group flow
is an important aspect. That flow corresponds to the vector
field $E+\oooo E$. The set $\Exp (F_t)$ is the union of the exponents
(the spectral numbers $+1$) of the singularities of $F_t$; that is a tuple
of $\mu$ numbers in $(0,n+1)\cap \Q$ which are symmetric around 
$\frac{n+1}{2}$ and which come from Hodge filtrations and monodromies.

The endomorphism $\QQ$ is perhaps the most fascinating piece of a 
CDV-structure. It is the `new supersymmetric index' from \cite{CFIV}.
If $h$ is positive definite then $\QQ$ is a hermitian
endomorphism, it is semisimple with real eigenvalues symmetric
around 0. The eigenvalues vary real analytically on $M$. But on the
set $\{t\in M\ |\ \UU|_t=0\}$ they are constant. There $e^{2\pi i \QQ}$
is an automorphism of the variation of Hodge structures, and the 
eigenspaces of $\QQ$ are subspaces of the Hodge decomposition.
In the singularity case there the set of eigenvalues of $\QQ$ is
$\Exp (F_t)-\frac{n+1}{2}$. That fits to the conjecture, because the set
$\{t\in M\ |\ \UU|_t=0\}$ is the fixed point set of $E$.
The following is the main result of chapter 8.

\begin{theorem}\label{t1.5}
Conjecture \ref{t1.4} is true for points $t\in M$ such that either
$\UU|_t$ is semisimple with $\mu$ different eigenvalues or $\UU|_t$ is 
nilpotent.
\end{theorem}

The semisimple case is a direct consequence of a result of Dubrovin
\cite{Du2}\cite{CV2} (cf. chapter \ref{s6}). The nilpotent case uses
polarized mixed Hodge structures and their relations to nilpotent orbits
(theorem \ref{t7.5}) and to Euler field orbits of $(TERP)$-structures
(theorem \ref{t7.20}). Associated to the $E+\oooo E$ orbit of a point 
$t$ in the $\mu$-constant stratum is an orbit 
$e^{rN}F^\bullet(t)$ of Hodge like filtrations.
Because the filtration $F(t)$ at $t$ is part of a polarized 
mixed Hodge structure they become Hodge filtrations of polarized
pure Hodge structures if one goes sufficiently far along 
$E+\oooo E$. This implies with some rather long technical arguments
(theorem \ref{t7.20}) that there $h$ is positive definite.

\bigskip
This paper grew out of the attempt to understand the physicists'
results in \cite{CV1}\cite{CV2} and other papers about singularities
and to recover them. This is only a first step, much more is in these 
papers, and many wishes and questions are open.

The first part of the paper benefitted a lot from Dubrovin's paper
\cite{Du2}.  The philosophy to generalize variations of Hodge structures
by considering connections with poles of Poincar\'e rank 1 is also
present in the work of Barannikov, the real structure is treated in
\cite[Remark 5.3]{Ba2}.

I am grateful to C. Sabbah for pointing me to the work of Simpson.
I thank V. Cort\'es, Yu. Manin, M. Rosellen, K. Saito, D. van Straten,
and A. Takahashi for interest in progress of this work 
and the Max-Planck-Institut f\"ur Mathematik for good working conditions.

\clearpage

\section{$tt^{*}$ geometry on abstract bundles}\label{s2}
\setcounter{equation}{0}

\noindent
In section \ref{s2.2} some basic notions in the work of Simpson are recalled,
Higgs bundles, $\lambda$-connections, and $(T)$-structures, which are
families of $\lambda$-connections.
In section \ref{s2.3} $(\www T)$-structures, the `antiholomorphic twins' of 
$(T)$-structures, are defined. This allows to give in section \ref{s2.4}
Simpson's notion of variation of twistor structures and to recover in
lemma \ref{t2.11} the correspondence between them and harmonic bundles
without metric ($(DC\www C)$-structures).
In section \ref{s2.5} and section \ref{s2.6} the structures on both sides
of the correspondence are enriched. The notion of $(TERP)$-structures
is central in the whole paper. Finally, in section \ref{s2.7}
the correspondence between these enriched structures is presented.

\subsection{Some notations}\label{s2.1}

Throughout the whole paper, $M$ will denote a complex manifold
of dimension $m$.
The sheaf of holomorphic $k$-forms is $\oomm_M^k$, the sheaf of 
$\ciii$ $(p,q)$-forms is $\ampq$, and $\AAA_M^k:=\bigoplus_{p+q=k}
\ampq$. Especially $\OO_M=\oomm_M^0$ and $\ciii_M=\am^0$.
The sheaf of holomorphic vector fields is $\tm$, the sheaf of
$\ciii$-vector fields is $\tm^\C = \tmen \oplus \tmne$, with
$\tmen=\tm \otimes_\om \ciii_M$, and $\tmne = \oooo{\tmen}$.

Vector fields denoted by latin letters $($X,Y,Z$)$ will always mean
vector fields in $\tmen$, sometimes even in $\tm$.
For vector fields not necessarily of type $(1,0)$ greek letters will be used,
$\xi,\xi_1,\xi_2\in \tm^\C$.

If not said otherwise, vector bundles will be complex 
$\ciii$-bundles, sometimes with additional structure (real analytic,
holomorphic, antiholomorphic).

Let $K\to M$ be a vector bundle with sheaf $\ciik$ of $\ciii$-sections.
A connection is a map $D:\ciik\to \ame\otimes \ciik$ which satisfies
the Leibniz rule, $D(a\cdot s)= \ddd a\otimes s + a\cdot Ds$ 
for $a\in \ciii_M,\ s\in \ciik$.
With the extension
$$D: \ame\otimes \ciik \to \amz\otimes \ciik,\ 
\omega\otimes s\mapsto \ddd\omega\otimes s - \omega\land D s
$$
the curvature of $D$ is the map
$D^2 : \ciik \to \amz\otimes \ciik$. The connection is flat if 
$D^2=0$.

A $(1,0)$-connection is a map $D':\ciik \to \amen\otimes\ciik$ 
which satisfies the Leibniz rule for $(1,0)$-vector fields,
$D'(a\cdot s)= \paa a\otimes s + a\cdot D's$. 
It is extended to 
$$D': \ame\otimes \ciik \to \amz\otimes \ciik,\ 
\omega\otimes s\mapsto \paa\omega\otimes s - \omega\land D' s\ .
$$
It is called flat if ${D'}^2=0$.
A $(0,1)$-connection is defined analogously. A connection $D$ 
is a sum $D= D^{(1,0)}+D^{(0,1)}$ of a $(1,0)$-connection and 
a $(0,1)$-connection.

Giving a holomorphic structure on $K\to M$ is equivalent to giving
a flat $(0,1)$-connection: If a holomorphic structure with sheaf
$\OO(K)$ of holomorphic sections is given one extends the map
$D'':\OO(K)\to 0$ with the Leibniz rule to a flat $(0,1)$-connection.
If a flat $(0,1)$-connection $D''$ is given then, because of the 
Newlander-Nirenberg theorem, $\ker D''$ is the sheaf of 
holomorphic sections of a holomorphic structure on $K\to M$.

A $\ciii_M$-linear map $C:\ciik \to \ame\otimes \ciik$ has the 
canonical extension
$$C: \ame\otimes \ciik \to \amz\otimes \ciik,\ 
\omega\otimes s\mapsto - \omega\land C s
$$

For operators $D, D', D'', C$ as above,the compositions $D^2, \ {D'}^2, 
\ {D''}^2,\ D'D'' + D''D',\ D(C):=DC+CD,\ D'(C):=D'C+CD',\ D''(C):=D''C+CD''$ 
and $C^2$ are all maps from $\ciik$ to $\amz\otimes \ciik$.
Compatibility conditions will be written subsequently as vanishing
of such conditions. Therefore it is useful to rewrite these
compositions after inserting vector fields 
$X,Y\in \tmen$ or $\xi_1,\xi_2\in \tm^\C$. One easily sees
\begin{eqnarray} 
D^2(\xi_1,\xi_2) &=& [D_{\xi_1},D_{\xi_2}]-D_{[\xi_1,\xi_2]}\ ,
\label{2.1}\\
{D'}^2(X,Y) &=& [D'_X,D'_Y]-D'_{[X,Y]}\ ,
\label{2.2}\\
{D''}^2(\oooo X,\oooo Y) &=& [D''_{\oooo X},D''_{\oooo Y}]-
D''_{[\oooo X,\oooo Y]}\ , \label{2.3}\\
(D'D'' + D''D')(X,\oooo Y) &=& [D'_{X},{D''}_{\oooo Y}]-
(D'+D'')_{[X,\oooo Y]}\ ,    \label{2.4}\\
D(C)(\xi_1,\xi_2) &=& D_{\xi_1}(C_{\xi_2})- D_{\xi_1}(C_{\xi_2})
-C_{[\xi_1,\xi_2]}\ ,\label{2.5}\\
D'(C)(X,Y) &=& D'_X(C_Y)-D'_Y(C_X)-C_{[X,Y]}\ ,\label{2.6}\\
D''(C)(\oooo X,Y) &=& D''_{\oooo X}(C_Y)
\mbox{ \ \ for } X,Y\in \tm\ ,\label{2.7}\\
C^2(\xi_1,\xi_2) &=& C_{\xi_1}C_{\xi_2}- C_{\xi_2}C_{\xi_1}\ ,\label{2.8}
\end{eqnarray}
and if $C$ and $\www C$ are two $\ciii_M$-linear maps as above then
\begin{eqnarray}
(C\www C+\www C C)(\xi_1,\xi_2)&=& [C_{\xi_1},\www C_{\xi_2}]
   - [C_{\xi_2},\www C_{\xi_1}]\ .\label{2.9}
\end{eqnarray}

\subsection{Higgs bundles, $\lambda$-connections and $(T)$-structures}
\label{s2.2}

\begin{definition}\label{t2.1}\cite{Si4}
A {\it Higgs bundle} is a holomorphic vector bundle $K\to M$ together
with an $\OO_M$-linear map $C:\OO(K)\to \oomm_M^1\otimes \OO(K)$ which
satisfies
\begin{eqnarray} \label{2.10}
C_XC_Y=C_YC_X \mbox{ \ \  for } X,Y\in \tm\ ,
\end{eqnarray}
that is, $C^2=0$. The map $C$ is called a {\it Higgs field}.
\end{definition}

Equivalently, a Higgs bundle is a triple $(K,C,D'')$ where $K\to M$ is 
a complex $\ciii$-vector bundle, $C$ is a $\ciii_M$-linear  map 
$C:\ciik\to \amen\otimes\ciik$, and $D''$ is a $(0,1)$-connection,
which satisfy $(D''+C)^2=0$, that is,
\begin{eqnarray} \label{2.11}
{D''}^2=0,\ D''(C)=0,\ C^2=0\ .
\end{eqnarray}
$D''$ gives a holomorphic structure on $K\to M$ with sheaf of holomorphic
sections $\OO(K)$, and $D''(C)=0$ says that 
$C_X:\OO(K)\to \OO(K)$ for $X\in \tm$.

The notion of a $\lambda$-connection generalizes the notions of a 
holomorphic connection and of a Higgs bundle.
It is due to Deligne \cite{Si3}.

\begin{definition}\label{t2.2}
Fix $\lambda\in \C$ and a holomorphic vector bundle $K\to M$.
A {\it $\lambda$-connection} is a map $D:\OO(K)\to \oomm_M^1\otimes \OO(K)$
which satisfies $D(a\cdot s)= \lambda\ddd a\otimes s +a\cdot Ds$.
It has the natural extension
\begin{eqnarray} \label{2.12}
D:\oomm_M^1\otimes \OO(K)\to \oomm_M^2\otimes \OO(K),\ 
\omega\otimes s \mapsto \lambda \ddd \omega \otimes s - \omega\land Ds\ .
\end{eqnarray}
A $\lambda$-connection is {\it flat} if $D^2=0$.
\end{definition}

For $\lambda\neq 0$, $D$ is a (flat) $\lambda$-connection if and only
if $\frac{1}{\lambda}D$ is a (flat) holomorphic connection.
A Higgs bundle is a flat $0$-connection.

A holomorphic family of flat $\lambda$-connections, $\lambda\in \C$, is 
equivalent to the following structure. It will be a basic playing
character and will be equipped with additional structure in section
\ref{s2.5}.

\begin{definition}\label{t2.3}
Let $H\to \C\times M$ be a holomorphic vector bundle.
A $(T)$-structure is a pair $(H,\nabla)$. Here $\nabla$ is a map
\begin{eqnarray} \label{2.13}
\nabla:\OO(H)\to \eezz \OO_{\C\times M}\cdot \oomm_M^1\otimes \OO(H)
\end{eqnarray}
such that for any $z\in \C^*$ the restriction of $D$ to 
$H|_\zmmm$ is a flat connection
('T' as in 'Twistor' \cite{Si5} or in 'Topological' \cite{CV1}).
\end{definition}

\begin{lemma}\label{t2.4}
Let $(H,\nabla)$ be a $(T)$-structure. Define $K:=H|_\nmmm$. The map 
$z\cdot \nabla : \OO(H)\to \OO_{\C\times M}\cdot \oomm_M^1\otimes
\OO(H)$ restricts to a flat $z$-connection $\nabla^{(z)}$ on 
$H|_\zmmm$ for any $z\in \C$ (including $z=0$).
The holomorphic bundle $K$ together with $C:= \nabla^{(0)}$ is a Higgs
bundle, the Higgs bundle of the $(T)$-structure.
In view of $\OO(K)= \OO(H)/z\OO(H)|_\nmmm$ one has for $X\in \tm,
\ a \in \OO(H),[a]\in\OO(K)$
\begin{eqnarray} \label{2.14}
C_X ([a]) = (z\cdot \nabla_X (a))\mod z\OO(H)\ .
\end{eqnarray}
\end{lemma}

{\it Proof.} trivial. \hfill $\qed$

\begin{remarks}\label{t2.5}
(i) The support of a Higgs bundle $(K\to M,C)$ is the set 
$$
\bigcup_{t\in M}\{\lambda\in T^*_tM\ |\ \forall\ X\in T_tM\ 
\ker (C_X-\lambda(X)\id:T_tM\to T_tM)\neq 0\}
$$
of simultaneous eigenvalues of all $C_X$ at all $t\in M$.
It is a subvariety of dimension $m$ in the holomorphic cotangent bundle
$T^*M$, and it is finite over $M$.

(ii) The Higgs bundle of a $(T)$-structure has the special property that
its support is Lagrange in $T^*M$. This follows from 
\cite[Proposition 1.25]{Sab6} and the fact that in the notation of
\cite{Sab6} a $(T)$-structure is a locally free $\OO_\XX$-module
and a strict holonomic $\RR_\XX$-module.

(iii) Usually a $(T)$-structure turns up as part of a meromorphic 
connection $\nabla$ on $H$ with a pole of Poincar\'e rank 1 
along $\nmmm$, see definition \ref{t2.12}.

(iv) In \eqref{2.13} $\oomm_M^1$ is written for $\pi^{-1}\oomm_M^1$, 
where $\pi:\C\times M \to M$.
Similar notations will be used in the formulas \eqref{2.15}, \eqref{2.16},
\eqref{2.29}.
Because of these formulas we write $\OO_{\C\times M}\cdot \oomm_M^1$
instead of $\pi^*\oomm_M^1$ or $\oomm^1_{\C\times M/M}$.
\end{remarks}

\subsection{$C^\infty$-considerations}\label{s2.3}

In order to define an antiholomorphic twin of a $(T)$-structure, one has
to start with a mixture of $\ciii$ and holomorphic data.

For $U\subset \P^1$ open let $\OO_U\ciii_M$ be the sheaf of 
$\ciii$-functions on $U\times M$ which are holomorphic in $z$-direction
($z$ is the coordinate on $\C$).
A complex $\ciii$-vector bundle $H\to U\times M$ is said to 
carry a $C^{h\infty}$ structure if in the sheaf $\ciii(H)$ of all
$\ciii$-sections a subsheaf $C^{h\infty}(H)$ is fixed which is a locally
free $\OO_U\ciii_M$-module of rank $\rk H$. Such a bundle $H$ can be 
seen as a $\ciii$-family over $M$ of holomorphic vector bundles on $U$.

It will also be useful to consider everything with $\ciii$ replaced
by real analytic; we will simply replace $\ciii$ and $\chiii$
by $C^\omega $ and $C^{h\omega}$.

Now a $(T)$-structure is a $\chiii$-vector bundle $H\to \C\times M$ 
with a map
\begin{eqnarray} \label{2.15}
\nabla: \chiii (H)\to (\eezz\OO_\C\ciii_M\cdot \amen + 
\OO_\C\ciii_M\cdot \amne) \otimes \chiii(H)
\end{eqnarray}
such that for any $z\in \C^*$ the restriction of $\nabla$ to 
$H|_\zmmm$ is a flat connection.
The $(0,1)$-part of $\nabla$ together with the $\chiii$-structure
on $H$ induces a holomorphic structure on $H$. An
antiholomorphic twin of a $(T)$-structure is defined as follows.

\begin{definition}\label{t2.6}
A {\it $(\www T)$-structure} is a pair $(\www H,\www \nabla )$. 
Here $\www H\to (\P^1-\{0\})\times M$ is a 
$\chiii$-vector bundle and $\www \nnn$ is a map
\begin{eqnarray} \label{2.16}
 \www \nnn : \chiii(H) \to (\OO_{\P^1-\{0\}}\ciii_M\cdot \amen +
z\cdot \OO_{\P^1-\{0\}}\ciii_M\cdot \amne)\otimes \chiii(H)
\end{eqnarray}
such that for any $z\in \C^*$ the restriction of $\www \nabla$ to
$H|_\zmmm$ is a flat connection.
\end{definition}

There is a natural correspondence between $(T)$-structures and
$(\www T)$-structures. One needs two operations,
the $\C$-antilinear map
\begin{eqnarray} \label{2.17}
\gamma:\P^1\to \P^1,\ z\mapsto \frac{1}{\oooo z}\ ,
\end{eqnarray}
and complex conjugation of vector bundles.
The complex conjugate $\oooo H\to U\times M$ of a vector bundle
$H\to U\times M$ is the same as $H$ except that the fibers are
equipped with the complex conjugate $\C$-linear structure.

\begin{lemma}\label{t2.7}
If $(H,\nabla)$ is a $(T)$-structure (respectively a $(\www T)$-structure)
then $\oooo{\gamma^*H}=:\www H$ carries a canonical $(\www T)$-structure
(respectively a $(T)$-structure).
This gives a one-to-one correspondence between $(T)$-structures
and $(\www T)$-structures.
\end{lemma}

{\it Proof.}
The equation $\oooo{\gamma^*(\eezz \OO_\C\ciii_M)}
=z\OO_{\P^1-\{0\}}\ciii_M$ shows that $\www H$ is a $\chiii$-vector bundle
with $\chiii(\www H)= \oooo{\gamma^*(\chiii (H))}$ and that the induced
map $\www \nabla$ on $\chiii(\www H)$ satisfies \eqref{2.16}.
\hfill $\qed$

\subsection{Variation of twistor structures}\label{s2.4}

Simpson \cite[ch. 3]{Si5} defined the notion of a variation of 
twistor structures. Sabbah called a variation of twistor structures 
a smooth twistor structure \cite[2.2.a]{Sab6} and greatly generalized this
notion.

\begin{definition}\label{t2.8}
(a) A {\it $(T\www T)$-structure} is a $\chiii$-vector bundle 
$H\to \pmmm$ and a map $\nabla$ such that $(H,\nabla)|_\cmmm$ is a 
$(T)$-structure and $(H,\nabla)|_\pnmmm$ is a $(\www T)$ structure.

(b) \cite[ch. 3]{Si5}\cite[2.2.a]{Sab6} A {\it variation of twistor 
structures }is a $(T\www T)$-structure such that for each 
$t\in M$ the restriction $H|_{\P^1\times \{t\}}$ is a trivial bundle.
We call it also a $(trT\www T)$-structure (`tr' for trivial).
\end{definition}

A $(T)$-structure $(H,\nabla)$ equips $K:=H|_\nmmm$ with the structure
of a Higgs bundle, that is, a flat $(0,1)$-connection $D''$ and a
$\ciii_M$-linear map
\begin{eqnarray}  \label{2.18}
C : \ciik\to \amen \otimes \ciik
\end{eqnarray}
with $(D''+ C)^2=0$. 
A $(\www T)$-structure $(H,\nnn)$ equips $H|_\immm$ with the antiholomorphic
twin of a Higgs bundle, which is by definition a flat 
$(1,0)$-connection $D'$ and a $\ciii_M$-linear map
\begin{eqnarray}  \label{2.19}
\www C : \ciii (H|_\immm)\to \amne \otimes \ciii (H|_\immm)
\end{eqnarray}
with $(D'+\www C)^2=0$. 
Here $\www C := \lim_{z\to \infty} (\eezz\nabla)|_\zmmm$.
In the case of a $(trT\www T)$-structure $K$ and $H|_\immm$ are canonically
isomorphic and thus all four operators $D'',\ C,\ D',\ \www C$ live
on $K$. Their properties are stated in definition \ref{t2.9} and
lemma \ref{t2.11}.

\begin{definition}\label{t2.9}
A $(DC\www C)$-structure is a $\ciii$ vector bundle $K\to M$ together
with a connection $D$ on it and two $\ciii_M$-linear maps
\begin{eqnarray} 
C:\ciik \to \amen \otimes \ciik \label{2.20}\\
\www C:\ciik \to \amne \otimes \ciik \label{2.21}
\end{eqnarray} 
with the following properties. Let $D'$ and $D''$ be the $(1,0)$-part
and the $(0,1)$-part of $D$. Then
\begin{eqnarray} 
(D''+C)^2=0,\ \ (D'+\www C)^2=0,\label{2.22}\\
D'(C)=0,\ \ D''(\www C)=0,\label{2.23}\\
D'D''+D''D' = -(C\www C + \www C C)\ . \label{2.24}
\end{eqnarray} 
\end{definition}

\begin{remarks}\label{t2.10}
(i) In formula \eqref{2.11} it was already stated that $(D''+C)^2=0$ 
decomposes into ${D''}^2=0,\ D''(C)=0,\ C^2=0$. An analogous statement
holds for $(D'+\www C)^2=0$.

(ii) The equations \eqref{2.23} and \eqref{2.24} are called $tt^*$-equations
in \cite{CV1}.
The equation $D'(C)=0$ says (cf. \eqref{2.6})
\begin{eqnarray} \label{2.25}
D'_X(C_Y)-D'_Y(C_X) = C_{[X,Y]} \mbox{ \ \ for } X,Y\in \tmen\ .
\end{eqnarray} 
We call $D'(C)=0$ potentiality condition, because it turns up in the
case of a Frobenius manifold as one formulation of potentiality.
The curvature condition \eqref{2.24} can be written as (cf. \eqref{2.4}
and\eqref{2.9}) 
\begin{eqnarray} 
[D'_X,D''_{\oooo Y}]-(D'+D'')_{[X,\oooo Y]} = -[C_X,{\www C}_{\oooo Y}]
\mbox{ \ \ for } X,Y\in \tmen\ .
\end{eqnarray} 
\end{remarks}

\begin{lemma}\label{t2.11}\cite[Lemma 3.1]{Si5}\cite[Lemma 2.2.2]{Sab6}
There is a one-to-one correspondence between variations of twistor structures
(i.e. $(trT\www T)$-structures) and $(DC\www C)$-structures. It is given
by the steps in (a) and (b). They are inverse to one another.

(a) Let $(K\to M,D,C,\www C)$ be a $(DC\www C)$-structure.
Let $\pi: \pmmm \to M$ be the projection. The lifted bundle
$H:=\pi^*K$ on $\pmmm$ is a $\chiii$-bundle. The operators
$D,C,\www C$ are lifted canonically to $H$. Then the bundle $H$ and the 
operator
\begin{eqnarray} \label{2.27}
\nabla := D+\eezz C+z\www C
\end{eqnarray} 
form a variation of twistor structure.

(b) In the case of a variation of twistor structures $(H,\nabla)$ the four
operators $D'',C,D',\www C$ on $K:=H|_\nmmm$ defined before definition
\ref{t2.9} yield a $(DC\www C)$-structure.
\end{lemma}

{\it Proof.} (a) The only thing to be shown is $\nabla^2=0$. 
This is equivalent to \eqref{2.22}--\eqref{2.24}.

(b) The operators $D'', C, D', \www C$ are lifted canonically to $H$,
using the canonical isomorphism $\pi^*K\to H$, respectively 
the isomorphism $\ciik \cong \pi_*\chiii (H)$.
One has to show
\begin{eqnarray} \label{2.28}
\nabla = D'+D''+ \eezz C+z\www C\ .
\end{eqnarray} 
Consider the operator $\nabla_X - (D'_X+D''_X+\eezz C_X+z\www C_X) = 
\nabla_X-D'_X-\eezz C_X$ for $X\in \tmen$.
It maps sections in $\pi_*\chiii (H)$ to $\chiii$-sections which have
no pole along $\nmmm$ because of the definition of $C$ and which vanish
along $\immm$ because of the definition of $D'$. Therefore they
vanish globally.
An analogous statement holds for $\oooo X$.
Therefore \eqref{2.28} holds. The flatness $\nabla^2=0$ gives the conditions
\eqref{2.22}--\eqref{2.24}. 
\hfill $\qed$

\subsection{$(TERP)$-structures}\label{s2.5}

We will now equip both sides in the correspondence in lemma \ref{t2.11}
with additional structures. In definition \ref{2.12} a $(T)$-structure
will be enriched to a $(TERP(w))$-structure. Here `E' stands for
Extension of $\nabla$ in $z$-direction, `R' stands for Real structure,
`P' stands for Pairing, and $w$ will be an integer.

Such a structure induces a $(T\www T)$-structure and additional geometry
on $K:=H|_\nmmm$ (lemma \ref{2.14}).  If the $(T\www T)$-structure is a
$(trT\www T)$-structure then one obtains a rich geometric structure on
$K$, a CV-structure (definition \ref{t2.16}). 
Theorem \ref{t2.19} extends lemma \ref{2.11} to an equivalence between
CV-structures and $(trTERP(w))$-structures.

They are really two sides of the same coin. But the formulation of the 
properties of real structure and pairings on both sides is rather lengthy,
and passing from one side to the other side is nontrivial.

$(TERP(w))$-structures turn up in \cite{CV1}\cite{CFIV}\cite{CV2} and 
in singularity theory (see chapter \ref{s8}). 
The discussion here and the equivalence in theorem \ref{t2.19} 
are obtained by unifying considerations in
\cite{CV1}\cite{CFIV}\cite{Du2}\cite{Si4}\cite{Si5}.

\begin{definition}\label{t2.12}
A $(TERP(w))$-structure is a tuple $(H,\nabla,H_\R,P)$ with the following
properties.

(a) The space $H$ is a holomorphic vector bundle $H\to \cmmm$.

(b) It is equipped with a flat holomorphic connection $\nabla$ on 
$H|_\csmmm$ with a pole of Poincar\'e rank 1 along $\nmmm$, that is,
\begin{eqnarray} \label{2.29}
\nabla :\OO(H) &\to& (\eezz\OO_\cmmm\cdot\oomm^1_M + 
\eezz\OO_\cmmm \dzz)\otimes \OO(H)\\
&& = \eezz \oomm^1_{\C\times M}(\log (\nmmm ))\otimes \OO(H) \nonumber 
\end{eqnarray} 
and $\nabla^2=0$.

(c) The bundle $H|_\csmmm$ contains a real bundle $H_\R\to \csmmm$ 
which is $\nabla$-flat and satisfies 
$H_{(z,t)}=(H_\R)_{(z,t)}\oplus i(H_\R)_{(z,t)}$.

(d) $w\in \Z$. There is a $\C$-bilinear pairing
\begin{eqnarray} \label{2.30a}
P:H_{(z,t)}\times H_{(-z,t)}\to \C \mbox{ \ \ for any }(z,t)\in\csmmm
\end{eqnarray} 
with the following properties. 
It is nondegenerate and $(-1)^w$-symmetric (symmetric for even $w$,
antisymmetric for odd $w$). It extends for an open subset 
$U_1\times U_2\subset \C\times M$ to a nondegenerate pairing
\begin{eqnarray} \label{2.30b}
\qquad P:\OO(H)(U_1\times U_2)\times \OO(H)((-U_1)\times U_2) \to 
z^w\OO_\cmmm (U_1\times U_2) \\
(a, b) \mapsto ((z,t)\mapsto P(a(z,t),b(-z,t)))\ .\nonumber
\end{eqnarray} 
Thus, $P$ is $\OO_M$-bilinear, but $z$-sesquilinear,
\begin{eqnarray} \label{2.31}
f(z,t)P(a,b) = P(f(z,t)a,b) = P(a,f(-z,t)b) \mbox{ \ \ for }
f\in \OO_\cmmm \ .
\end{eqnarray} 
It is $\nabla$-flat, that is, $\nabla(P)=0$, or explicitely for 
$X\in \tm(U_2),\ a\in \OO(H)(U_1\times U_2),\ b\in \OO(H)((-U_1)\times U_2)$,
\begin{eqnarray} \label{2.32}
X\, P(a,b) &=& P(\nabla_Xa,b) + P(a,\nabla_Xb)\ ,\\
z\dz \, P(a,b) &=& P(\nabla_\zdz a,b) + P(a,\nabla_\zdz b)\ .\label{2.33}
\end{eqnarray} 
The pairing $P$ takes values in $i^w\R$ on $H_\R$,
\begin{eqnarray} \label{2.34}
P:(H_\R)_{(z,t)}\times (H_\R)_{(-z,t)}\to i^w\R\ .
\end{eqnarray} 
\end{definition}

\begin{remarks}\label{t2.13}
(i) The definition of weaker structures is obvious, for example\\
$(TE)$: $(H,\nnn)$ with (a) and (b),\\
$(TER)$: $(H,\nnn,H_\R)$ with (a), (b), and (c),\\
$(TEP(w)$: $(H,\nnn,P)$ with (a), (b), (d) except for \eqref{2.34},\\
$(TP(w))$: a $(T)$-structure $(H,\nnn)$ with a pairing $P$ as in (d) except
for \eqref{2.33}, \eqref{2.34}.

In \cite{Si5} and \cite{Sab6} $(T)$-structures with hermitian pairings
$H_{(z,t)}\times H_{(-\gamma(z),t)}\to \C$ for $(z,t)\in \csmmm$ are 
considered ($\gamma:\P^1\to \P^1$ was defined in \eqref{2.17}).
Here such a pairing arises as $P(\cdot,\tau\cdot)$ for $\tau$ as in 
lemma \ref{t2.14} (d).

(ii) In chapter \ref{s8} $(TERP(w))$-structures will be discussed which come
from oscillating integrals for isolated singularities.
$(TERP(w))$-structures turn also up in \cite{CV1}.
The first structure connections of any holomorphic Frobenius manifold
(without Euler field) contain $(TEP(w))$-structures 
(respectively $(TP(w))$-structures), see section \ref{s5.3}.

Establishing such structures often is the first and most difficult part
in the construction of Frobenius manifolds. In quantum cohomology
and in the Barannikov--Kontsevich construction in most cases only formal
versions of $(TEP(w))$-structures are established. But one may expect
that they are convergent and carry real structures as in (c).
\end{remarks}

\begin{lemma}\label{t2.14}
Let $(H,\nnn,H_\R,P)$ be a $(TERP(w))$-structure. Then
$K:=H|_\nmmm$ together with the tensor $C$ in lemma \ref{t2.4} is a 
Higgs bundle. There is additional structure with the following 
properties.

(a) There is a symmetric nondegenerate holomorphic pairing $g$ on $K$,
\begin{eqnarray} 
g &:& K_t\times K_t \to \C \mbox{ \ \ for }t\in M\ , \nonumber\\
g &:& \OO(K)\times \OO(K) \to\OO_M \ \ \OO_M\mbox{-bilinear}\ ,\label{2.35}
\end{eqnarray} 
defined by
\begin{eqnarray} \label{2.36}
g([a],[b]):=(z^{-w}P(a,b))\mod z\OO(H)
\end{eqnarray} 
for $[a],[b]\in \OO(K)$ and lifts $a,b\in \OO(H)$ (in suitable open sets).
It satisfies
\begin{eqnarray} \label{2.37}
g(C_X[a],[b])=g([a],C_X[b])\mbox{ \ \ for }X\in \tm\ .
\end{eqnarray} 

(b) There is a holomorphic endomorphism $\UU$ of $K$,
\begin{eqnarray} 
\UU &:& K_t\to K_t \mbox{ \ \ for }t\in \UU\ ,\nonumber\\
\UU &:& \OO(K) \to \OO(K) \ \ \OO_M\mbox{-linear}\ ,\label{2.38}
\end{eqnarray} 
defined by
\begin{eqnarray} \label{2.39}
\UU([a]):= (z^2\nnn_\dz (a))\mod z\OO(H)
\end{eqnarray} 
for $[a]\in \OO(K),\ a\in \OO(H)$. It satisfies
\begin{eqnarray} \label{2.40}
&& \UU C_X = C_X\UU \mbox{ \ \ for }X\in \tm\ ,\\
&& g(\UU[a],[b]) = g([a],\UU[b]) \mbox{ \ \ for }[a],[b]\in \OO(K)\ .
\label{2.41}
\end{eqnarray} 

(c) Let $\kappa_H:H|_\csmmm \to H|_\csmmm$ be the fiberwise 
$\C$-antilinear involution with $\kappa_H|_{H_\R}=\id$.
The points $z\in \C^*$ and $\gamma(z):=\frac{1}{\oooo z}$ (cf. \eqref{2.17})
are contained in the real halfline
$\{\zeta\in\C^*\ |\ \arg\zeta=\arg z\}$. 
For $(z,t)\in \csmmm$ let 
$\gamma_\nnn : H_{(z,t)}\to H_{(\gamma(z),t)}$ 
be the isomorphism which one obtains from the $\nnn$-flat shift along this 
half line. Then the map
\begin{eqnarray} \label{2.42}
 \trrr := \kappa_H \circ \gamma_\nnn = \gamma_\nnn \circ \kappa_H : 
H_{(z,t)}\to H_{(\gamma (z),t)} \mbox{ \ \ for }(z,t)\in \csmmm
\end{eqnarray} 
is $\C$-antilinear and $\nnn$-flat and satisfies $\trrr^2=\id$ and
$\trrr (z\cdot a)=\eezz\trrr (a)$ (this last equation has to be understood
with $z$ as a function on $\C^*$, not pointwise).

(d) The map
\begin{eqnarray} \label{2.43}
\tau:H_{(z,t)}&\to& H_{(\gamma (z),t)}\mbox{ \ \ for }(z,t)\in \csmmm\\
a&\mapsto & \trrr (z^{-w}a)\nonumber
\end{eqnarray} 
is $\C$-antilinear and satisfies $\tau^2=\id$, $\tau(z\cdot a)=\eezz\tau(a)$,
and 
\begin{eqnarray} \label{2.44}
\nnn_\xi \circ \tau = \tau \circ \nnn_{\oooo \xi} \mbox{ and }
\nnn_\zdz \circ \tau = \tau \circ (\nnn_{-\zdz}+w\cdot \id)\ .
\end{eqnarray} 
It is an isomorphism
\begin{eqnarray} \label{2.45}
\tau: H|_\csmmm \to \oooo{(\gamma^* H)|_\csmmm} \ ,
\end{eqnarray} 
which respects the $\chiii$-structures and the restrictions of the flat
connections to the slices $\zmmm$ for $z\in \C^*$.

(e) The bundles $H$ and $\oooo{\gamma^*H}$ are glued with $\tau$ to a bundle
$\hat H$. Then $(\hat H,\nnn)$ is a $(T\www T)$-structure.
\end{lemma}

{\it Proof.}
(a) The pairing $g$ is symmetric because $P$ is $(-1)^w$-symmetric and
$z$-sesquilinear. The pairing $g$ is nondegenerate and holomorphic 
because of \eqref{2.30b}. It satisfies \eqref{2.37}, because for
$X\in \tm,\ a,b,\in \OO(H),\ [a],[b]\in \OO(K)$ one has
\begin{eqnarray} 
0& = & zX\, z^{-w}P(a,b) \mod z\OO_\cmmm\nonumber \\
& = & z^{-w}P(z\nnn_X a,b) - z^{-w}P(a,z\nnn_X b)
\mod z\OO_\cmmm\label{2.46}\\
&=& g(C_X[a],[b])-g([a],C_X[b])    \ .\nonumber 
\end{eqnarray} 

(b) One has for $X\in \tm,\ a,b,\in \OO(H),\ [a],[b]\in \OO(K)$
\begin{eqnarray} 
0 &=& ([\nnn_{z^2\dz},\nnn_{zX}]-\nnn_{z^2X})a \nonumber \\
& = & [\nnn_{z^2\dz},\nnn_{zX}] a \mod z\OO(H) \label{2.47}\\
&=& (\UU C_X - C_X\UU)[a]  \nonumber
\end{eqnarray} 
and
\begin{eqnarray} 
0 & = & z^2\dz \, z^{-w}P(a,b) \mod z\OO_\cmmm \nonumber \\
& = & (-w)z\cdot z^{-w}P(a,b)+ z^{-w}P(\nnn_{z^2\dz}a,b) 
- z^{-w}P(a,\nnn_{z^2\dz}b)\nonumber\\
&& \mod z\OO_\cmmm \label{2.48} \\
& = & g(\UU[a],[b])-g([a],\UU[b]) \ .\nonumber
\end{eqnarray}

(c) is obvious.

(d) Because of $\oooo \gamma (z) =\eezz$ and ${\oooo\gamma}^*(\frac{\ddd z}{z})
= - \frac{\ddd z}{z}$ the flatness $\nnn (\trrr )=0$ means
$\nnn_\xi \circ \trrr = \trrr \circ \nnn_{\oooo \xi}$ and
$\nnn_\zdz \circ \trrr = \trrr \circ \nnn_{-\zdz}$.
The rest is obvious.

(e) Use lemma \ref{t2.7}.
\hfill $\qed$

\begin{definition}\label{t2.15}
A $(trTERP(w))$-structure is a $(TERP(w))$-structure such that the
$(T\www T)$-structure defined in lemma \ref{2.14} (e) is a 
$(trT\www T)$-structure.
\end{definition}

A $(trTERP(w))$-structure gives rise to a harmonic bundle with additional 
structure. This structure is discussed in the  next section, the 
correspondence is presented in section \ref{s2.7}.

\subsection{CV-structures}\label{s2.6}

The structures in \cite{CV1} motivated the following definition.

\begin{definition}\label{t2.16}
A CV-structure is a tuple $(K\to M,D,C,\www C,\kappa,h,\UU,\QQ)$
such that $(K\to M,D,C,\www C)$ is a $(DC\www C)$-structure 
(definition \ref{t2.9}) and the
other objects have the following properties.

(a) $\kappa$ is a fiberwise $\C$-antilinear automorphism of $K$ as
$\ciii$-bundle with
\begin{eqnarray} \label{2.49}
\kappa^2=\id\ ,\\
D(\kappa)= 0\ ,\label{2.50}\\
\kappa C \kappa =\www C\ .\label{2.51}
\end{eqnarray} 

(b) $h$ is a hermitian pseudo-metric on $K$; that means, it is 
linear on the left, semilinear on the right, nondegenerate and 
satisfies $h(b,a)=\oooo{h(a,b)}$.
It also has the three properties:\\
it takes real values on the real subbundle 
$K_\R:= \ker (\kappa -\id)\subset K$;
\begin{eqnarray} \label{2.52}
D(h) &=& 0\ ;\\
h(C_Xa,b) &=& h(a,\www C_{\oooo X} b) \label{2.53}
\end{eqnarray} 
for $a,b\in \ciik,\ X\in \tmen$.

(c) $\UU$ and $\QQ$ are $\ciii_M$-linear endomorphisms of $K$ with 
\begin{eqnarray} 
[C,\UU]=0\ ,\label{2.54}\\
D'(\UU)-[C,\QQ]+C=0\ ,\label{2.55}\\
D''(\UU)=0\ ,\label{2.56}\\
D'(\QQ)+[C,\kappa \UU \kappa]=0\ ,\label{2.57}\\
\QQ +\kappa \QQ \kappa =0\ ,\label{2.58}\\
h(\UU a,b)=h(a,\kappa\UU\kappa b)\ ,\label{2.59}\\
h(\QQ a,b)=h(a,\QQ b)\ .\label{2.60}
\end{eqnarray} 
\end{definition}

\begin{remarks}\label{t2.17}
(i) If the hermitian pseudo-metric $h$ is positive definite, we call the 
CV-structure a CV$\oplus$-structure.

(ii) A CV-structure on $K\to M$ is given by the data 
$(DC\www C\kappa h\UU\QQ)$. One can easily define weaker structures,
for example with the data
$(DC\www C\kappa)$, $(DC\www Ch)$, $(DC\www C\kappa h)$, 
$(DC\www C\kappa\UU\QQ)$ and with those 
properties which are formulated in these
data. Then Simpson's notion of a harmonic bundle \cite{Si4} is a 
$(DC\www C h)$-structure with positive definite $h$.

(iii) $\kappa$ yields a real structure on $K$, a real subbundle 
$K_\R := \ker (\kappa -\id)$ with
\begin{eqnarray} \label{2.61}
K = K_\R \oplus iK_\R \mbox{ \ and \ } 
iK_\R = \ker (\kappa +\id)\ .
\end{eqnarray} 

(iv) $\kappa$ and $h$ induce a $\ciii_M$-bilinear nondegenerate pairing
$g:=h(\cdot ,\kappa \cdot)$ on $K$. It satisfies $D(g)=0$.
It is symmetric because $h$ takes real values on $K_\R$.
These two properties are equivalent.

(v) The pairing $g$ satisfies
\begin{eqnarray} 
g(C_Xa,b) &=& g(a,C_Xb) \mbox{ \ \ for }X\in \tmen\ ,\label{2.62}\\
g(\www C_{\oooo X}a,b) &=& g(a,\www C_{\oooo X}b) 
\mbox{ \ \ for }X\in \tmen\ ,\label{2.63}\\
g(\UU a,b) &=& g(a,\UU b) \ ,\label{2.64}\\
g(\kappa \UU\kappa a,b) &=& g(a,\kappa\UU\kappa b) \ ,\label{2.65}\\
g(\QQ a,b) &=& - g(a,\QQ b) \ .\label{2.66}
\end{eqnarray} 
Two of the three properties \eqref{2.51}, \eqref{2.53}, \eqref{2.62} 
imply the third. Two of the three properties 
\eqref{2.58}, \eqref{2.60}, \eqref{2.66} 
imply the third.

(vi) With \eqref{2.51} one recovers $\www C$ from $C$ and $\kappa$,
with \eqref{2.53} one recovers $\www C$ from $C$ and $h$.
With \eqref{2.52} one recovers $D$ as the metric connection from
$h$ and from the holomorphic structure $\OO(K)$ defined by $D''$.

(vii) Let $\cooo_M$ be the sheaf of real analytic functions on $M$
and define the real analytic structure on $K$ with sheaf of sections
$\cook := \OO(K)\otimes_{\OO_M} \cooo_M$. It turns out that the whole
CV-structure is automatically compatible with this real analytic 
structure, see lemma \ref{t2.18} (a). 
But regarding the weaker structures in (ii), this property is only clear
for $(DC\www C\kappa\UU\QQ)$.
Therefore definition \ref{2.16} was not given
rightaway in the real analytic category.

(viii) The endomorphism $\QQ$ turns up in \cite{CFIV} as a
`new supersymmetric index'. In the singularity case its real analytically
varying eigenvalues are generalizations of the (shifted) spectral numbers.
This will be discussed in chapter \ref{s8}.
In general a CV-structure is a grand generalization of a variation
of Hodge structures, see chapter \ref{s3}, and the eigenspaces of $\QQ$ 
are generalizations of the subspaces of the Hodge decomposition.

(ix) Often there exists a global vector field $E\in \tm$ with 
$\UU= -C_E$. Then \eqref{2.51}, \eqref{2.57} and \eqref{2.58} show
\begin{eqnarray}
D'_E(\QQ)= D''_{\oooo E}(\QQ) = [\UU,\kappa\UU\kappa] \mbox{ \ and \ }
D_{E-\oooo E}(\QQ)=0\ .\label{2.66b}
\end{eqnarray}
When a CV-structure comes together with a
Frobenius manifold (see section \ref{s5.4}) the Euler field has this
property. Then the behaviour of $\QQ$ along $E+\oooo E$ is very interesting.
\end{remarks}

\begin{lemma}\label{t2.18}
Consider a CV-structure.

(a) The tensors $\www C,\kappa, h, \QQ$ and the connection $D$ 
respect the real analytic structure $\cook$. The other tensors 
$C,g,\UU$ respect the holomorphic structure $\OO(K)$.
Therefore a CV-structure is a real analytic structure.

(b) Suppose that $h$ is positive definite.
Then $\QQ$ is a hermitian endomorphism, it is semisimple with (real
analytically varying) real eigenvalues which are distributed
symmetrically around 0. Define
$K_{t,\lambda}:=\ker (\QQ -\lambda\id:K_t\to K_t)$ for $\lambda\in \R$.
Then
\begin{eqnarray} 
&& \kappa :K_{t,\lambda}\to K_{t,-\lambda} \ , \label{2.67}\\
&& h(K_{t,\lambda},K_{t,\lambda'})=0 \mbox{ \ \ for }\lambda\neq \lambda' \ ,
\label{2.68}\\
&& g(K_{t,\lambda},K_{t,\lambda'})=0 \mbox{ \ \ for }\lambda+\lambda'\neq 0\ .
\label{2.69}
\end{eqnarray} 
\end{lemma}

{\it Proof.}
(a) The statement on $C,g,\UU$ follows from 
$D''(C)=0$, $D''(\UU)=0$, $D(g)=0$. 
The statement on $\www C, \kappa,  h , \QQ$ and $D$ will follow from the proof
of theorem \ref{t2.19}, see remark \ref{t2.20} (iii).

(b) This follows from \eqref{2.58}, \eqref{2.60}, \eqref{2.66}.
\hfill $\qed$

\subsection{A correspondence}\label{s2.7}

\begin{theorem}\label{t2.19}

Fix $w\in\Z$. There is a one-to-one correspondence between 
$(trTERP(w))$-structures and CV-structures. It extends that in
lemma \ref{t2.11}. It is given by the steps in (a) and (b).
They are inverse to one another.

(a) Let $(K\to M,D,C,\www C,\kappa,h,g,\UU,\QQ)$ be a CV-structure.
Let $\pi:\P^1\times M\to M$ be the projection. 
Define $\hat H:=\pi^*K$ and $H:=\hat H|_\cmmm$. 
Extend $C,\www C,\kappa,\UU,\QQ$ canonically to $\hat H$.
Extend $D$ to a connection on $\hat H$ such that $D_\dz$
and $D_{\oooo \dz}$ vanish on sections in $\pi^{-1}\ciik$.

$(\alpha)$ Define
\begin{eqnarray} \label{2.70}
\nnn := D+\eezz C + z\www C + (\eezz \UU -\QQ +\frac{w}{2}\id
- z\kappa \UU\kappa)
\dzz\ .
\end{eqnarray} 
Then $\nnn$ is a flat connection on $H|_\csmmm$, the $(0,1)$-part
$\nnn^{(0,1)}=D^{(0,1)}+z\www C$ defines a holomorphic structure
$\OO_\nnn(H)$ on $H$, the connection $\nnn$ 
is compatible with it and has a pole of Poincar\'e rank 1 along
$\nmmm$.

$(\beta)$ Define $\gamma :\P^1\to\P^1$ and 
$\gamma_\nnn:H_{(z,t)}\to H_{(\gamma(z),t)}$ for $z\in \C^*$ as in 
lemma \ref{t2.14} (c). 
Let $\psi_z:\hat H_{(z,t)}\to K_t$ for $z\in \P^1$ be the canonical
projection. Define a $\C$-antilinear map
\begin{eqnarray} \label{1.71a}
\tau:\hat H_{(z,t)}\to \hat H_{(\gamma(z),t)}\ , \ 
\tau(a):=\psi_{\gamma(z)}^{-1}\circ \kappa \circ \psi_z(a)
\mbox{ \ for }z\in \P^1\ .
\end{eqnarray} 
Then $\tau^2=\id$ and $\tau(za)=\eezz \tau(a)$ 
(this has to be understood with $z$ as function on $\C^*$, not pointwise).
For $(z,t)\in \csmmm$ define the maps
\begin{eqnarray} \label{2.71b}
&& \trrr :H_{(z,t)}\to H_{(\gamma(z),t)}\ ,\ 
\trrr (a) = \tau(z^w a)\ ,\\
&& \kappa_H:=\gamma_\nnn \circ \trrr : H_{(z,t)}\to H_{(z,t)}\ .\label{2.72}
\end{eqnarray} 
Then $\trrr$ and $\kappa_H$ are $\nnn$-flat, fiberwise $\C$-antilinear
and satisfy $\trrr^2=\id$ and $\kappa_H^2=\id$.
Thus $\kappa_H$ defines a $\nnn$-flat real subbundle 
$H_\R:=\ker (\kappa_H-\id)\subset H|_{\csmmm}$.

$(\gamma)$ Define a pairing
\begin{eqnarray} \label{2.73a}
P:H_{(z,t)}\times H_{(-z,t)}&\to& \C \mbox{ \ \ for }(z,t)\in\csmmm\\
(a,b) &\mapsto& z^w\cdot g(\psi_za,\psi_{-z}b)\ .
\end{eqnarray} 
It satisfies all properties in definition \ref{2.12} (d).
The tuple $(H,\nnn,H_\R,P)$ is a $(trTERP(w))$-structure.

(b) Let $(H,\nnn,H_\R,P)$ be a $(trTERP(w))$-structure.
Let $K:=H|_\nmmm$ and $D'',C,D',\www C$ be the $(DC\www C)$-structure
in lemma \ref{t2.11} (b).
Define $g,\ \UU,\ \tau$ and $\hat H$ as in lemma \ref{2.14}.
Because $\hat H|_{\P^1\times \{t\}}$ is a trivial bundle for any
$t\in M$, there is a canonical projection $\psi:\hat H\to K$.
It restricts to canonical isomorphisms $\psi_z:H_{(z,t)}\to K_t$
for $(z,t)\in \pmmm$.
Then the map
\begin{eqnarray} \label{2.73b}
\kappa:K_t\to K_t\ , \ \kappa(a):=\psi_{\gamma(z)}\circ \tau \circ
\psi_z^{-1}(a) \mbox{ for some }z\in \C^*
\end{eqnarray} 
is independent of the choice of $z$ and is a $\C$-antilinear
involution. Define $h:=g(\cdot,\kappa \cdot)$ on $K$. 
Finally, lift 
$D, C,\www C,\kappa,\UU$ canonically to $\hat H$, using $\psi$.
Then there exists a unique endomorphism $\QQ$ of $K$ such that 
\eqref{2.70} holds. 
The tuple $(K\to M,D,C,\www C,\kappa,h,g,\UU,\QQ)$ is a CV-structure.
\end{theorem}

{\it Proof.}
(a) $(\alpha)$ $\nnn^2$ maps sections in $\pi_*\chiii(\hat H)$ to 
sections in 
\begin{eqnarray} \label{2.74}
&& (\frac{1}{z^2}\amzn + \eezz\amzn + \eezz\amee + \amz + z\amee
+ z\amnz + z^2\amnz  \\
&& + (\frac{1}{z^2}\amen + \eezz\ame + \ame + z\ame + z^2\amne)
\land \frac{\ddd z}{z})\otimes \pi_*\chiii(\hat H)\ .
\nonumber
\end{eqnarray} 
$\nnn^2=0$ holds because it splits into \eqref{2.22}--\eqref{2.24},
\eqref{2.54}--\eqref{2.57}, and the four identities
\begin{eqnarray} 
[\www C,\kappa\UU\kappa]=0\ ,\label{2.75}\\
D''(\kappa\UU\kappa)+[\www C,\QQ]+\www C = 0\ ,\label{2.76}\\
D'(\kappa\UU\kappa) = 0\ ,\label{2.77}\\
D''(\QQ)-[\www C,\UU]=0\ ,\label{2.78}
\end{eqnarray} 
which follow from \eqref{2.54}--\eqref{2.58} and \eqref{2.49}--\eqref{2.51}.
For the summand $C$ in \eqref{2.55} remark that $D$ in \eqref{2.70}
contains a covariant derivative $D_{\dz}$.

$(\beta)$ 
The only nontrivial claim is the flatness $\nnn(\trrr)=0$.
It is equivalent to 
\begin{eqnarray} \label{2.79}
\nnn_\xi\circ \tau = \tau \circ \nnn_{\oooo\xi} \mbox{ and }
\nnn_\zdz \circ \tau = \tau \circ (\nnn_{-\zdz} + w\cdot \id)\ .
\end{eqnarray} 
Both $\tau $ and $\kappa$ act on the sheaf $\pi_*\chiii(\hat H)$ of fiberwise
global and holomorphic $\ciii$-sections, and these actions coincide.
Therefore for $a\in \pi_*\chiii(\hat H)$
\begin{eqnarray} \label{2.80}
&& (\nnn_\xi\circ\tau - \tau\circ\nnn_{\oooo \xi})(a) \\
&=& (D_\xi+\eezz C_\xi+z\www C_\xi)\tau(a) - 
    \tau(D_{\oooo\xi}+\eezz C_{\oooo\xi}+z\www C_{\oooo\xi})(a)\nonumber\\
&=& (D_\xi\circ \kappa - \kappa\circ D_{\oooo\xi})(a) 
+\eezz(C_\xi\circ\kappa - \kappa\circ \www C_{\oooo \xi})(a) 
+z(\www C_\xi\circ \kappa - \kappa\circ C_{\oooo \xi})(a)\nonumber\\
&=&0 \nonumber
\end{eqnarray} 
and
\begin{eqnarray} \label{2.81}
&& (\nnn_\zdz\circ\tau - \tau\circ\nnn_{-\zdz} -w\id)(a)\\
&=& (\eezz\UU - \QQ +\frac{w}{2}\id -z\kappa\UU\kappa)\tau(a) 
+ \tau(\eezz\UU-\QQ+\frac{w}{2}\id-z\kappa\UU\kappa)(a) -w\cdot a\nonumber\\
&=& -(\QQ+\kappa \QQ \kappa )\kappa (a)  = 0\ . \nonumber
\end{eqnarray} 

$(\gamma)$ 
The pairing $P$ is $(-1)^w$-symmetric because $g$ is symmetric.
The only nontrivial claims are $\nnn(P)=0$ and that $P$ takes
values in $i^w\R$ on $H_\R$. Consider $a,b\in \ciik$ and 
their canonical lifts $\www a,\www b \in \pi_*\chiii(\hat H)$ to $\hat H$.
The following calculation shows $\nnn (P)=0$. 
It uses the $z$-sesquilinearity of $P$ and in the last step 
$D(g)=0$ and \eqref{2.62}--\eqref{2.66}.
\begin{eqnarray} \label{2.82}
&& \nnn (P)(\www a,\www b) \\
&=& P(\nnn \www a, \www b) + P(\www a,\nnn \www b)
-\ddd (P(\www a, \www b)) \nonumber\\
&=& z^w g((D + \eezz C + z\www C)a,b) + z^w g(a,(D-\eezz C -z\www C)b)
\nonumber \\
&& +  z^w g((\eezz \UU-\QQ+\frac{w}{2}\id -z\kappa \UU\kappa)a,b)\dzz 
\nonumber \\
&& + z^w g(a,(-\eezz \UU-\QQ+\frac{w}{2}\id +z\kappa \UU\kappa)b)\dzz 
\nonumber \\
&& - z^w \ddd (g(a,b)) - wz^w g(a,b)\dzz \nonumber \\
&=& 0\ .  \nonumber
\end{eqnarray} 

Because $H_\R$ and $P$ are $\nnn$-flat it is now sufficient to show
$P((H_\R)_{(1,t)}\times (H_\R)_{(-1,t)})\subset i^w\R$.
But for $a,b\in K_t$ and $\psi_{\pm1}^{-1}(a)\in H_{(\pm1,t)}$
\begin{eqnarray} \label{2.83}
 \kappa_H(\psi_{\pm1}^{-1}(a)) =\trrr (\psi_{\pm1}^{-1}(a)) 
= (\pm1)^w\tau(\psi_{\pm1}^{-1}(a)) = (\pm1)^w \psi_{\pm1}^{-1}\kappa(a)\ .
\end{eqnarray} 
Therefore
\begin{eqnarray} \label{2.84a}
(H_\R)_{(\pm1,t)} = \psi_{\pm1}^{-1}(\ker (\kappa-(\pm1)^w\id):K_t\to K_t) \ ,
\end{eqnarray} 
and with \eqref{2.61}
\begin{eqnarray}\label{2.84b}
(H_\R)_{(1,t)} = \psi_1^{-1}(K_{\R,t}) \mbox{ and } 
(H_\R)_{(-1,t)} =\psi_{-1}^{-1}(i^wK_{\R,t}) \ ,
\end{eqnarray}
and for $a\in K_{\R,t}$, $b\in i^w K_{\R,t}$
\begin{eqnarray} \label{2.85}
P(\psi_1^{-1}(a),\psi_{-1}^{-1}(b)) = g(a,b) = (-1)^w h(a,b) \in i^w \R\ .
\end{eqnarray} 

(b) The map $\tau$ from lemma \ref{t2.14} acts on $\hat H$ by definition
of $\hat H$ and thus also on $\pi_*\chiii(\hat H)\cong \ciik$. Therefore
$\kappa$ in \eqref{2.73b} is independent of the choice of $z$.

The lifts with the projection $\psi:\hat H\to K$ of $D,C,\www C,\kappa,\UU$ to
$\hat H$ act on $\pi_*\chiii(\hat H)\cong \ciik$, and by definition
of $\kappa$
\begin{eqnarray} \label{2.86}
\kappa=\tau\mbox{ \ \ on }\pi_*\chiii(\hat H)\ .
\end{eqnarray} 
Define on $\chiii(H|_\csmmm)$ the map
\begin{eqnarray} \label{2.87}
\QQ := -\nnn_\zdz +D_\zdz + \eezz \UU+\frac{w}{2}\id - z\kappa\UU\kappa\ .
\end{eqnarray} 
Then $\tau(zb)=\eezz \tau(b)$, \eqref{2.44} and \eqref{2.86} show
for $a\in \pi_*\chiii(\hat H)$
\begin{eqnarray} \label{2.88}
\QQ a = -\tau \QQ \tau a\ .
\end{eqnarray} 
$\QQ a$ extends as a $\chiii$-section to $\nmmm$ because of the definition
of $\UU$ and to $\immm$ because of \eqref{2.88}.
Therefore $\QQ a\in \pi_*\chiii(\hat H)$, and $\QQ$ descends to a 
$\ciii$-linear endomorphism of $K$.
Together with the proof of lemma \ref{t2.11} this shows \eqref{2.70}.
Now \eqref{2.88} implies $\QQ=-\kappa\QQ\kappa$.
The calculation \eqref{2.80} and \eqref{2.44} show
$\www C=\kappa C\kappa$ and $D(\kappa)=0$. By definition $\kappa^2=\id$.

Splitting $\nnn^2=0$ into pieces with \eqref{2.70} and \eqref{2.74}
gives \eqref{2.22}--\eqref{2.24} and \eqref{2.54}--\eqref{2.57}.

Next we have to show
\begin{eqnarray}\label{2.89}
P(a,b)= z^wg(\psi_z a,\psi_{-z}b)\ .
\end{eqnarray}
Consider a neighborhood $U_1\times U_2 
\subset \pmmm$ of a point in $\nmmm$ such that
$U_1$ is invariant with respect to $z\mapsto -z$ and consider sections
$a,b\in \OO(\hat H)(U_1\times U_2)$. Then
$P(a,b)\in z^w\OO_{U_1\times U_2}$ by \eqref{2.30b}, 
and the definitions of $\tau_{real}$ and $\tau$ show
\begin{eqnarray}
P(\tau_{real}(a),\tau_{real}(b))\in 
{\oooo \gamma}(z)^w \OO_{\gamma(U_1)}\ciii_M\ ,\label{2.90}\\
P(\tau(a),\tau(b))\in 
z^w \OO_{\gamma(U_1)}\ciii_M\ .\label{2.91}
\end{eqnarray}
Therefore
\begin{eqnarray}\label{2.92}
P(a,b)\in z^w\cdot \ciii_M \mbox{ \ \ for \ }a,b\in \pi_*\chiii (\hat H)
\end{eqnarray}
and \eqref{2.89} holds. 

Now one can repeat the calculation in 
\eqref{2.82}. Together with $\nnn(P)=0$ it shows 
$D(g)=0$ and \eqref{2.62}--\eqref{2.66}. By definition of $h$ these imply
\eqref{2.59} and \eqref{2.60}. The pairing $P$ takes values in $i^w\R$ 
on $H_\R$. Define $K_\R:= \ker (\kappa -\id)$. 
Now \eqref{2.86} shows \eqref{2.83}--\eqref{2.84b}.
Reading \eqref{2.85} backwards shows that $h$ takes real values on $K_\R$.
This and the symmetry of $g$ imply that $h$ is hermitian.
$D(h)=0$ follows from $D(\kappa)=0$ and $D(g)=0$.
\hfill $\qed$

\begin{remarks}\label{t2.20}
(i) One can formulate and show many intermediate correspondences
for structures stronger than those in lemma \ref{t2.11} and
weaker than those in 
theorem \ref{t2.19}. For example in \cite[Lemma 3.1]{Si5}
and \cite[Lemma 2.2.2]{Sab6} a correspondence is given between
$(DC\www Ch)$-structures and $(trT\www T)$-structures
with a sesquilinear pairing, which turns up as $P(\cdot ,\tau \cdot )$
in the case of a CV-structure.

(ii) Harmonic bundles \cite{Si4} are $(DC\www Ch)$-structures with a
positive definite hermitian metric $h$. 
In Simpson's nonabelian Hodge theorem \cite{Si3}
real structures are not considered. But in \cite{Si5} $(T)$-structures
with real structures with properties as $\trrr$ for $w=0$ are discussed
under the name (variation of) `circular real structures'. 
But the $(trTERP(w))$-structures and CV-structures here are more 
motivated by the work of Cecotti and Vafa and by singularity theory
than by Simpson's work.

(iii) The proof of lemma \ref{2.18} (a):
The bundle $H$ in lemma \ref{t2.14} carries a $\chooo$-structure
with sheaf $\chooo(H)$ of real analytic sections which are holomorphic
in $z$-direction. The map $\tau$ in lemma \ref{t2.14} respects this
structure and the induced structure on $\hat H$. Therefore $\kappa$
in theorem \ref{t2.19} (b) respects the real analytic structure
of $K$. Because all other nonholomorphic operators in the CV-structure
are obtained from $\kappa$ and holomorphic operators, the whole
CV-structure is real analytic.
This shows lemma \ref{t2.18} (a).
\end{remarks}

\clearpage

\section{Variation of Hodge structures}\label{s3}
\setcounter{equation}{0}

\noindent
Simpson \cite{Si1}--\cite{Si4} defined harmonic bundles as a generalization
of variations of polarized Hodge structures. But from the harmonic
bundle of a variation of Hodge structures one cannot recover the 
Hodge filtrations without additional information. The tensors
$\UU$ and $\QQ$ of a CV-structure provide such information.
They allow to formulate the correspondence in theorem \ref{t3.1}.
Of course, most of it is due to Simpson. Also in \cite{CV1}\cite{CV2}
it is stated that the structures there generalize variations
of Hodge structures.

\begin{theorem}\label{t3.1}
Fiz an integer $w$. There is a natural correspondence between the 
following structures.

$(\alpha)$ CV$\oplus$-structures with $\UU=0$.

$(\beta)$ Sums of two variations of polarized Hodge structures;
one is of weight $w$, the other is of weight $w-1$ and is equipped
with an automorphism with eigenvalues $\neq 1$.
\end{theorem}

In the first section the notion of a variation of Hodge structures
is recalled. Subsequently the correspondence is made
precise and is proved. The last section contains some remarks
about special K\"ahler manifolds.

\subsection{Variation of polarized Hodge structures}\label{s3.1}

A {\it variation of Hodge structures of weight $w$} is a holomorphic
vector bundle $H\to M$ with the additional structure:

$(\alpha)$ a flat holomophic connection 
$\nnn:\OO(H)\to \oomm^1_M\otimes \OO(H)$,

$(\beta)$ a $\nnn$-flat real subbundle $H_\R\to M$ with 
$H_t=H_{\R,t}\oplus iH_{\R,t}$ for $t\in M$; the corresponding
complex conjugation is denoted by $\oooo\cdot$ or by $\kappa$.

$(\gamma)$ an exhaustive decreasing filtration $F^\bullet$ by holomorphic
subbundles $F^p\subset H$ such that
\begin{eqnarray} \label{3.1}
H_t = F^p_t\oplus \oooo{F^{w+1-p}_t} \mbox{ \ \ for any }t\in M
\end{eqnarray} 
and (Griffiths transversality)
\begin{eqnarray} \label{3.2}
\nnn :\OO(F^p)\to \oomm^1_M\otimes \OO(F^{p-1})\ .
\end{eqnarray} 

The Hodge subbundles are the complex $\ciii$-subbundles $H^{p,w-p}$
(in fact, real analytic) with
\begin{eqnarray} \label{3.3}
H^{p,w-p}_t = F^p_t\cap \oooo{F^{w-p}_t} \mbox{ \ \ for }
t\in M\ .
\end{eqnarray} 
They satisfy $\oooo{H^{p,w-p}} = H^{w-p,p}$. Property \eqref{3.1}
is equivalent to the Hodge decomposition
\begin{eqnarray} \label{3.4}
F^p_t = \oplus_{q\geq p}H_t^{q,w-q} \mbox{ \ \ for }t\in M\ .
\end{eqnarray} 

A {\it polarization} of such a variation of Hodge structures is a 
$\nnn$-flat $\OO_M$-bilinear pairing $S$ on $H$ which takes real values
on $H_\R$, is nondegenerate and $(-1)^w$-symmetric and satisfies
\begin{eqnarray} \label{3.5}
S(F^p_t,F^q_t)=0 \mbox{ \ \ for }p+q>w\ ,\\
i^{p-(w-p)}S(v,\oooo v)>0 \mbox{ \ \ for }v\in H^{p,w-p}_t-\{0\}\ .
\label{3.6}
\end{eqnarray} 
\eqref{3.5} is equivalent to 
\begin{eqnarray} \label{3.7}
S(H^{p,w-p}_t,H^{q,w-q}_t)=0 \mbox{ \ \ for }
p+q\neq w\ .
\end{eqnarray} 
\eqref{3.6} motivates the definition of a sesquilinear form
\begin{eqnarray}  
h&:&H_t\times H_t \to \C \nonumber\\
h&:&H^{p,w-p}_t\times H^{q,w-q}_t \to 0 \mbox{ \ \ for } p\neq q\ ,
\label{3.8}\\
h(a,b) &:=& \frac{(-1)^p}{(2\pi i)^w} S(a,\oooo b) \mbox{ \ \ for }a,b\in 
H^{p,w-p}_t\ . \label{3.9} 
\end{eqnarray} 
The factor $1/(2\pi )^w$ is inserted in order to make $h$ and $S$
consistent with the singularity case, see section \ref{s8.1}.
The form $h$ is hermitian, positive definite, and takes real values
on $H_\R$.

An automorphism of a variation of polarized Hodge structures is a 
$\nnn$-flat bundle automorphism which respects $H_\R, F^\bullet,S$,
and $h$. Because $h$ is positive definite, the automorphism is semisimple
with eigenvalues on the unit circle.

\subsection{Correspondence with special CV-structures}\label{s3.2}

\begin{lemma}\label{t3.2}
Let $(H\to M,\nnn,H_\R,S,A)$ be a variation of polarized Hodge
structures of weight $w$ with an automorphism $A$ with eigenvalues
$\neq(-1)^{w+1}$.

Define $\kappa = \oooo\cdot$ and $h$ as above.
There exist unique maps $D',D'',C,\www C:\ciii(H)\to \ame\otimes\ciii(H)$
with
\begin{eqnarray} 
D':\ciii(H^{p,w-p})&\to& \amen\otimes\ciii(H^{p,w-p})\ ,\label{3.11}\\
D'':\ciii(H^{p,w-p})&\to& \amne\otimes\ciii(H^{p,w-p})\ ,\label{3.12}\\
C:\ciii(H^{p,w-p})&\to& \amen\otimes\ciii(H^{p-1,w+1-p})\ ,\label{3.13}\\
\www C:\ciii(H^{p,w-p})&\to& \amne\otimes\ciii(H^{p+1,w-1-p})\label{3.14}\ ,
\end{eqnarray} 
and
\begin{eqnarray} \label{3.15}
\nnn  = D'+D''+C+\www C\ .
\end{eqnarray} 
Define $\UU:=0$ and define the endomorphism $\QQ:H_t\to H_t$ by
\begin{eqnarray} 
\QQ| \ (H_t^{p,w-p}\cap\ker (A-\lambda\id)):=\alpha\id \label{3.16}\\
\mbox{for }\alpha\in \R \mbox{ with }e^{2\pi i\alpha}=\lambda,\ 
\left[\alpha+\frac{w+1}{2}\right] =p\ .\label{3.17}
\end{eqnarray}
Then $e^{2\pi i \QQ}=A$, and $(H\to M, D',D'',C,\www C,\kappa,h,\UU,\QQ)$ is
a CV$\oplus$-structure.
\end{lemma}

{\it Proof.}
Griffiths transversality, holomorphicity of the subbundles $F^p$
and $\nnn$-flatness of $H_\R$ show
\begin{eqnarray} 
\nnn :\ciii(F^p) \to \amen\otimes \ciii(F^{p-1} )
   + \amne\otimes \ciii(F^p) &,&\label{3.18}\\
\nnn :\ciii(\oooo{F^p}) \to \amne\otimes \ciii(\oooo{F^{p-1}}) 
   + \amen\otimes \ciii(\oooo{F^p}) &,&\label{3.19}\\
\nnn :\ciii(H^{p,w-p}) \to \amen\otimes \ciii(H^{p,w-p}) 
+\amne\otimes \ciii(H^{p,w-p}) && \nonumber \\
+ \amen\otimes \ciii(H^{p-1,w+1-p}) 
  +\amne\otimes \ciii(H^{p+1,w-1-p})&.& \label{3.20}
\end{eqnarray} 
The maps $D',D'',C,\www C$ with \eqref{3.11}--\eqref{3.15} are well
defined. $D:=D'+D''$ is a connection, $C$ and $\www C$
are $\ciii_M$-linear. Splitting $\nnn^2$ into pieces according to
the pieces ${\mathcal A}_M^{r,2-r}\otimes \ciii(H^{p,w-p})$ 
of $\amz\otimes \ciii(H)$ shows that
$(H\to M,D,C,\www C)$ is a $(DC\www C)$-structure (definition \ref{t2.9}).

The complex conjugation respects the Hodge decomposition. Therefore
the flatness $\nnn(\kappa)=0$ of $H_\R$ splits into $D(\kappa)=0$ and
$(C+\www C)(\kappa)=0$. The last equation means $\kappa C\kappa =\www C$.
Because of \eqref{3.7} the flatness $\nnn(S)=0$ of the pairing $S$
splits into $D(S)=0$ and
\begin{eqnarray} 
S(C_Xa,b)+S(a,C_Xb)&=&0\ ,\label{3.22}\\
S(\www C_{\oooo{X}}a,b)+S(a,\www C_{\oooo{X}}b) &=&0\ .\label{3.23}
\end{eqnarray} 
The definition of $h$, $D(S)=0$, $D(\kappa)=0$,
and \eqref{3.11}+\eqref{3.12} show $D(h)=0$ and
\begin{eqnarray}  \nonumber
h(C_Xa,\oooo b) &=& \frac{(-1)^{p-1}}{(2\pi i)^w}S(C_Xa,\oooo b)
= \frac{(-1)^p}{(2\pi i)^w}S(a,C_X\oooo b) \\
 &=& h(a,\kappa C_X\kappa b)=h(a,\www C_{\oooo X}b) \label{3.24}
\end{eqnarray} 
for $a\in \ciii(H^{p,w-p})$. It rests to show the equations for $\UU$ and
$\QQ$ in a CV-structure. The automorphism $A$ is $A=e^{2\pi i\QQ}$
and commutes with $C_X$, because $C_X$ comes from the variation of Hodge
structures. The definition of $\QQ$ shows $D(\QQ)=0$
and 
\begin{eqnarray} \label{3.25}
C_X:\ker (\QQ-\alpha\id)
\to \ker (\QQ-(\alpha-1)\id) \ .
\end{eqnarray} 
Therefore $[C,\QQ]-C=0$. The automorphism $A$ respects the real structure.
Therefore
\begin{eqnarray} \label{3.26}
\kappa(H^{p,w-p}\cap\ker(A-\lambda\id))=
H^{w-p,p}\cap \ker(A-\oooo\lambda\id) .
\end{eqnarray} 
This together with the identity for 
$\alpha\in \R-(\frac{w+1}{2}+\Z)$
\begin{eqnarray} \label{3.27}
\left[-\alpha+\frac{w+1}{2}\right]-(w-p) 
= -(\left[\alpha+\frac{w+1}{2}\right]-p)
\end{eqnarray} 
shows $\kappa \QQ \kappa =-\QQ$. Finally, $h$ is invariant under $A$, so
\begin{eqnarray} \label{3.28}
\qquad h:\left(H^{p,w-p}_t\cap\ker(A-\lambda_1\id)\right)\times 
\left(H^{q,w-q}_t\cap\ker(A-\lambda_2\id)\right)\to 0 \\
\mbox{ \ \ for }(p,\lambda_1)\neq (q,\lambda_2)\ , \nonumber
\end{eqnarray}
and $h(\QQ a,b)=h(a,\QQ b)$. 
\hfill $\qed$

\begin{remark}\label{t3.3}
If one starts just with a variation of polarized Hodge structures of weight
$w$, one can define an automorphism $A:=(-1)^w\id$ and one obtains a 
CV$\oplus$-structure. If one has already an automorphism, but with 
$(-1)^{w+1}$ as eigenvalue. one can change $A$ to the sum of $A$ on the
subbundle $\bigoplus_{\lambda\neq (-1)^{w+1}}\ker(A-\lambda\id)$ 
and to $-A$ on the subbundle $\ker(A-(-1)^{w+1}\id)$ and then proceed.
\end{remark}

The inverse of lemma \ref{t3.2} is straightforward. We formulate it.

\begin{lemma}\label{t3.4}
Let $(H\to M,D,C,\www C,\kappa,h,\UU,\QQ)$ be a CV$\oplus$-structure with 
$\UU=0$ and such that $\QQ$ has no eigenvalues in 
$\frac{w+1}{2}+\Z$. 

Define a connection $\nnn:=D+C+\www C$ and define
\begin{eqnarray} \label{3.29}
H^{p,w-p}_t&:=&\bigoplus_{\alpha:\ [\alpha+\frac{w+1}{2}]=p}
\ker(\QQ-\alpha\id:H_t\to H_t)\ ,\\
F^p_t&:=&\bigoplus_{q\geq p}H^{q,w-q}_t\ ,\label{3.30}\\
S&:&H_t\times H_t \to \C \mbox{ \ \ with}\label{3.31}\\
S(a,b)&:=& (2\pi i)^w(-1)^ph(a,\oooo b) \mbox{ \ \ for }
a\in H^{p,w-p}_t,\ b\in H_t\ ,\nonumber\\
A &:=& e^{2\pi i\QQ}\ .\label{3.32}
\end{eqnarray} 
Then $(H\to M,\nnn,H_\R,S,F^\bullet,A)$ is a variation of polarized
Hodge structures of weight $w$ with an automorphism $A$.
\end{lemma}

{\it Proof.}
$\nnn$ is flat because of lemma \ref{t2.11} (a). $\QQ$ is semisimple with
real eigenvalues by lemma \ref{t2.18} (b). It satisfies the equations
$D(\QQ)=0$, $[C,\QQ]-C=0$ and $[\www C,\QQ]+\www C=0$.
They show $[C,A]=0=[\www C,A]$, $\nnn(A)=0$, and \eqref{3.11}--\eqref{3.14}.
Therefore Griffiths transversality holds and the bundles $F^p$ are
holomorphic subbundles with respect to the holomorphic structure
on $H$ defined by $\nnn^{(0,1)}$. The equation $\kappa \QQ \kappa=-\QQ$
shows $\kappa A\kappa =A$ and $\kappa\ker(\QQ-\alpha\id)=
\ker(\QQ+\alpha\id)$. With \eqref{3.27} this gives 
$\oooo{H^{p,w-p}}=H^{w-p,p}$, \eqref{3.3} and \eqref{3.1}.

The pairing $S$ satisfies \eqref{3.5}, \eqref{3.6} and $D(S)=0$
by definition. It is $(-1)^w$-symmetric and takes real values on $H_\R$
because $h$ is hermitian and takes real values on $H_\R$.
The equation $h(C_Xa,b)=h(a,\www C_{\oooo X}b)$ shows \eqref{3.22} and
\eqref{3.23}. With $D(S)=0$ this implies $\nnn(S)=0$.
Finally, because of $h(\QQ a,b)=h(a,\QQ b)$ the pairings $h$ and $S$
are invariant under the automorphism $A$.
\hfill $\qed$

\bigskip
{\it Proof of theorem \ref{t3.1}.} 
From $(\alpha)$ to $(\beta)$: On the subbundle 
$\bigoplus_{\alpha\neq (-1)^w}\ker(\QQ-\alpha\id)$ one applies lemma
\ref{t3.4} with $w-1$ instead of $w$ and equips the resulting 
variation of polarized Hodge structures of weight $w-1$ with the
automorphism $(-1)^wA$. On the subbundle $\ker(\QQ-(-1)^w\id)$
one applies lemma \ref{t3.4} with $w$ and obtains a variation of 
polarized Hodge structures of weight $w$.

From $(\beta)$ to $(\alpha)$: One equips the variation of Hodge structures
of weight $w$ with the automorphism $A:=(-1)^w\id$, and the one with
weight $w-1$ with the automorphism $A:=(-1)^w\cdot \ (the \ given \ 
automorphism)$. Then one applies lemma \ref{t3.2}, with $w-1$ instead
of $w$ in the case of the variation of Hodge structures of weight
$w-1$. One adds the resulting CV$\oplus$-structures.

It is clear that these two procedures are inverse to one another.
\hfill $\qed$

\begin{remarks}\label{t3.5}

(i) A CV$\oplus$-structure obtained as in lemma \ref{t3.2} decomposes
into CV$\oplus$-structures on the subbundles
$\ker(A-\lambda\id)+\ker(A-\oooo\lambda\id)$. For $\lambda\neq \pm 1$
one can apply lemma \ref{t3.4} 
with $w-1$ instead of $w$ and one gets a variation
of Hodge structures of weight $w-1$ on the subbundle, having started with
one of weight $w$.
But for $\lambda =\pm 1$ this is not possible; the smallest shift of
the weight of Hodge structures on $\ker (A-(\pm 1)\id)$ is 2, using
the Tate twist.

(ii) The formulation and the proof of theorem \ref{3.1} fit
to the case of families $f_t,t\in M$, of quasihomogeneous singularities
$f_t:\C^w\to \C$. There one has a flat bundle $H\to M$ of middle
cohomologies of Milnor fibers with a semisimple monodromy $(-1)^wA$;
on the subbundle 
$\ker((-1)^wA-\id)\subset H$ one has a variation of Hodge structures 
of weight $w$ and on the subbundle
$\bigoplus_{\lambda\neq 1}\ker((-1)^w A-\lambda\id)\subset H$
a variation of Hodge structures of weight $w-1$.
See step 4 in section \ref{s8.1}.

(iii) The point of view of $(TERP)$-structures for variations of Hodge
structures will be discussed in section \ref{s7.4}.
\end{remarks}

\subsection{Special geometry}\label{s3.3}

In \cite{CV1} and \cite{Du2} it is mentioned that CV$\oplus$-structures
are generalizations of special geometry. One can consider
`special geometry' roughly as a synonym for variations of 
polarized Hodge structures. Then theorem \ref{t3.1} makes this remark
in \cite{CV1} and \cite{Du2} precise. But one can also give
`special geometry' a more precise and narrower meaning.

In \cite{Fr} and \cite{ACD} two versions are considered: 
affine special K\"ahler manifolds and projective special K\"ahler manifolds.
The notion of projective special K\"ahler manifolds is the richer one.
It axiomatizes the geometry on moduli spaces of Calabi--Yau threefolds,
which is induced from the corresponding variations of Hodge structures
of weight 3 on the middle cohomologies \cite{BrG}\cite{Co}\cite{Fr}\cite{BCOV}.
In fact, it is equivalent to such a variation of Hodge structures.
We refer to these references for a discussion of this.

But in the definition in \cite{Fr} and \cite{ACD} of an affine special K\"ahler
manifold such a description in terms of variations of Hodge structures
is not given. The only point in this section is proposition \ref{t3.7},
which gives such a description.

By definition in \cite{ACD}, an {\it affine special complex manifold}
$(M,J,\nnn)$ is a complex manifold $(M,J)$, where $J:T^\R M\to T^\R M $ with
$J^2=-\id$ gives the complex structure, together with a torsion free
flat connection $\nnn$ on $T^\R M$ with 
\begin{eqnarray} \label{3.33}
(\nnn_{\xi_1}J)(\xi_2) = (\nnn_{\xi_2}J)(\xi_1) \mbox{ \ \ for }
\xi_1,\xi_2\in \tm^\C\ .
\end{eqnarray} 
By definition in \cite{ACD} and \cite{Fr},
an {\it affine special K\"ahler manifold}
$(M,J,\nnn,\omega)$ is an affine special complex manifold
$(M,J,\nnn)$ with a $\nnn$-flat symplectic form $\omega$ on $T^\R M$
which is $J$-invariant and such that the pseudo K\"ahler metric
$h:=\omega(J\cdot,\oooo\cdot)$ is positive definite.

Each fiber of the complex tangent bundle $T^\C M=T^{1,0}M\oplus T^{0,1}M$
of a complex manifold carries a natural Hodge structure of weight
$1$ with
\begin{eqnarray} \label{3.34}
0=F^2_t\subset F^1_t=T^{1,0}_tM\subset F^0_t=T^\C_t M\ .
\end{eqnarray}

\begin{lemma}\label{t3.6}
Let $\nnn$ be a flat torsion free connection on the complex tangent bundle
$T^\C M$ of a complex manifold $M$.
Then $T^{1,0}M=F^1$ is a holomorphic subbundle of $T^\C M$ with the 
holomorphic structure defined by $\nnn$ if and only if \eqref{3.33}
holds.
\end{lemma}

{\it Proof.} The first condition is equivalent to 
\begin{eqnarray} \label{3.35}
\nnn_{\oooo Y}X=0 \mbox{ \ \ for all } X,Y\in \tm\ .
\end{eqnarray} 
The second condition is equivalent to 
\begin{eqnarray} \label{3.36}
(\nnn_X J)(\oooo Y)= (\nnn_{\oooo Y}J)(X) \mbox{ \ \ for all }
X,Y\in \tm\ .
\end{eqnarray} 
Consider $X,Y\in \tm$. Then $\Lie_X(J)=0$ and 
\begin{eqnarray} \nonumber
0 &=& [X,J(\oooo Y)]-J([X,\oooo Y]) \\ 
&=& \nnn_X J(\oooo Y)-\nnn_{J(\oooo Y)}X - J\nnn_X\oooo Y + 
J\nnn_{\oooo Y}X       \nonumber\\
&=& (\nnn_X J)(\oooo Y) - (\nnn_{\oooo Y}J)(X) 
  + \nnn_{\oooo Y}J(X)-\nnn_{J(\oooo Y)}X \nonumber \\
&=& \left((\nnn_X J)(\oooo Y) - (\nnn_{\oooo Y}J)(X) \right)
+2i(\nnn_{\oooo Y}X)\ .\label{3.37}
\end{eqnarray} 
This shows the lemma.
\hfill $\qed$

\begin{proposition}\label{t3.7}
Let $(M,J)$ be a complex manifold and $F^\bullet$ as in \eqref{3.34}.

(a) $(M,J,\nnn)$ is an affine special complex manifold if and only if
$\nnn$ is a torsion free connection on $T^\R M$ 
and $\nnn$ and $F^\bullet$ give a variation of Hodge structures on $T^\C M$.

(b) $(M,J,\nnn,\omega)$ is an affine special K\"ahler manifold if and 
only if $\nnn$ is a torsion free connection on $T^\R M$ 
and $\nnn$, $F^\bullet$ and
$S:=2\pi \omega$ give a variation of polarized Hodge structures of weight
$1$ on $T^\C M$.
\end{proposition}

{\it Proof.} (a) Lemma \ref{t3.6}

(b) Start with a polarizing form $S$; define $\omega:=\frac{1}{2\pi}S$
and $h:=\omega(J\cdot,\oooo\cdot)$. Then $h$ and $S$ are related as in
\eqref{3.9}, and $h$ is positive definite.
\hfill $\qed$

\bigskip
So indeed, affine special complex manifolds and affine special 
K\"ahler manifolds are very natural objects.
By Lemma \ref{t3.2} they give $(DC\www C\kappa\UU\QQ)$-structures
(defined in remark \ref{t2.17} (ii)) respectively
CV$\oplus$-structures on the complex tangent bundle.

\clearpage

\section{$tt^*$ geometry on the tangent bundle}\label{s4}
\setcounter{equation}{0}

\noindent
Suppose that  a CV-structure on a vector bundle $K\to M$ is given such
that $\rk K=\dim M$ and $\OO(K)$ is a free $\tm$-module of rank 1.
Let $C$ be the Higgs field on $K$. Choose a generator $\zeta\in \OO(K)$
of $\OO(K)$ as a $\tm$-module. Then one can shift the whole CV-structure
with the isomorphism $C_\bullet \zeta :\tm \to \OO(K)$ to the 
holomorphic tangent bundle $TM$. Most of the induced structure depends
on the choice of $\zeta$, but not all.

In section \ref{s4.1} some induced structure, which does not depend on
the choice of $\zeta$, is studied, a multiplication on the tangent bundle
(lemma \ref{t4.1}), a unit field $e$ and an Euler field $E$. 
They form an F-manifold (definition \ref{t4.2} and lemma \ref{t4.3}).

In section \ref{t4.2} it is discussed how the induced CV-structure on $TM$
behaves with respect to the flows of unit field $e$, Euler field $E$, 
and their complex conjugates $\oooo e$ and $\oooo E$.
Theorem \ref{t4.5} leads to some wishes about the CV-structures
on $TM$. In section \ref{s5.4} it will be shown how a careful choice
of $\zeta$ satisfies these wishes and gives a Frobenius manifold structure.

\subsection{Multiplication on the tangent bundle}\label{s4.1}

Let $M$ be a complex manifold of dimension $m\geq 1$. A multiplication
$\circ$ on the holomorphic tangent bundle $TM$ is an $\OO_M$-bilinear
commutative and associative map $\circ:\tm\times \tm \to\tm$.
A global holomorphic vector field $e$ is called a unit field if 
$e\circ=\id$. 

A multiplication provides $TM$ with the structure of a Higgs bundle
with the Higgs field $C:\tm\to \oomm_M^1\otimes \tm$ defined by
\begin{eqnarray}\label{4.1}
C_XY:=-X\circ Y\ .
\end{eqnarray}
The minus sign is chosen only for the sake of the discussion of the
singularity case in section \ref{s8.1} (cf. also \cite[0.13.d]{Sab4}).

A multiplication can arise from an abstract Higgs bundle in the 
following way.

\begin{lemma}\label{t4.1}
Let $(K\to M,C)$ be a Higgs bundle with $\rk K=\dim M$ such that
$\OO(K)$ is a locally free $\tm$-module of rank 1.
Then there is a unique multiplication $\circ$ on $M$ with
$C_XC_Y=-C_{X\circ Y}$ and there is a unique unit field $e$.
It satisfies $C_e=-\id$.
\end{lemma}

{\it Proof.}
Choose a (local) generator $\zeta\in \OO(K)$ of $\OO(K)$ as a 
$\tm$-module. The map $C_\bullet \zeta :\tm \to \OO(K)$ is an
isomorphism. Define $e\in \tm$ by $C_e\zeta =-\zeta$ and define
$X\circ Y$ for $X,Y\in \tm$ by $C_XC_Y\zeta =-C_{X\circ Y}\zeta$.
Then $(M,\circ,e)$ is a manifold with commutative and associative
multiplication $\circ$ and unit field $e$. Because of 
$C_eC_Z\zeta = C_ZC_e\zeta=-C_Z\zeta$ and \
$C_XC_YC_Z\zeta= ... = -C_{X\circ Y}C_Z\zeta$ the unit field $e$
and the multiplication are independent of the choice of the 
generator $\zeta$ and satisfy $C_e=-\id$ and $C_XC_Y=-C_{X\circ Y}$.
\hfill $\qed$

\bigskip
If a Higgs bundle as in lemma \ref{t4.1} is part of a richer structure,
a $(DC\www C)$-structure or even a CV-structure, one can choose a 
generator $\zeta$ of $\OO(K)$ as $\tm$-module and pull the richer
structure to the tangent bundle $TM$. Lemma \ref{t4.3} says that then
the multiplication and the Lie bracket for vector fields satisfy the
compatibility condition \eqref{4.2}.

\begin{definition}\label{t4.2}\cite{HM}\cite[I.5]{Man}

(a) A manifold $M$ with multiplication $\circ$ on the holomorphic tangent
bundle and unit field $e$ is 
an {\it F-manifold} if the multiplication satisfies
\begin{eqnarray}\label{4.2}
\Lie_{X\circ Y}(\circ)= X\circ\Lie_Y (\circ )+Y\circ\Lie_X (\circ )\ .
\end{eqnarray}

(b) Let $(M,\circ,e)$ be an F-manifold. A vector field $E\in \tm$ is an
{\it Euler field} if it satisfies
\begin{eqnarray}\label{4.3}
\Lie_E(\circ)=\circ \ .
\end{eqnarray}
\end{definition}

F-manifolds with Euler fields are studied in \cite[part 1]{He2}.
It seems that all interesting manifolds with multiplication on the tangent
bundle are F-manifolds.

\begin{lemma}\label{t4.3}
Let $(M,\circ,e)$ be a manifold with multiplication and unit field.
Define a Higgs field $C$ by \eqref{4.1}. Let $D'$ be a $(1,0)$-connection on 
$TM$ as $\ciii$-bundle with $D'(C)=0$.

(a) Then $(M,\circ,e)$ is an F-manifold.

(b) Let $E\in \tm$ be a holomorphic vector field and $\QQ$ a 
$\ciii$-endomorphism of $TM$ with 
\begin{eqnarray}\label{4.4}
D'(-C_E)-[C,\QQ]+C=0\ .
\end{eqnarray}
Then $E$ is an Euler field of the F-manifold, that is, 
$\Lie_E(\circ)=\circ$.
\end{lemma}

{\it Proof.} 
(a) Formula \eqref{2.6} shows that $D'(C)=0$ serves as a replacement
for torsion freeness of $D'$,
\begin{eqnarray}\label{4.5}
[X,Y]\circ Z=D'_X(Y\circ)(Z)-D'_Y(X\circ)(Z)\ .
\end{eqnarray}
Using this, one simply calculates for $X,Y,Z,W\in \tmen$ (or $\tm)$
\begin{eqnarray}\label{4.6}
&&(\Lie_{X\circ Y}(\circ)-X\circ\Lie_Y (\circ )-Y\circ\Lie_X )(Z,W)\\
&=& [X\circ Y,Z\circ W]-[X\circ Y,Z]\circ W 
                            - Z\circ [X\circ Y,W]\nonumber\\
&& - X\circ [Y,Z\circ W]+ X\circ[Y, Z]\circ W 
                             + X\circ Z\circ [Y,W]\nonumber\\
&& - [X,Z\circ W]\circ Y + [X,Z]\circ Y\circ W
                             + Z\circ[X,W]\circ Y\nonumber\\
&=& D_{X\circ Y}(Z\circ W\circ )(e) - D_{Z\circ W}(X\circ Y\circ )(e) 
      \nonumber\\
&& - D_{X\circ Y}(Z\circ )(W) + D_Z(X\circ Y\circ )(W)      \nonumber\\
&& -Z\circ D_{X\circ Y}(W\circ)(e) + Z\circ D_W(X\circ Y\circ )(e)    
   \nonumber\\
&-& X\circ D_Y (Z\circ W\circ )(e) + X\circ D_{Z\circ W}(Y\circ )(e) 
      \nonumber\\
&& + X\circ D_Y(Z\circ )(W) - X\circ D_Z (Y\circ )(W)      \nonumber\\
&& + X\circ Z\circ D_Y(W\circ )(e) - X\circ Z\circ D_W(Y\circ )(e)   
   \nonumber\\
&-& D_X(Z\circ W\circ )(Y) + D_{Z\circ W}(X\circ )(Y)     \nonumber\\
&& + D_X(Z\circ )(Y\circ W) -D_Z (X\circ )(Y\circ W)       \nonumber\\
&& +Z\circ D_X(W\circ )(Y) - Z\circ D_W(X\circ )(Y)      \nonumber\\
&=& 0\ .\nonumber
\end{eqnarray}

(b) One calculates with \eqref{2.6} in the second and third step and 
\eqref{4.4} in the fifth step
\begin{eqnarray}\label{4.7}
&& C_{\Lie_E(\circ)(X,Y)-X\circ Y} \\
&=& C_{[E,X\circ Y]} - C_{X\circ[E,Y]} - C_{[E,X]\circ Y} - C_{X\circ Y}
    \nonumber \\
&=& D'_E (C_{X\circ Y}) - D'_{X\circ Y}(C_E) + C_XC_{[E,Y]} 
    + C_{[E,X]} C_Y + C_XC_Y\nonumber \\
&=& - D'_E(C_XC_Y) -D'_{X\circ Y}(C_E) + C_X(D'_E(C_Y)-D'_Y(C_E))
    \nonumber\\
&& + (D'_E(C_X)-D'_X(C_E))C_Y + C_XC_Y\nonumber\\
&=& -D'_{X\circ Y}(C_E) - C_XD'_Y(C_E) - D'_X(C_E)C_Y+C_XC_Y\nonumber\\
&=& [C_{X\circ Y},\QQ]+C_X[C_Y,\QQ]+[C_X,\QQ]C_Y\nonumber\\
&=& 0\ . \nonumber 
\end{eqnarray}
\hfill $\qed$

\subsection{CV-structures on the tangent bundle}\label{s4.2}

In the situation assumed at the beginning of the chapter, one obtains an 
induced CV-structure on $TM$ and the structure of an F-manifold
$(M,\circ,e,E)$ with unit field $e$ and Euler field $E$.
On the tangent bundle one has additionally to the induced connection $D$
the Lie derivative of vector fields. In theorem \ref{t4.5}
the actions of $D_e, \Lie_e, D_E,\Lie_E$ and also 
$D_{\oooo e},\Lie_{\oooo e}, D_{\oooo E}, \Lie_{\oooo E}$ on the 
CV-structure on $TM$ will be discussed.

\begin{remarks}\label{t4.4}
(i) An $\OO_M$-linear map $\tm^{\otimes k}\to \tm^{\otimes l}$ 
is called a holomorphic $(k,l)$-tensor; a $\ciii_M$-linear map
$(\tmen)^{\otimes k}\to (\tmen)^{\otimes l}$ is called a $\ciii$
$(k,l)$-tensor. Vector fields in $\tmen$ can be identified with
$\ciii$ $(0,1)$-tensors by $\ciii_M\to \tmen$, $ a\mapsto a\cdot X$,
for $X\in \tmen$. We will also consider non $\ciii$ tensors; for example
a sesquilinear pairing on $TM$ is a $\ciii_M$-sesquilinear $(2,0)$-tensor.

(ii) Let $D$ be a $\ciii$-connection on $TM$ with $D^{(0,1)}\tm=0$.
Then one has for $X\in \tm$
\begin{eqnarray}\label{4.8}
\Lie_{\oooo X}(T)=D_{\oooo X}(T) 
\mbox{ \ \ for }T \mbox{ a }\ciii\mbox{ tensor}\ ,\\
\Lie_{\oooo X}(T) =0=D_{\oooo X}(T) \mbox{ \ \ for } T 
\mbox{ a holomorphic tensor.}\label{4.9}
\end{eqnarray}
But $D_{\oooo X}-\Lie_{\oooo X}$ acts nontrivially on non $\ciii$ tensors.
\end{remarks}

\begin{theorem}\label{t4.5}
Let $(TM,D,C,\www C,\kappa,h,g,\UU,\QQ)$ be a CV-structure on $TM$
such that $\tm$ is a locally free module over itself with respect to
the Higgs field $C$.

(a) There is a unique multiplication $\circ$ with $C_XC_Y=-C_{X\circ Y}$,
and there is a unique unit field $e$. Define $E:=\UU(e)$. Then 
$(M,\circ,e,E)$ is an F-manifold with Euler field $E$. Furthermore
\begin{eqnarray}\label{4.10}
\UU=E\circ \ ,\\
\Lie_e(\circ)=0\ ,\label{4.11}\\
\Lie_{\oooo e}(\QQ) = D_{\oooo e}(\QQ)= D_e(\QQ)=0\ .\label{4.12}
\end{eqnarray}

(b) We have the equivalences
\begin{eqnarray}\label{4.13}
&& D_e-\Lie_e=0 \mbox{ \ on }\ciii\mbox{ tensors }\\
&\iff & D_ee=0 \label{4.14}\\
&\iff & \Lie_e(h)=0 \iff \Lie_{\oooo e}(h)=0 \label{4.15}\\
&\iff &  \Lie_e(\kappa)=0 \iff \Lie_{\oooo e}(\kappa )=0 \ .\label{4.16}
\end{eqnarray}
If the six statements in the equivalences are true then
\begin{eqnarray}\label{4.17}
\Lie_e(g)=0 \mbox{ \ and \ } \Lie_e(\QQ)=0
\end{eqnarray}
and then all the tensors $\circ,e,\kappa,g, h,\QQ$ are invariant under the
flows of $e$ and $\oooo e$.

(c) Fix some $d\in\R$ and define the map
\begin{eqnarray}\label{4.18}
\www\QQ := D'_E - \Lie_E - \frac{2-d}{2}\id :\tmen\to \tmen\ .
\end{eqnarray}
Then 
\begin{eqnarray}\label{4.19}
D'(\UU)-[C,\www\QQ] + C=0\ ,\\
\Lie_E(\www\QQ)=D_E(\www\QQ)\ , \ \ 
\Lie_{\oooo E}(\www\QQ)=D_{\oooo E}(\www\QQ)\ , \label{4.20}
\end{eqnarray}
and the two equivalences
\begin{eqnarray}\label{4.21}
\Lie_{E-\oooo E}(h)=0 &\iff& 
h(\www\QQ X,Y)=h(X,\www\QQ Y) \mbox{ for }X,Y\in\tmen\ ,\\
  \Lie_{E}(g)=(2-d)g &\iff& 
g(\www\QQ X,Y)=-g(X,\www\QQ Y) \mbox{ for }X,Y\in\tmen\label{4.22} \qquad
\end{eqnarray}
hold.  If $\Lie_{E-\oooo E}(h)=0$ and $\Lie_{E}(g)=(2-d) g$
are true then $(TM, D, C, \www C, \kappa, h, g, \UU, \www\QQ)$ 
is a CV-structure and the following equations hold,
\begin{eqnarray}\label{4.23}
\Lie_{E-\oooo E}(\www\QQ)=D_{E-\oooo E}(\www\QQ)=0 = D_{E-\oooo E}(\QQ)&,&\\
\Lie_E(\www\QQ)=D_E(\www\QQ)=\frac{1}{2}\Lie_{E+\oooo E}(\www\QQ)
=\frac{1}{2}D_{E+\oooo E}(\www\QQ)&& \label{4.24} \\
= [\UU,\kappa \UU\kappa] = D_E(\QQ) = \frac{1}{2}D_{E+\oooo E}(\QQ)&.&
\nonumber
\end{eqnarray}
Especially, then
the tensors $h$ and $\www\QQ$ are invariant under the flow of $E-\oooo E$.
\end{theorem}

{\it Proof.}
(a) Lemma \ref{4.1} and lemma \ref{4.3} show everything except for
\eqref{4.10}--\eqref{4.12}. $\Lie_e(\circ)=0$ follows from \eqref{4.2}
for $X=Y=e$. The calculation
\begin{eqnarray}\nonumber
\UU(X)=-\UU C_X(e) = -C_X\UU(e)=X\circ E
\end{eqnarray}
shows $\UU=E\circ$. The identity 
\eqref{4.8} together with \eqref{2.58}, \eqref{2.50} , 
\eqref{2.57} gives
\begin{eqnarray}\nonumber
\Lie_{\oooo e}(\QQ) = D''_{\oooo e}(\QQ) = -\kappa D'_e(\QQ)\kappa
= \kappa [C_e,\kappa\UU\kappa]\kappa =0\ .
\end{eqnarray}
(b) Fix a holomorphic vector field $X\in \tm$. Then for $Y,Z\in \tmen$
\begin{eqnarray}\nonumber
\Lie_X(h)(Y,Z) &=& X\, h(Y,Z)-h(\Lie_XY,Z)-h(Y,\Lie_{\oooo X}Z)\\
&=& h((D_X-\Lie_X)Y,Z)\ ,\label{4.25}
\end{eqnarray}
using $D(h)=0$ and $\Lie_{\oooo X}Z = D_{\oooo X}Z$. Similarly
\begin{eqnarray}
 \Lie_{\oooo X}(h)(Y,Z) &=& h(Y,(D_X-\Lie_X)Z) \ ,\label{4.26} \\
\Lie_{X}(g)(Y,Z) &=& g((D_X-\Lie_X)Y,Z) + 
         g(Y,(D_X-\Lie_X)Z) \ ,\label{4.27} \\
\Lie_{X}(\kappa)(Y) &=& (\Lie_X-D_X)(\kappa(Y))\ ,\label{4.28} \\
\Lie_{\oooo X}(\kappa)(Y) &=& \kappa ((D_X-\Lie_X) Y) \ .\label{4.29}
\end{eqnarray}
These statements for $X=e$ and the calculation 
\begin{eqnarray}\nonumber
(D_e-\Lie_e)(Y) &=& D_e(Y\circ e) - [e,Y]\circ e \\
&=& D_e(Y\circ e)-D_e(Y\circ )(e) + D_Y(e\circ )(e) \nonumber\\
&=& Y\circ D_ee \label{4.30}
\end{eqnarray}
show the equivalences in \eqref{4.13}--\eqref{4.16}.
For \eqref{4.17} use \eqref{4.27} and $D_e(\QQ)=0$.

\noindent
(c) The following calculation shows \eqref{4.19}. It uses 
$\Lie_E(\circ )=\circ$ and $D'(C)=0$.
\begin{eqnarray}\label{4.31}
&& D'_X(\UU)-[C_X,\www\QQ] +C_X \\
&=& D'_X(E\circ ) + [X\circ ,D'_E-\Lie_E]-X\circ \nonumber \\
&=& D'_X(E\circ) - D'_E(X\circ ) + \Lie_E(X\circ ) - X\circ \nonumber \\
&=& [X,E]\circ + \Lie_E(X) \circ + \Lie_E(\circ)(X,\cdot ) -X\circ \nonumber\\
&=& 0\ .\nonumber
\end{eqnarray}
One can rewrite $\www\QQ$ as $\www\QQ= D_{E\pm\oooo E} -\Lie_{E\pm\oooo E}
-\frac{2-d}{2}\id$. Then
\begin{eqnarray}\label{4.32}
D_{E\pm \oooo E}(\www\QQ)= [\Lie_{E\pm\oooo E},D_{E\pm \oooo E}] 
= \Lie_{E\pm \oooo E}(\www\QQ)\ .
\end{eqnarray}
This shows \eqref{4.20}. The two equivalences \eqref{4.21} and \eqref{4.22}
follow from \eqref{4.25}--\eqref{4.27} for $X=E$.

Now suppose that $\Lie_{E-\oooo E}(h)=0$ and $\Lie_E(g)=(2-d)g$ hold.
Then $g=h(\cdot ,\kappa \cdot)$ shows
\begin{eqnarray}\label{4.33}
\Lie_{E-\oooo E}(\kappa) = (d-2)\kappa\ ,
\end{eqnarray}
and $\www\QQ +\kappa\www\QQ \kappa =0$ follows from the right hand sides of
\eqref{4.21} and \eqref{4.22}. 
In order to see that one obtains a CV-structure with 
$\www\QQ$ instead of $\QQ$ it rests to make the following calculation
for $X\in \tm$. 
It uses curvature properties of the $(DC\www C)$-structure,
$D(\kappa)=0$, \eqref{4.8}, \eqref{4.33} and
$[\Lie_Y,\Lie_Z]=\Lie_{[Y,Z]}$ for $Y,Z\in \tmen$. 
\begin{eqnarray}\label{4.34}
&& D'_X(\www\QQ) + [C_X,\kappa \UU\kappa] \\
&=& [D'_X,D'_E-\Lie_E] - [C_X,\www C_{\oooo E}] \nonumber\\
&=& D'_{[X,E]} - [D'_X,\Lie_E] + [D'_X,D''_{\oooo E}]  \nonumber\\
&=& D'_{[X,E]} - [D'_X,\Lie_{E-\oooo E}]  \nonumber\\
&=& D'_{[X,E]} - [\kappa D''_{\oooo X}\kappa,\kappa\Lie_{\oooo E-E}\kappa]
  \nonumber\\
&=& D'_{[X,E]} - \kappa [\Lie_{\oooo X},\Lie_{\oooo E-E}]\kappa  \nonumber\\
&=& D'_{[X,E]} - \kappa\Lie_{[\oooo X,\oooo E-E]}\kappa  \nonumber\\
&=& D'_{[X,E]} - \kappa D''_{\oooo{[X,E]}}\kappa  \nonumber\\
&=& 0\ . \nonumber
\end{eqnarray}
The properties of $\QQ$ in \eqref{4.23} and \eqref{4.24} were already 
discussed in \eqref{2.66b}. The same arguments together with 
\eqref{4.20} give the properties of $\www\QQ$ in \eqref{4.23} and \eqref{4.24}.
\hfill $\qed$

\bigskip
Of course one wishes that the equivalent statements in 
\eqref{4.13}--\eqref{4.16} and \eqref{4.21}, \eqref{4.22} are all satisfied
and that $\QQ=\www\QQ$. Furthermore one may wish that the holomorphic
data $(M,\circ,e,E,g)$ yield a Frobenius manifold.
If one starts with a CV-structure on a vector bundle $K\to M$ such that
$\OO(K)$ is a free $\tm$-module of rank 1, one has to choose a 
generator $\zeta\in \OO(K)$ very carefully. In section \ref{s5.4}
it is discussed how to do this.
Here we make a definition, motivated by theorem \ref{t4.5}
and by \cite{CV1}\cite{Du2}.

\begin{definition}\label{t4.6}
A CDV-structure on a manifold $M$ is a CV-structure
$(TM,D,C,\www C,\kappa,h,g,\UU,\QQ)$ together with a Frobenius manifold
structure $(M,\circ,e,E,g)$ such that $C_X Y=-X\circ Y$, 
$E=\UU (e)$, $\Lie_E(g)=(2-d)g$, and
\begin{eqnarray}\label{4.35}
\QQ=D_E-\Lie_E-\frac{2-d}{2}\id
\end{eqnarray}
for some $d\in \R$ and such that the equivalent statements in 
\eqref{4.13}--\eqref{4.16} hold.
\end{definition}

\begin{remarks}\label{t4.7}
(i) The definition of a Frobenius manifold is recalled in definition
\ref{t5.10}.
For a CDV-structure one needs besides the CV-structure with
\eqref{4.35}, $C_XC_Y=-C_{X\circ Y}$, $C_Xe=-X$, 
$E=\UU(e)$, $D_ee=0$, only that
the metric $g$ is flat and that its Levi--Civita connection 
$\nnn^g$ satisfies the potentiality $\nnn^g (C)=0$.
Then $\nnn^ge=0$ and $\Lie_E(g)=(2-d)g$ will follow automatically,
as well as $\Lie_{E-\oooo E}(h)=0$.

(ii) If one starts with a Frobenius manifold one needs for a 
CDV-structure only one new ingredient, a real structure $\kappa$
(or a hermitian pairing $h$). All the other structures,
$D,\www C,h$ (or  $\kappa$), $\QQ$ can be derived from it.
But it has to satisfy a lot of properties.
They are given in definition \ref{t1.2}. 
In order to see the equivalence of definition \ref{t1.2} and
definition \ref{t4.6} one needs \eqref{4.21} and \eqref{4.22}.
\end{remarks}

\clearpage

\section{Frobenius manifolds and $tt^*$ geometry}\label{s5}
\setcounter{equation}{0}

\noindent
In chapter \ref{s2} interesting geometry on a vector bundle
$K\to M$ arose from a vector bundle $H\to \pmmm$ with a flat
connection on $\csmmm$. This connection induced additional structure
on $H|_\nmmm=K$ and on $H|_\immm$; and $K$ and $H|_\immm$ could
be identified under the assumption that $H$ is a family of 
trivial bundles on $\P^1$. The same ideas will be used here.
But instead of $(T\www T)$-structures here we will consider
`$(TL)$-structures'. The antiholomorphic twin at $\immm$ of a 
pole of Poincar\'e rank 1 will be replaced by a (holomorphic)
logarithmic pole. Therefore the structure here is simpler in
two aspects: it is holomorphic, and the logarithmic pole
induces a simpler structure on $K$ than (an antiholomorphic twin of)
a pole of Poincar\'e rank 1. Still the whole discussion is very similar
and many pieces coincide. 

In section \ref{s5.1} logarithmic poles and the induced structures
are treated. In section \ref{s5.2} the ideas above are carried out.
In section \ref{s5.3} the resulting structures are lifted to the
tangent bundle and give Frobenius manifolds. In section \ref{s5.4}
CDV-structures are obtained by combining the discussion here
with chapter \ref{s4}.

\subsection{$(L)$-structures and logarithmic poles}\label{s5.1}

In chapter \ref{s2} meromorphic connections were considered on 
$\pmmm$ with poles along $\nmmm$ and $\immm$. It will be useful
to separate inessential choices of coordinates from essential 
choices.

\begin{definition}\label{t5.1}
Let $M_1$ be a complex manifold, $M_2\subset M_1$ a submanifold of
codimension 1 and $f:M_1\to \C$ a holomorphic submersion with 
$M_2=f^{-1}(0)$. Let $H\to M_1$ be a holomorphic vector bundle.

(a) A holomorphic family of flat connections $\nnn$ on $H|_{f^{-1}(t)}$,
$t\in \C^*$, is said to have a {\it pole of Poincar\'e rank $r\in \Z_{\geq 0}$}
if
\begin{eqnarray}\label{5.1}
\nnn : \OO(H) \to \frac{1}{f^r}\oomm^1_{M_1/\C} \otimes \OO(H)\ .
\end{eqnarray}
If $r=0$ then the pair $(H,\nnn)$ is called an $(L)$-structure
(`L' for Logarithmic); if $r=1$ then it is called a $(T)$-structure
(cf. definition \ref{t2.3}).

(b) A flat connection $\nnn$ on $H|_{M_1-M_2}$ has a {\it pole of Poincar\'e
rank $r$} along $M_2$ if 
\begin{eqnarray}\label{5.2}
\nnn : \OO(H) \to \frac{1}{f^r}\oomm^1_{M_1}(\log M_2) \otimes \OO(H)\ .
\end{eqnarray}
If $r=0$ then the pair $(H,\nnn)$ is called an $(LE)$-structure.
A pole of Poincar\'e rank 0 is called a logarithmic pole.

(c) Let $M_1=\C\times M, \ M_2=\{0\}\times M$, and $f=z$. Fix $w\in \Z$.
The triple $(H,\nnn,P)$ is an $(LEP(w))$-structure if $(H,\nnn)$ is an
$(LE)$-structure and $P$ is a $\C$-bilinear pairing
\begin{eqnarray}\label{5.3}
P:H_{(z,t)}\times H_{(-z,t)}\to\C \mbox{ \ for any }(z,t)\in \csmmm
\end{eqnarray}
with all properties in definition \ref{t2.12} (d) except \eqref{2.34}.
That means, $P$ is nondegenerate, $(-1)^w$-symmetric, $\nnn$-flat, and
on an open subset $U_1\times U_2\subset\cmmm$ it extends to a nondegenerate
pairing
\begin{eqnarray} \label{5.4}
\qquad P:\OO(H)(U_1\times U_2)\times \OO(H)((-U_1)\times U_2) \to 
z^w\OO_\cmmm (U_1\times U_2) \\
(a, b) \mapsto ((z,t)\mapsto P(a(z,t),b(-z,t)))\ .\nonumber
\end{eqnarray} 

(d) Let $M_1=\C\times M, \ M_2=\{0\}\times M$, and $f=z$.
The triple $(H,\nnn,P)$ is an $(LP(w))$-structure if $(H,\nnn)$ is
an $(L)$-structure and $P$ is a pairing as in (c).
\end{definition}

\begin{remarks}\label{t5.2}
(i) The sheaf $\oomm^k_{M_1}(\log M_2)$ of logarithmic $k$-forms is 
$\oomm^k_{M_1}+\oomm^{k-1}_{M_1}\land \frac{\ddd f}{f}$.
Part (b) is independent of the choice of $f$, but in part (a) 
the function $f$ is essential.

(ii) We will also consider $LEP(w)$-structures with $M_1=\pnmmm$,
$M_2=\immm$, $f=\eezz$. They are obtained from $(LEP(w))$-structures
on $\cmmm \supset \nmmm$ by the map $\pmmm\to\pmmm$, 
$(z,t)\mapsto (\eezz,t)$. The structures in lemma \ref{t5.3} are shifted
with it from $\nmmm$ to $\immm$.

(iii) In the same way $(LP(w))$-structures with
$M_1=\pnmmm$, $M_2=\immm$, and $f=\eezz$ are defined.
\end{remarks}

\begin{lemma}\label{t5.3}
Let $M_1,M_2,f$ and $H\to M_1$ be as in definition \ref{t5.1}.

(a) Let $(H,\nnn)$ be an $(L)$-structure. Then $\nnn$ restricts to 
a flat connection $\nnn^{res}$ on $H|_{M_2}$, the residual connection.

(b) Let $(H,\nnn)$ be an $(LE)$-structure. It contains an $(L)$-structure;
let $\nnn^{res}$ be its residual connection  on $H|_{M_2}$. 
There is a $\nnn^{res}$-flat endomorphism $\Nu^{res}$ on $H|_{M_2}$,
the residue endomorphism. Locally it is defined by
\begin{eqnarray}\label{5.4b}
\Nu^{res}([a]) := \nnn_{f\cdot \zeta}a\mod f\cdot\OO(H)\ ,
\end{eqnarray}
where $a\in \OO(H)$, $[a]\in \OO(H|_{M_2})$, and $\zeta\in \TT_{M_1}$ is a 
vector field with $\zeta(f)=1$.

(c) Let $(H\to\cmmm,\nnn,P)$ be an $(LEP(w)$-structure.
On $K:=H|_\nmmm$ one has $\nnn^{res}$ and $\Nu^{res}$ as in (b).
There is a symmetric nondegenerate holomorphic pairing $g$ on $K$
defined by
\begin{eqnarray}\label{5.5}
g([a],[b]) := (z^{-w}P(a,b)) \mod z\OO(H)
\end{eqnarray}
for  $[a],[b]\in \OO(K)$, and lifts $a,b\in \OO(H)$. It is $\nnn^{res}$-flat
and satisfies
\begin{eqnarray}\label{5.6}
g(\Nu^{res}[a],[b])+g([a],\Nu^{res}[b]) = w g([a],[b])\ .
\end{eqnarray}
\end{lemma}

{\it Proof.}
(a) Trivial.

(b) Suppose for simplicity $M_1=\cmmm,M_2=\nmmm$ and $f=z$. Then one can
choose $\zeta=\dz$ and lift vector fields $X\in \tm$ canonically to $M_1$
such that $[X,\dz]=0$. For $a\in\OO(H)$, $[a]\in \OO(H|_\nmmm)$, $X\in \tm$
\begin{eqnarray}\label{5.7}
(\nnn^{res}_X\Nu^{res})([a]) = [\nnn_X,\nnn_{\zdz}](a) \mod z\OO(H)
= 0\ .
\end{eqnarray}
Therefore $\Nu^{res}$ is $\nnn^{res}$-flat.

(c) The pairing $g$ is symmetric, nondegenerate and holomorphic as in 
lemma \ref{t2.14} (a). For $a,b\in \OO(H)$ one finds with $\nnn(P)=0$
\begin{eqnarray}\label{5.8}
&& g(\Nu^{res}[a],[b]) + g([a],\Nu^{res}[b]) -wg([a],[b]) \\
&=& z^{-w}P(\nnn_\zdz a,b) + z^{-w}P(a,\nnn_\zdz b) 
-w\cdot z^{-w}P(a,b) \mod z\OO_\cmmm \nonumber\\
&=& \zdz (z^{-w}P(a,b))\mod z\OO_\cmmm =0 \ .\nonumber
\end{eqnarray}
A similar calculation shows $\nnn^{res}(g)=0$.
\hfill $\qed$

\begin{remarks}\label{t5.4}
(i) The residue endomorphism of a connection with logarithmic pole is 
independent of a function $f$ as in definition \ref{5.1} (b).
But the residual connection depends on it.

(ii) Given a flat bundle $\www H\to\csmmm$, the extensions to
vector bundles $H\to \cmmm$ whose connections have logarithmic poles
along $\nmmm$ are rather simple objects and can be classified and
encoded by certain filtrations. See section \ref{s7.5}.
\end{remarks}

\subsection{$(TLEP)$-structures and Frobenius type structures}\label{s5.2}

\begin{definition}\label{t5.5}
Let $H\to \pmmm$ be a holomorphic vector bundle such that 
$H|_{\P^1\times\{t\}}$ is a trivial bundle for any $t\in M$.

(a) A $(trTL)$-structure is a pair $(H,\nnn)$ whose restriction to
$\cmmm$ is a $(T)$-structure (definition \ref{2.3} (i)) and whose
restriction to $\pnmmm$ is an $(L)$-structure.

(b) A $(trTLP(w))$-structure is a triple $(H,\nnn,P)$ whose restriction to
$\cmmm$ is a $(TP(w))$-structure (remark \ref{2.13} (i)) and whose
restriction to $\pnmmm$ is an $(LP(-w))$-structure (remark \ref{5.2} (iii)).

(c) A $(trTLEP(w))$-structure is a triple $(H,\nnn,P)$ whose restriction to
$\cmmm$ is a $(TEP(w))$-structure (remark \ref{2.13} (i)) and whose
restriction to $\pnmmm$ is an $(LEP(-w))$-structure (remark \ref{5.2} (ii)).
\end{definition}

Theorem \ref{t5.7} and corollary \ref{t5.8} will say that these structures
induce on $K:=H|_\nmmm$ the structures in definition \ref{t5.6}
and that they are equivalent to them.

\begin{definition}\label{t5.6}
(a) A $(\nnn^rC)$-structure is a holomorphic vector bundle $K\to M$ 
together with a flat holomorphic connection $\nnn^r$ and a Higgs field
$C:\OO(K)\to \oomm^1_M\otimes \OO(K)$ such that $\nnn^r(C)=0$;
that means (cf. \eqref{2.6})
\begin{eqnarray}\label{5.9a}
\nnn^r_X(C_Y)-\nnn^r_Y(C_X)-C_{[X,Y]} =0 \mbox{ \ for }X,Y\in\tm
\end{eqnarray}
(the `r' in $\nnn^r$ stands for `residual' or `restriction').

(b) A $(\nnn^rCg)$-structure is a $(\nnn^rC)$ structure together with a
symmetric nondegenerate $\nnn^r$-flat pairing $g$ on $K$ with
\begin{eqnarray}\label{5.9b}
g(C_Xa,b)=g(a,C_Xb)
\end{eqnarray}
for $X\in \tm,a,b\in \OO(K)$.

(c) A Frobenius type structure on a holomorphic vector bundle $K\to M$
consists of a flat holomorphic connection $\nnn^r$ on $K$, a Higgs field
$C$ on $K$ with $\nnn^r(C)=0$, an $\OO_M$-linear endomorphism $\UU$
of $K$ with $[C,\UU]=0$, an $\OO_M$-linear $\nnn^r$-flat endomorphism
$\Nu$ of $K$ with 
\begin{eqnarray}\label{5.10}
\nnn^r(\UU)-[C,\Nu]+C=0\ ,
\end{eqnarray}
and a symmetric nondegenerate $\nnn^r$-flat pairing $g$ on $K$ with
\begin{eqnarray}\label{5.11}
g(C_Xa,b)&=&g(a,C_Xb) \  ,\\
g(\UU a,b)&=&g(a,\UU b)\  ,\label{5.12}\\
g(\Nu a,b)&=& -g(a,\Nu b) \label{5.13}
\end{eqnarray}
for $X\in \tm ,a,b\in \OO(K)$.
\end{definition}

\begin{theorem}\label{t5.7}
Fiz $w\in \Z$. There is a one-to-one correspondence between 
$(trTLEP(w))$-structures and Frobenius type structures on holomorphic
vector bundles. It is given by the steps in (a) and (b).
They are inverse to one another.

(a) Let $(K \to M,\nnn^r,C,g,\UU,\Nu)$ be a Frobenius type structure
on $K$. Let $\pi :\pmmm \to M$ be the projection. 
Define $H:=\pi^*K$, and let $\psi_z:H_{(z,t)}\to K_t$ for $z\in\P^1$
be the canonical projection.
Extend $\nnn^r,C,g,\UU,\Nu$ canonically to $H$. Define
\begin{eqnarray}\label{5.14}
\nnn:= \nnn^r +\eezz C + (\eezz\UU-\Nu+\frac{w}{2}\id )\frac{\ddd z}{z}\ .
\end{eqnarray}
Define a pairing
\begin{eqnarray}\label{5.15}
P:H_{(z,t)}\times H_{(-z,t)}&\to& \C \mbox{ \ \ \ for }(z,t)\in \csmmm\\
(a,b) &\mapsto & z^wg(\psi_z a,\psi_{-z}b)\ .\nonumber
\end{eqnarray}
Then $(H,\nnn,P)$ is a $(trTLEP(w))$-structure.

(b) Let $(H,\nnn,P)$ be a $(trTLEP(w))$-structure. Define $K:=H|_\nmmm$.
Let $g$ and $\UU$ on $K$ come from the $(TEP(w))$-structure on 
$H|_\cmmm$ as in lemma \ref{t2.14}.
Let $\nnn^{res}$ and $\Nu^{res}$ be the residual connection and the
residue endomorphism on $H|_\immm$.

Because $H|_{\P^1\times \{t\}}$ is a trivial bundle for any $t\in M$,
there is a canonical projection $\psi:H\to K$, and the bundles
$K$ and $H|_\immm$ are canonically isomorphic. Structure on 
$H|_\immm$ can be shifted to $K$. Let $\nnn^r$ on $K$ be the shift of
$\nnn^{res}$ and let $\Nu$ on $K$ be the shift of 
$\Nu^{res}+\frac{w}{2}\id$. 
Then $(K\to M,\nnn^r,C,g,\UU,\Nu)$ is a Frobenius type structure and
\eqref{5.14} holds.
\end{theorem}

{\it Proof.}
Most of the proof is similar to the proof of theorem \ref{t2.19},
but simpler.

(a) $\nnn^2=0$ holds because it splits into $(\nnn^r)^2=0$, 
$\nnn^r(C)=0$, $C^2=0$, $[C,\UU]=0$, $\nnn^r(\Nu)=0$, and \eqref{5.10}.
For the summand $C$ in \eqref{5.10} observe that $\nnn^r$ in \eqref{5.14}
contains a covariant derivative $\nnn^r_\dz$.

The pairing $P$ is $(-1)^w$-symmetric because $g$ is symmetric.
The flatness $\nnn(P)=0$ follows with a calculation similar to
\eqref{2.82} from $\nnn^r(g)=0$ and \eqref{5.11}--\eqref{5.13}.
By definition $P$ satisfies the limit behaviour along $\immm$ which is
required for an $(LEP(-w))$-structure on $\pnmmm$
(formula \eqref{5.4}). The same holds for the 
$(TEP(w))$-structure on $\cmmm$.

(b) Lift $\nnn^r,C,\UU$, and $\Nu$ canonically to $H$, using $\psi$.
Consider
\begin{eqnarray}\label{5.17}
\Delta := \nnn - (\nnn^r+\eezz C+(\eezz\UU-\Nu+\frac{w}{2}\id)\dzz)\ .
\end{eqnarray}
The $\OO_\pmmm$-linear map $\Delta_X$ for $X\in \tm$ maps sections
in $\pi_*\OO(H)$ to sections which have no pole along $\nmmm$ because
of the definition of $C$ and which vanish along $\immm$ because
of the definition of $\nnn^r$. Therefore $\Delta_X=0$. The same
holds for $\Delta_\zdz$, using the definitions of $\UU$ and $\Nu$.
Here observe $\dzz = -\frac{\ddd (1/z)}{1/z}$. This shows $\Delta=0$
and \eqref{5.14}.

Now the flatness $\nnn^2=0$ splits into $(\nnn^r)^2=0$,
$\nnn^r(C)=0$, $C^2=0$, $[C,\UU]=0$, $\nnn^r(\Nu)=0$, and \eqref{5.10}.
The pairing $g$ is symmetric and nondegenerate and satisfies
\eqref{5.11} and \eqref{5.12} by lemma \ref{t2.14}.
Consider sections $a,b\in \OO(K)$ and their canonical lifts 
$\www a,\www b\in \pi_*\OO(H)$ to $H$. Then by definition of a 
$(trTLEP(w))$-structure
\begin{eqnarray}\label{5.18}
P(\www a,\www b)\in z^w\cdot \pi_*\OO_\pmmm = z^w\cdot \OO_M\ .
\end{eqnarray}
Therefore $P(\www a,\www b)=z^wg(a,b)$. The lift of $g$ to $H$ restricts
on $H|_\immm$ to the pairing in lemma \ref{t5.3} (c) of the 
$(LEP(-w))$-structure at $\immm$. By lemma \ref{5.3} (c) $g$ is $\nnn^r$-flat
and satisfies \eqref{5.13}.
\hfill $\qed$

\begin{corollary}\label{t5.8}
There is a one-to-one correspondence between $(trTL)$-structures and
$(\nnn^rC)$-structures. For a fixed $w\in\Z$ there is a one-to-one
correspondence between $(trTLP(w))$-structures and $(\nnn^rCg)$-structures.
In both cases one can read off from theorem \ref{t5.7} how to pass
from one structure to the other.
\end{corollary}

{\it Proof.} A part of the proof of theorem \ref{t5.7} gives the
proof.
\hfill $\qed$

\begin{remark}\label{t5.9}
Usually a $(TEP(w))$-structure comes from geometry, e.g. from oscillating
integrals for singularities, and one wants to extend it to a 
$(trTLEP(w))$-structure.
The choice of some extension to $\immm$ with a logarithmic pole is cheap.
But it is nontrivial to find extensions such that the restrictions
$H|_{\P^1\times \{t\}}$ of the total bundle are trivial.
This is called a Birkhoff problem (with parameters). A distinguished
way to solve it takes the point of view of filtrations, which are 
associated to logarithmic extensions along $\immm$ and to restrictions
$(H|_{\C\times \{t\}},\nnn)$ for special parameters $t\in M$.
See section \ref{s7.5}.
\end{remark}

\subsection{Frobenius manifolds}\label{s5.3}

A Frobenius manifold $M$ is a Frobenius type structure on the tangent
bundle with some additional properties. By theorem \ref{t5.7}, for any
$w\in \Z$ 
it gives rise to a meromorphic connection and a pairing on the lifted
tangent bundle $\pi^*TM$, where $\pi:\P^1\times M\to M$, 
which together form a $(trTLEP(w))$-structure.

The connections are sometimes called first structures connections
\cite{Man}. They are most important. They are a gate to richer structures.
They also give a frame to construct Frobenius manifolds: One has
to construct $(trTLEP(w))$-structures on abstract bundles and shift
them with some diligence to the lifted tangent bundle.

Here we recall the definition of Frobenius manifolds and their
first structure connections and describe this frame for constructing
Frobenius manifolds. Of course, this is essentially well known. 
The presentation here is influenced by papers of Barannikov,
Dubrovin, Manin, Sabbah, and K. Saito.

\bigskip

\begin{definition}\label{t5.10}
(a) A {\it Frobenius manifold} $(M,\circ,e,E,g)$ is a complex manifold $M$
of dimension $\geq 1$ with a multiplication $\circ$ on the holomorphic
tangent bundle $TM$ , a unit field $e$, an {\it Euler field} $E$ and
a symmetric nondegenerate $\OO_M$-bilinear pairing $g$ on $TM$ with the
following properties.\\
The metric $g$ is multiplication invariant, that means,
\begin{eqnarray}\label{5.19}
g(X\circ Y,Z)=g(Y,X\circ Z)\mbox{ \ for }X,Y,Z\in \tm\ ;
\end{eqnarray}
the Levi--Civita connection $\nnn^g$ of the metric $g$ is flat;
together with the Higgs field 
$C:\tm\to \oomm^1_M\otimes \tm$ with $C_XY:=-X\circ Y$
it satisfies the potentiality condition $\nnn^g(C)=0$; 
the unit field $e$ is $\nnn^g$-flat;
and the Euler field satisfies $\Lie_E(\circ)=\circ$ and 
$\Lie_E(g)=(2-d)\cdot g$ for some $d\in\C$.

(b) A {\it Frobenius manifold without Euler field} consists of $(M,\circ,e,g)$
as in (a) with all properties without those involving $E$.

(c) A {\it Frobenius manifold without metric and Euler field}
is a manifold $(M,\circ,e)$ with multiplication $\circ$, unit field $e$, 
and Higgs field $C$ as in (a) together with a torsion free flat
connection $\nnn^g$ with $\nnn^g(C)=0$ and $\nnn^g(e)=0$.
\end{definition}

\begin{lemma}\label{t5.11}
Let $(M,\circ,e,E,g)$ be a Frobenius manifold with $d\in \C$, Higgs field
$C$, and Levi-Civita connection $\nnn^g$ as in definition \ref{t5.10}.
Define $\UU:=E\circ$ and 
\begin{eqnarray}\label{5.20}
\Nu:= \nnn^g_\bullet E - \frac{2-d}{2}\id = \nnn^g_E-\Lie_E-
\frac{2-d}{2}\id \ .
\end{eqnarray}

(a) Then $(TM,\nnn^g,C,g,\UU,\Nu)$ is a Frobenius type structure on 
$TM$.

(b) Fix $w\in \Z$. Let $(\pi^*TM,\nnn,P)$ be the $(trTLEP(w))$-structure
which corresponds to this Frobenius type structure by theorem 
\ref{t5.7} (a). Then the canonical isomorphism $TM\to \pi^*TM|_\immm$
maps $\nnn^g$ to the residual connection of $\nnn$ at $\immm$ and
$\Nu-\frac{w}{2}\id$ to the residue connection of $\nnn$ at $\immm$.
The unit field satisfies
\begin{eqnarray}\label{5.21}
\nnn^g e=0 \mbox{ \ and \ } \Nu e=\frac{d}{2}e\ .
\end{eqnarray}
\end{lemma}

{\it Proof.}
$\Lie_E(\circ)(e,e)=e$ shows $[e,E]=e$ and
\begin{eqnarray}\label{5.22}
\Nu(e) = \nnn^g_e E-\frac{2-d}{2}e = [e,E]-\frac{2-d}{2}e=\frac{d}{2}e\ .
\end{eqnarray}
The calculation \eqref{4.31} with $D'$ replaced by $\nnn^g$ and
$\www\QQ$ replaced by $\Nu$ shows \eqref{5.10}. Because of 
$\nnn^g g=0$ equation \eqref{5.13} is equivalent to $\Lie_E(g)=(2-d)\cdot g$.
In order to see the flatness $\nnn^g(\Nu)=0$ of $\Nu$ one has to
show the symmetry and the skew symmetry
\begin{eqnarray}\label{5.23}
(\nnn^g_X\Nu)(Y)=(\nnn^g_Y\Nu)(X) \mbox{ \ \ for } X,Y\in \tm\ ,\\
g((\nnn^g_X\Nu)(Y),Z)=-g((\nnn^g_X\Nu)(Z),Y) \mbox{ \ \ for }
X,Y,Z\in \tm \ , \label{5.24}
\end{eqnarray}
and use both three times in order to see
$g((\nnn^g_X\Nu)(Y),Z)=-g((\nnn^g_X\Nu)(Y),Z)=0$.
The other statements are clear.
\hfill $\qed$

\bigskip

These $(trTLEP(w))$-structures on $\pi^*TM$ are called
{\it first structure connections} of the Frobenius manifold.

The following gives an inverse to lemma \ref{t5.11} and a frame 
for constructing Frobenius manifolds.

\begin{theorem}\label{t5.12}
Let $(H\to \pmmm,\nnn,P)$ be a $(trTLEP(w))$-structure and let
$(K=H|_\nmmm,\nnn^r,C^K,g^K,\UU^K,\Nu^K)$ be the corresponding
Frobenius type structure on $K$.

Let $\zeta\in \OO(K)$ be a global section such that 
$v:=-C^{K}_\bullet \zeta:\tm\to \OO(K)$ is an isomorphism and such that
$\nnn^r\zeta=0$ and $\Nu^K\zeta=\frac{d}{2}\zeta$ for some $d\in \C$.
Shift the Frobenius type structure on $K$ with the isomorphism
$v^{-1}$ to a Frobenius type structure 
$(TM,\nnn^g,C,g,\UU,\Nu)$ on $TM$. 

Then this is the Frobenius type
structure of a Frobenius manifold. Especially, $\nnn^g$ is torsion free,
$e=v^{-1}(\zeta)$ is the unit field,
$\Lie_E(g)=(2-d)\cdot g$, and $\Nu$ satisfies \eqref{5.20}.
\end{theorem}

{\it Proof.}
Define multiplication $\circ$ and unit field $e$ as in lemma \ref{t4.1}
and define $E:=\UU(e)$. By lemma \ref{4.3} for $D'=\nnn^{(1,0)}$ 
and $\QQ=\Nu$ the tuple $(M,\circ,e,E)$ is an F-manifold with Euler
field. The unit field $e$ is $\nnn^g$-flat because of $\nnn^r\zeta=0$.

With the isomorphism $\OO(K)=\pi_*\OO(H)$ elements of $\OO(K)$ will also
be considered as sections in $H$.
By \eqref{5.14} for $X,Y\in \tm$
\begin{eqnarray}\label{5.24b}
\nnn_Xv(Y)&=& \nnn^r_Xv(Y)+\eezz C_X^Kv(Y)
=v(\nnn^g_XY)-\eezz v(X\circ Y)\ ,\\
\nnn_\zdz v(Y) &=& \eezz v(E\circ Y)-v((\Nu-\frac{w}{2}\id)Y)\ .\label{5.25}
\end{eqnarray}
The connection $\nnn^g$ is torsion free because of the calculation
\begin{eqnarray}\label{5.26}
&& v(\nnn^g_{[X,Y]}Z)-\eezz v([X,Y]\circ Z)\\
&=& \nnn_{[X,Y]}v(Z) = [\nnn_X,\nnn_Y]v(Z) \nonumber \\
&=& v([\nnn^g_X,\nnn^g_Y]Z) - 
\eezz v((\nnn^g_X(Y\circ)-\nnn^g_Y(X\circ))(Z))\ .
\nonumber
\end{eqnarray}
Therefore it is the Levi--Civita connection of $g$.
The following calculation shows \eqref{5.20}. It uses $\nnn^ge=0$ and 
$\Nu e=\frac{d}{2}e$, which follows from $\Nu^K\zeta=\frac{d}{2}\zeta$.
\begin{eqnarray}\label{5.27}
&& v(\nnn^g_X E) = \nnn_Xv(E) +\eezz v(X\circ E) \\
&=& \nnn_X z\left(\nnn_\zdz v(e)+v(\frac{d-w}{2}e)\right) +
   \left(\nnn_\zdz v(X) + v((\Nu-\frac{w}{2}\id)X)\right) \nonumber \\
&=& -\nnn_\zdz v(X) + v(X) - v(\frac{d-w}{2}X) + \nnn_\zdz v(X) 
+ v((\Nu-\frac{w}{2}\id)X)\nonumber \\
&=& v((\Nu+\frac{2-d}{2}\id)X)\ .\nonumber
\end{eqnarray}
Now $\Lie_E(g) = (2-d)\cdot g$ follows from \eqref{5.20} and 
\eqref{5.13}.
\hfill $\qed$

\bigskip

The following two weaker cases are immediate consequences of the proof
of theorem \ref{t5.12}.

\begin{corollary}\label{t5.13}
Consider the two cases:

(a) a $(trTL)$-structure and the corresponding $(\nnn^rC)$-structure
on a vector bundle $K\to M$;

(b) a $(trTLP(w))$-structure and the corresponding $(\nnn^rCg)$-structure
on a vector bundle $K\to M$.

Suppose in both cases that $\zeta\in \OO(K)$ is a global section such
that $-C^K_\bullet \zeta:\tm\to \OO(K)$ is an isomorphism and 
$\nnn^r\zeta$=0.

In case (a) one obtains a Frobenius manifold without metric and Euler
field. In case (b) one obtains a Frobenius manifold without Euler field.
\end{corollary}

\begin{remarks}\label{t5.14}
(i) In the construction of Frobenius manifolds in singularity theory
\cite{SK3}\cite{SK4}\cite{SM2}, a partial Fourier transformation
maps the Gau{\ss}-Manin connection to a $(TERP(w))$-structure.
A primitive form of K. Saito yields an extension to a $(trTLEP(w))$-structure
and is mapped itself to the lift to the total bundle of a section
$\zeta$ as in theorem \ref{5.12}. Therefore one may call such 
sections also primitive forms. In a first structure connection
the unit field plays the role of a primitive form.

(ii) In the case of M-tame functions on affine manifolds one can by
work of Sabbah establish $(trTLEP(w))$-structures and
sections $\zeta\in\OO(K)$ with all properties in theorem \ref{5.12}
except for the property $\Nu^K\zeta =\frac{d}{2}\zeta$.
This {\it homogeneity property} is only clear for special cases
(Newton nondegenerate polynomials).
In the general case one obtains in the moment only Frobenius manifolds
without Euler field.
\end{remarks}

\subsection{$tt^*$ geometry on a Frobenius manifold}\label{s5.4}

Now we will show how the wishes after theorem \ref{t4.5}
can be satisfied. The construction in theorem \ref{t5.15} equips
the Frobenius manifold with a real structure $\kappa$ such that it
becomes a CDV-structure, if the $(trTLEP(w))$-structure carries a
suitable real structure. More precisely, one has the following.

\begin{theorem}\label{t5.15}
Let $(H\to \pmmm,\nnn,H_\R,P)$ be a tuple such that $(H,\nnn,P)$ is 
a $(trTLEP(w))$-structure (definition \ref{t5.5})
and $(H|_\cmmm,\nnn,H_\R,P)$ is a 
$(trTERP(w)$-structure (definition \ref{t2.15}). \\
By theorem \ref{t2.19} the bundle $K:=H|_\nmmm$ carries a CV-structure;
by theorem \ref{5.7} it carries a  Frobenius type structure
$(K,\nnn^K,g^k,\UU^K,\Nu^K)$.

Let $\zeta\in \OO(K)$ be a global section such that 
$v:=-C^K_\bullet \zeta:\tm\to\OO(K)$ is an isomorphism and such that
$\nnn^K\zeta =0$ and $\Nu^K\zeta=\frac{d}{2}\zeta$ for some $d\in \R$.

Then the CV-structure and the Frobenius type structure on $K$ are
shifted with $v^{-1}$ to a CDV-structure on the holomorphic tangent
bundle $TM$.
\end{theorem}

{\it Proof.}
The tensors $C^K,g^K,\UU^K$ are also part of the CV-structure on $K$.
Let $(TM,D,C,\www C,\kappa,h,g,\UU,\QQ)$ be the induced CV-structure
on $TM$ and let $(TM,\circ,e,E,g,\nnn^g,C,\UU,\Nu)$ be the induced
Frobenius type structure on $TM$ and Frobenius manifold structure on
$M$.

Let $\hat H\to \pmmm$ be the bundle constructed in lemma \ref{2.14} (e).
The bundles $H$ and $\hat H$ are different extensions of $H|_\cmmm$.
Both are fiberwise trivial with respect to 
$\pi:\pmmm \to M$. Let $\psi:\OO(K)\to \pi_*\OO(H)$ 
and $\hat \psi :\OO(K)\to \pi_*\OO(\hat H)$ be the canonical isomorphisms.
Then for $a\in \OO(K)$
\begin{eqnarray}\label{5.28}
\psi(a)\equiv \hat \psi(a) \mod z\chiii (H|_\cmmm)
\end{eqnarray}
and for $Y\in \tm$
\begin{eqnarray}\label{5.29}
\nnn \hat \psi v(Y) &=& \hat \psi v (DY)+\eezz\hat \psi v(CY)+
z \hat \psi v(\www C Y) \\
&+&(\eezz \hat\psi v(\UU Y)-\hat \psi v((\QQ-\frac{w}{2}\id)Y)
-z\hat\psi v(\kappa\UU\kappa Y))\dzz\ ,\nonumber 
\end{eqnarray}
and
\begin{eqnarray}
\nnn \psi v(Y) &=& \psi v (\nnn^g Y)+\eezz\psi v(CY)\label{5.30}\\
&+&(\eezz \psi v(\UU Y)- \psi v((\Nu-\frac{w}{2}\id)Y))\dzz\ .\nonumber
\end{eqnarray}
Especially
\begin{eqnarray}\label{5.31}
\nnn_{\zdz+E} \hat \psi v(Y) &=& \hat \psi v 
\left(D_EY-(\QQ-\frac{w}{2}\id)Y\right)
-z\hat\psi v(\kappa\UU\kappa Y)\ ,\\
\nnn_{\zdz+E} \psi v(Y) &=& \psi v \left([E,Y]+\frac{w+2-d}{2}Y\right)\ .
\label{5.32}
\end{eqnarray}
$\nnn_{\zdz+E}$ maps $z\chiii(H|_\cmmm)$ to itself. Therefore 
\eqref{5.28},\eqref{5.31}, and \eqref{5.32} show
\begin{eqnarray}\label{5.33}
\QQ=D_E-\Lie_E-\frac{2-d}{2}\id\ .
\end{eqnarray}
Write $\hat\psi v(e) = \psi v(e)+z\cdot a$ for some 
$a\in \chiii(H|_\cmmm)$. Then
\begin{eqnarray}\label{5.34}
\hat \psi v(D_ee) &=& \nnn_e\hat\psi v(e) +\eezz \hat\psi v(e) \\
&=& \nnn_e\psi v(e) +\eezz \psi v(e) +z \nnn_e a + a \nonumber\\
&=& \psi v (\nnn^g_ee) + z\nnn_e a +a = z\nnn_e a + a \nonumber \\
&\in & z\chiii (H|_\cmmm)\nonumber
\end{eqnarray}
because of $C_e=-\id$ and $\nnn^g_ee=0$.
This shows $D_ee=0$ and \eqref{4.13}-\eqref{4.16}. 
Finally, $\Lie_{E-\oooo E}(h)=0$ follows now from \eqref{5.33},
\eqref{4.21}, and \eqref{2.60}.
\hfill $\qed$

\clearpage

\section{Remarks on the semisimple case}\label{s6}
\setcounter{equation}{0}

\noindent
In \cite{CV1}\cite{CV2} the $tt^*$-equations turn up in the study of
moduli spaces of {\it massive} deformations of topological field
theories. The corresponding CV$\oplus$-structures are semisimple
at generic points of the moduli spaces.

A CV-structure on a bundle $K\to M$ is called {\it semisimple}
at $t\in M$, if the endomorphism $\UU|_t$ is semisimple with eigenvalues
with multiplicity one.

In \cite{CV1}\cite{CV2} the CV-structures are not semisimple at the 
(most interesting) points of conformal field theories;
at those points $\UU|_t=0$. A very intriguing moral in \cite{CV1}\cite{CV2}
is to study the CV-structures at semisimple points and get results
for the points with $\UU|_t=0$ and their conformal field theories.

Here we will just make some remarks about the semisimple points.

Consider a $(TERP(w))$-structure $(H\to \cmmm,\nnn, H_\R,P)$, that is,
(definition \ref{t2.12}) a holomorphic vector bundle
$H$ with a flat connection $\nnn$ on $H|_\csmmm$, which has a pole of
Poincar\'e rank 1 along $\nmmm$, a real $\nnn$-flat subbundle
$H_\R\subset H|_\csmmm$, and a pairing $P$ as in definition \ref{t2.12}.
In lemma \ref{t2.14} several objects are associated to it,
the bundle $K:=H|_\nmmm$ over $M$ with Higgs field $C$, holomorphic
metric $g$ and endomorphism $\UU$, and an extension $\hat H\to\pmmm$
of the bundle $H$ with a fiberwise $\C$-antilinear automorphism $\tau$,
which maps $\hat H_{z,t}$ to $\hat H_{(1/\oooo z,t)}$.

The $(TERP(w))$-structure is called semisimple at $t\in M$ if 
$\UU|_t$ is semisimple with eigenvalues of multiplicity one.

At semisimple points a connection with a pole of Poincar\'e rank 1 can be
described in terms of Stokes data, at least up to tensoring it
with a regular singular connection of rank 1.
We refer to \cite{Mal2}\cite{Mal3}\cite{Sab4} for definitions and 
discussions of this and just state some consequences.

Suppose the $(TERP(w))$-structure is semisimple everywhere.
Its sections take the form $z^{w/2}$ times sections of a connection
which is described by Stokes data. The Stokes data are given a priori
by two complex upper triangular $\rk H \times \rk H$ matrices $S$ and
$\www S$ with $\id$ on the diagonals. The monodromy of the 
$(TERP(w))$-structure is then given by $(-1)^wS\cdot {\www S}^{-1}$.
But the pairing $P$ implies $\www S = S^{tr}$, and the real subbundle
implies that $S$ is a real matrix. The Stokes matrix $S$ together with
the number $w$ and the eigenvalues of $\UU|_t$ for all $t\in M$ 
determines the $(TERP(w))$-structure completely \cite{Mal3}\cite{Sab4}.

In fact, for any matrix $S$ as above one obtains a universal semisimple
$(TERP(w))$-structure on 
\begin{eqnarray} \label{6.1}
 M_{univ}:=\mbox{ universal covering of }
\{(u_1,...,u_m)\in \C^m\ |\ u_i\neq u_j\, \forall i\neq j\}\ .
\end{eqnarray}
Here the eigenvalues of $\UU|_t$ are $u_1(t),...,u_m(t)$ for any 
$t\in M_{univ}$. Consider such a universal $(TERP(w))$-structure
$(H,\nnn,\H_\R,P)$. The functions $u_1,...,u_m$ on $M=M_{univ}$ are
local coordinates, the canonical coordinates.
$M$ is a semisimple F-manifold $(M,\circ,e)$ (section \ref{s4.1})
with $e_i:=\frac{\paa}{\paa u_i}, \ e=\sum_i e_i, 
\ e_i\circ e_j=\delta_{ij}e_i$.
The Higgs bundle $(K=H|_\nmmm,C)$ is a free $\tm$-module of rank 1.
There are unique global sections $v_i\in \OO(K), i=1,...,m$, with 
\begin{eqnarray} \label{6.2}
\OO_M\cdot v_i = \Im (C_{e_i}:\OO(K)\to \OO(K))
\end{eqnarray}
and $g(e_i,e_i)=1$. They satisfy
\begin{eqnarray} \label{6.3}
\bigoplus \OO_M\cdot v_i =\OO(K)\ ,\\
g(v_i,v_i)=\delta_{ij}\ .\label{6.4}
\end{eqnarray}
Now, being interested in CV-structures, one has to ask for which points
$t\in M$ the bundle $\hat H|_{\P^1\times \{t\}}$ is trivial

\begin{theorem}\label{t6.1}
Let $(H\to \C\times M_{univ},\nnn,H_\R,P)$ be a universal
semisimple $(TERP(w))$-structure as above.

(a) \cite[Proposition 2.2]{Du2} 
There exists a bound $b>0$ such that the restriction to
$\{t\in M_{univ}\ |\ |u_i-u_j|>b\ \forall\ i\neq j\}$ 
is a $(trTERP(w))$-structure.
The corresponding CV-structure there is a CV$\oplus$-structure, that means,
the hermitian metric $h$ is positive definite.

(b) \cite{CV1}\cite{CFIV}\cite[(4.5), (4.6)]{CV2} Let $\QQ_{mat}$ and
$\kappa_{mat}=h_{mat}$ be the matrices of the maps $\QQ$ and $\kappa$
and the metric $h$ with respect to the sections $v_i$. Then along a path
with $|u_i-u_j|\to \infty$ for $i\neq j$ one has $\QQ_{mat}\to 0$ and 
$\kappa_{mat}=h_{mat}\to \id$.
\end{theorem}

\begin{remarks}\label{t6.2}
(i) In \cite[ch. 2]{Du2} and \cite[4.2]{CV2} the coefficients of the lifts
of the $v_i$ to $\hat H$ with respect to special $\nnn$-flat bases
of $H|_{\{\Re z>0\}}$ and $H|_{\{\Re z<0\}}$ are written down and solved
in terms of integral equations.

(ii) The estimates for $\QQ_{mat}$ and $\kappa_{mat}=h_{mat}$ in 
\cite[(4.5), (4.6)]{CV2} are more precise. From the leading vanishing
terms one can recover the Stokes matrix $S$.

(iii) The subset of $M_{univ}$ where one does not have a CV-structure
is $R:=\{t\in M_{univ}\ |\ \hat H|_{\P^1\times \{t\} } \mbox{ is not trivial}
\}$. It is empty or a real analytic subvariety of $M_{univ}$, because
$\hat H$ is a real analytic family of bundles on $\P^1$
(lemma \ref{t2.18} (a)). In section 
\ref{s8.3} examples will be given where $R$ is not empty.
\end{remarks}

The monodromy of a $(TERP(w))$-structure with Stokes matrix $S$ is
$(-1)^wS(S^{-1})^{tr}$ in a suitable basis. The Stokes matrices
and $(TERP(w))$-structures which turn up in 
\cite{CV1}\cite{CV2} and in singularity theory have special properties:
the matrix $S$ has integer entries, and the monodromy is quasi\-uni\-po\-tent.

In \cite{CV2} a classification of such Stokes matrices was started.
For $\rk H=2$ and 3 it was carried out and corresponding field theories
were found and discussed.

The manifold $M_{univ}$ has the Euler field $E=\sum u_ie_i$.
The renormalization group flow in \cite{CV1}\cite{CV2} corresponds
on $M_{univ}$ to the flow of $E$ or of $E+\oooo E$. Recall that 
$D_{E-\oooo E}\QQ=0$ (formula \eqref{2.66b}). Of course, 
the moduli spaces of massive deformations of conformal field theories
contain non semisimple points.
The points $t_0$ with $\UU|_{t_0}=0$ are the fixed points of $E$ and
correspond to conformal field theories. The eigenvalues of
$\QQ|_{t_0}$ at such points are the {\it charges}, certain rational 
numbers $\alpha_j$ such that $(-1)^we^{-2\pi i\alpha_j}$ are the
eigenvalues of the monodromy. In the singularity case they are up to 
a shift the spectral numbers (= the exponents $-1$).

In \cite{CV2} an intriguing conjecture is made about these numbers.
It is a beautiful, although incomplete, example of how to use the 
geometry at semisimple points for points with $\UU|_{t_0}=0$.
Consider an $E$-orbit $M_{orb}$ of semisimple points whose closure
contains a point $t_0$ with $\UU|_{t_0}=0$. Suppose that there exist
a $\ciii$ one parameter family $S(r), r\in [0,1]$, of real Stokes matrices
with $S(1)=S$ and $S(0)=\id$, such that all monodromies
$(-1)^wS(r)(S(r)^{-1})^{tr}$ have all eigenvalues on the unit circle.
Suppose further that
for each $r\in [0,1]$ the $(TERP(w))$-structure on $M_{orb}$ extends to a
$(trTERP(w))$-structure on $\oooo {M_{orb}}$. Then for each $r$ the
eigenvalues of $\QQ(r)|_{t_0}$ are logarithms of the eigenvalues of 
$S(r)(S(r)^{-1})^{tr}$. For $r=0$ they are all 0. Now passing from 
$r=0$ to $r=1$ and observing how the eigenvalues of 
$S(r)(S(r)^{-1})^{tr}$ move around on the unit circle, one can 
determine the charges $\alpha_i$. Of course, these arguments are
quite incomplete. Especially neither existence nor essential uniqueness
of such a one parameter family of Stokes matrices is known.
Still the argument is very intriguing.

\clearpage

\section{Regular singular $(TERP)$-structures and Hodge filtrations}\label{s7}
\setcounter{equation}{0}

\noindent
In chapter \ref{s3} a correspondence between CV$\oplus$-structures with
$\UU=0$ and variations of polarized pure Hodge structures was presented.
In the case of hypersurface singularities $(TERP(w))$-structures and
CV-structures turn up such that at many points $t\in M$ the endomorphism
$\UU|_t$ is nilpotent, but not zero, the $(TERP(w))$-structure at 
$\C\times \{t\}$ is regular singular, but not logarithmic, and
such that to these points polarized mixed Hodge structures are 
associated.

The purpose of this chapter is to study such $(TERP(w))$-structures
and the arising filtrations. Sections \ref{s7.2}--\ref{s7.4} provide a
general frame and a systematic discussion. Section \ref{s7.5} describes
an important step in the construction of Frobenius manifolds in singularity
theory. Section \ref{s7.1} presents a characterization of polarized mixed
Hodge structures by nilpotent orbits which will be crucial in section 
\ref{s7.6}. 
There it will be used together with the results of sections 
\ref{s7.2}--\ref{s7.4} to study CV-structures along an Euler field.
The main result is theorem \ref{t7.20}. Together with the conjectural
inverse statement (remark \ref{t7.21}) it shows 
that polarized mixed Hodge structures turn up most naturally in CV-structures.
It will be used in the case of singularities (theorem \ref{t8.4} (b)).

\subsection{Polarized mixed Hodge structures and nilpotent orbits}\label{s7.1}

Throughout the whole section, $w$ will be an integer, $H$ a complex vector
space of finite dimension, $H_\R$ a real subspace with 
$H=H_\R\oplus iH_\R$, and $S$ a nondegenerate $(-1)^w$-symmetric pairing
on $H$ with real values on $H_\R$.

The definition below of a polarized mixed Hodge structure is motivated by 
Schmid's limit mixed Hodge structure \cite[Theorem 16.6]{Sch}.
It requires a nilpotent endomorphism $N$ of $H_\R$ and of $H$ which is an
infinitesimal isometry of $S$. From it one obtains a weight filtration
$W_\bullet$ in the following way.

\begin{lemma} \label{t7.1} \cite[Lemma 6.4]{Sch} 
Let $(H,H_\R,S,N,w)$ be as above.

(a) There exists a unique exhaustive increasing filtration $W_\bullet$ on
$H_\R$ such that $N(W_l)\subset W_{l-2}$ and such that 
$N^l:\Gr_{w+l}^W\to \Gr_{w-l}^W$ is an isomorphism.

(b) The filtration satisfies
$S(W_l,W_{l'})=0$ for $l+l'<w$.

(c) A nondegenerate $(-1)^{w+l}$-symmetric bilinear form $S_l$ is well 
defined on $\Gr_{w+l}^W$ for $l\geq 0$ by
$S_l(a,b):=S(\www a,N^l \www b)$ for $a,b\in \Gr_{w+l}^W$ 
with representatives $\www a,\www b\in W_{w+l}$.

(d) The primitive subspace $P_{w+l}\subset \Gr_{w+l}^W$ is defined by
\begin{eqnarray}  \label{7.1}
P_{w+l}:= \ker (N^{l+1}:\Gr_{w+l}^W\to \Gr_{w-l-2}^W)
\end{eqnarray}
for $l\geq 0$ and by $P_{w+l}:=0$ for $l<0$. Then
\begin{eqnarray} 
\Gr_{w+l}^W=\bigoplus_{i\geq 0}N^iP_{w+l+2i}\ ,
\end{eqnarray}
and this decomposition is orthogonal with respect to $S_l$ if $l\geq 0$.
\end{lemma}

\begin{definition}\label{t7.2} \cite{CaK1}\cite{He1}
A {\it polarized mixed Hodge structure of weight $w$} (abbreviation PMHS)
consists of data $H,H_\R,S,N$ and $W_\bullet$ as above and an exhaustive
decreasing {\it Hodge filtration} $F^\bullet$ on $H$ with the following
properties.

(a) $F^\bullet\Gr_k^W$ gives a pure Hodge structure of weight $k$,
that means, $\Gr_k^W = F^p \Gr_k^W\oplus \oooo{F^{k+1-p}\Gr_k^W}.$

(b) $N(F^p)\subset F^{p-1}$, i.e. $N$ is a $(-1,-1)$-morphism of mixed Hodge
structures.

(c) $S(F^p,F^{w+1-p})=0$.

(d) The pure Hodge structure $F^\bullet P_{w+l}$ of weight $w+l$ on 
$P_{w+l}$ is polarized by $S_l$, that means,
\begin{eqnarray} \label{7.3}
&& S_l(F^pP_{w+l},F^{w+l+1-p}P_{w+l})=0\ ,\\
&& i^{p-(w+l-p)}S_l(a,\oooo a)>0 \mbox{ \ for \ }a\in F^pP_{w+l}\cap
\oooo{F^{w+l-p}P_{w+l}}-\{0\}\ .\label{7.4}
\end{eqnarray}
\end{definition}

A PMHS comes equipped with a generalized Hodge decomposition, with
Deligne's $I^{p,q}$.

\begin{lemma}\label{t7.3} \cite[(1.2.8)]{De}\cite{CaK1}\cite[2.3]{He1}
For a PMHS as above define
\begin{eqnarray} 
I^{p,q} &:=& (F^p\cap W_{p+q})\cap (\oooo{F^q}\cap W_{p+q}+
\sum_{j>0}\oooo{F^{q-j}}\cap W_{p+q-j-1})\ ,\label{7.5}\\
I_0^{p,q} &:=& \ker (N^{p+q-w+1}:I^{p,q}\to I^{w-q-1,w-p-1})\ .\label{7.6}
\end{eqnarray}
Then
\begin{eqnarray} 
&& F^p=\bigoplus_{i,q:\ i\geq p}I^{i,q}\ ,\ \ 
W_l = \bigoplus_{p+q\leq l}I^{p,q}\ ,\label{7.7}\\
&& N(I^{p,q})\subset I^{p-1,q-1}\ ,\ \ I^{p,q}=
\bigoplus_{j\geq 0}N^jI_0^{p+j,q+j}\ ,\label{7.8}\\
&& S(I^{p,q},I^{r,s})=0 \mbox{ for \ }(r,s)\neq (w-p,w-q)\ ,\label{7.9}\\
&& S(N^iI_0^{p,q},N^jI_0^{r,s})=0 \mbox{ for \ } 
(r,s,i+j)\neq (q,p,p+q-w)\ ,\label{7.10}\\
&& I^{q,p}\cong \oooo {I^{p,q}} \mod W_{p+q-2}\ ,\label{7.11}\\
&& I_0^{q,p}\cong \oooo {I_0^{p,q}} \mod W_{p+q-2}\ .\label{7.12}
\end{eqnarray}
\end{lemma}

A polarized Hodge structure is a polarized mixed Hodge structure with 
$N=0$. Let us fix one reference polarized Hodge structure
$(H,H_\R,S,F^\bullet_0)$ of weight $w$. The space
\begin{eqnarray} 
\check{D} := \{\mbox{filtrations }F^\bullet\subset H\ |\ 
\dim F^p=\dim F^p_0,\ S(F^p,F^{w+1-p})=0\}\label{7.13}
\end{eqnarray}
is a complex homogeneous space and a projective manifold, and the subspace
\begin{eqnarray} 
&& D:= \{ F^\bullet\in \check{D}\ |\ F^\bullet\mbox{ is part of a polarized
Hodge structure}\}\label{7.14}
\end{eqnarray}
is an open submanifold and a real homogeneous space \cite{Sch}.
It is a classifying space for polarized Hodge structures with fixed
Hodge numbers.

\begin{definition}\label{t7.4}
A pair $(F^\bullet,N)$ is said {\it to give rise to a nilpotent orbit}
if $F^\bullet \in \check{D}$, the endomorphism $N$ of $H_\R$ is nilpotent
and an infinitesimal isometry with $N(F^p)\subset F^{p-1}$, and there
exists a bound $b\in \R$ such that 
\begin{eqnarray} \label{7.15}
e^{zN}F^\bullet \in D \mbox{ for \ } \Im z>b\ .
\end{eqnarray}
Then the set $\{e^{zN}F^\bullet\ |\ z\in \C\}$ is called {\it a nilpotent
orbit}.
\end{definition}

Nilpotent orbits play a fundamental role in Schmid's work \cite{Sch}.
They give rise to and in fact characterize PMHS's. This will be crucial
for the applications to CV-structures in section \ref{s7.6}.

\begin{theorem}\label{t7.5} \cite{Sch}\cite{CaK1}\cite{CaKS}\cite{CaK2}
Let $(H,H_\R,S)$ be as above.
\begin{list}{}{}
\item[(a)] 
The tuple $(H,H_\R,S,F^\bullet,N)$ is a PMHS of weight $w$ \\
\qquad $\iff$ the pair $(F^\bullet,N)$ gives rise to a nilpotent orbit.
\item[(b)]
If $(H,H_\R,S,F^\bullet,N)$ is a PMHS with $I^{q,p}=\oooo{I^{p,q}}$
then $e^{zN}F^\bullet \in D$ for $\Im z>0$.
\end{list}
\end{theorem}

{\it Proof.}
`$\Longleftarrow$' is \cite[Theorem 16.6.]{Sch}. It is a consequence of the
$SL_2$-orbit theorem. 
`$\Longrightarrow$' is \cite[Corollary 3.13]{CaKS}.
Short proofs of both directions are given in \cite[Theorem 3.13]{CaK2}.
The special case $(b)$ is proved in \cite[Proposition 2.18]{CaK1}
and in \cite[Lemma 3.12]{CaKS}.

Theorem \ref{t7.20} will generalize `$\Longrightarrow$'. As a preparation
for its proof, we will now sketch a rather elementary proof of  
`$\Longrightarrow$', which uses (b).

The first and longest part of the proof is devoted to showing that
$e^{zN}F^\bullet$ gives a pure Hodge structure for $z\in \C$ with
$|\Im z|$ sufficiently large, forgetting the polarization.
Equivalently, we want to show
\begin{eqnarray} \label{7.16}
e^{(z-\oooo z)N}F^p \oplus \oooo{F^{w+1-p}}=H \mbox{ for }|z-\oooo z|
\mbox{ large.}
\end{eqnarray}
We choose a Jordan block basis of $H$ with respect to $N$, which respects
the $I_0^{p,q}$ and the real structure as good as possible.
Namely, we choose a basis $\{ v_a\}_{a\in A}$ of 
$\bigoplus_{p,q}I_0^{p,q}$ such that 
\begin{eqnarray} \label{7.17}
v_a\in I_0^{p(a),q(a)}\ ,\\
\oooo{v_a}\equiv v_{k(a)}\mod W_{p(a)+q(a)-2}\label{7.18}\\
N^{n(a)+1}v_a=0\ ,\ N^{n(a)}v_a\neq 0\ .\label{7.19}
\end{eqnarray}
Then 
\begin{eqnarray} 
p(a)+q(a)-w=n(a)\ ,\ k(k(a))=a\ ,\ n(k(a))=n(a)\ ,\label{7.20}\\
p(k(a))=q(a) ,\ q(k(a))=p(a)\ .\label{7.21}
\end{eqnarray}
The set
\begin{eqnarray} \label{7.22}
\bigcup_{a\in A}\{ N^jv_a\ |\ j=0,1,..,n(a)\}
\end{eqnarray}
is a basis of $H$, the set
\begin{eqnarray} \label{7.23}
\bigcup_{a\in A}\{e^{(z-\oooo z)N}N^jv_a\ |\ j=0,1,..,p(a)-p\}
\end{eqnarray}
is a basis of $e^{(z-\oooo z)N}F^p$, and the set
\begin{eqnarray} \label{7.24}
\bigcup_{a\in A}\{N^j\oooo{v_{k(a)}}\ |\ j=0,1,...,p(k(a))-(w+1-p)\}
\end{eqnarray}
is a basis of $\oooo{F^{w+1-p}}$.

Consider the coefficient matrix $B$ with which the bases in \eqref{7.23}
and \eqref{7.24} are written in terms of that in \eqref{7.22}.
We claim
\begin{eqnarray} \label{7.25}
\det B = c_0\cdot (z-\oooo z)^{c_1} +\mbox{ smaller powers of }(z-\oooo z)
\end{eqnarray}
where $c_0\neq 0$ depends only on the numbers $p,p(a),q(a),n(a)$ ($a\in A$)
and
\begin{eqnarray} \label{7.26}
c_1=\sum_{a:\ p(a)\geq p} c_1(a) 
\end{eqnarray}
with $c_1(a):= (p(a)-p+1)(n(a)-p(a)+p)$ if $0\leq p(a)-p\leq n(a)$ and
$c_1(a):=0$ else.

First consider the case of a PMHS with $I^{q,p}=\oooo{I^{p,q}}$. 
Then $\oooo{v_a} = v_{k(a)}$. The matrix $B$ consists of blocks in the
diagonal, one for each $a\in A$, of size $(n(a)+1)\times(n(a)+1)$.
One easily sees that each diagonal block contributes a monomial
$c_0(a)(z-\oooo z)^{c_1(a)}$ to $\det B$. Because of 
\cite[Proposition 2.18]{CaK1} (or by elementary calculation) 
$c_0=\prod_{a\in A}c_0(a)\neq 0$.

More precisely, for a diagonal block of a fixed $a\in A$, the subsquare
of coefficients of the part of $a$ in \eqref{7.23} with respect to the
largest powers of $N$ in $\{ N^jv_a\ |\ j=0,1,..,n(a)+1\}$ contributes
$c_0(a)(z-\oooo z)^{c_1(a)}$, the dual subsquare in the diagonal block
is the identity matrix and contributes 1.

Going over to the general case with $I^{q,p}\neq \oooo{I^{p,q}}$, 
only the coefficients in $B$ which express \eqref{7.24} in terms of 
\eqref{7.22} change. The leading monomial 
$\prod_{a\in A}c_0(a)(z-\oooo z)^{c_1(a)}$ is contributed by the same
subsquares of the diagonal blocks. But now a priori more pieces than
only the dual subsquares are involved in the coefficients of this
monomial.

Nevertheless, by an induction on the size of the Jordan blocks one sees that
still the only contribution comes from the dual subsquares and is 1.
The details are left to the reader.

Now the claim \eqref{7.25} is proved. It shows \eqref{7.16}. It rests
to see $e^{zN}F^\bullet \in D$ for large $\Im z$. One can choose a 
continuous one parameter family of PMHS's with Hodge filtrations
$F^\bullet(r),$ $r\in [0,1]$, such that $F^\bullet(1)=F^\bullet$ and
$F^\bullet(0)$ satisfies $I^{q,p}(0)=\oooo{I^{p,q}(0)}$. For 
$z\in \C$ with $\Im z$ sufficiently large $\det B\neq 0$ for all
$r\in [0,1]$, and $e^{zN}F^\bullet(r)$ is a Hodge filtration and has
a nondegenerate hermitian form $h(r)$ as in \eqref{3.8}+\eqref{3.9}.
By \cite[Proposition 2.18]{CaK1} $e^{zN}F^\bullet(0)\in D$ and 
$h(0)$ is positive definite.
Therefore all $h(r)$ are positive definite and all
$e^{zN}F^\bullet(r)\in D$.
\hfill $\qed$

\subsection{Structures around connections on $\C^*$}\label{s7.2}

The following definitions are modelled after the singularity case.
In section \ref{s8.1} the objects will be identified with topological
data of a singularity. Section \ref{s7.3} will discuss objects which
correspond to nontopological  data of a singularity.

Throughout this section we will fix the following data
$(H,\nnn,H_\R,P,w)$ and develop some useful structures around them:

(a) A flat complex vector bundle $H\to \C^*$ of rank $\mu\in \N$ 
with connection $\nnn$
and monodromy $M_{mon}=M_s\cdot M_u$ with semisimple
part $M_s$, unipotent part $M_u$ and nilpotent part $N:=\log M_u$
on each fiber $H_z$, $z\in \C^*$, such that the eigenvalues of 
$M_{mon}$ have absolute value 1.

(b) A real $\nnn$-flat subbundle $H_\R\to \C^*$ with 
$H_z=H_{\R,z}\oplus iH_{\R,z}$ for $z\in \C^*$.

(c) An integer $w$ and a $\nnn$-flat nondegenerate bilinear pairing
\begin{eqnarray} \label{7.27}
P:H_z\times H_{-z}\to \C\mbox{ for \ }z\in \C^*
\end{eqnarray}
which takes values in $i^w\R$ on $H_\R$.

First we derive from this structure a tuple
$(H^\infty,H^\infty_\R,M_s,N,S)$ of `topological' data which together
with certain filtrations may give PMHS's.

(d) $H^\infty$ is the space of manyvalued global $\nnn$-flat sections
in $H$. More precisely, define
\begin{eqnarray} \label{7.28}
{\bf e}:\C\to \C^*  ,\ \zeta\mapsto e^{2\pi i\zeta}=z\ ,
\end{eqnarray}
and let $pr_H:{\bf e}^*H\to H$ be the canonical projection. Then
\begin{eqnarray} \label{7.29}
&& H^\infty =\{ A:\C\to H\ |\ A=pr_H\circ \www A,\ \www A \mbox{ is a 
flat global section of }{\bf e}^*H\}.
\end{eqnarray}

(e) The space $\hiii$
is equipped with the real subspace $H^\infty_\R$ of manyvalued
$\nnn$-flat sections in $H_\R$ and with a monodromy which is also
denoted $M_{mon}$, with semisimple part $M_s$ and nilpotent part
$N$. In fact, all monodromy invariant structure on a fiber $H_z$ can be
shifted to $H^\infty$. Define
\begin{eqnarray} \label{7.30}
H^\infty_\lambda &:=& \ker (M_s-\lambda \id:H^\infty\to H^\infty)
\mbox{ for \ }\lambda\in S^1 \ ,\\
H^\infty_{\neq 1} 
&:=& \bigoplus_{\lambda\neq 1}H^\infty_\lambda \ .\label{7.31}
\end{eqnarray}
The decomposition $H^\infty = H^\infty_1\oplus H^\infty_{\neq 1}$ 
is already defined for $H_\R$.

(f) A {\it polarizing form} $S$ on $H^\infty$ is defined in the following
way.
It may look artificial, but it is precisely the 
right way (cf. remark \ref{t7.8}; for the singularity case see
section \ref{s8.1}).
Let $\gamma_\pi:H_z\to H_{-z}$ be the $\nnn$-flat shift in 
mathematical positive direction; so $\gamma_\pi\circ \gamma_\pi=M_{mon}$.
Define
\begin{eqnarray} \label{7.32}
L:H_z\times H_z\to \C,\ (a,b)\mapsto P(a,\gamma_\pi(b))\ .
\end{eqnarray}
Then
\begin{eqnarray} \nonumber
\qquad L(a,b) &=& P(a,\gamma_\pi(b))= 
P(\gamma_\pi(a),\gamma_\pi\circ\gamma_\pi(b))
= (-1)^w P(M_{mon}b,\gamma_\pi (a)) \\
&=& (-1)^w L(M_{mon}b,a) \ ,\label{7.33}\\
\qquad L(a,b)&=& L(M_{mon}a,M_{mon}b) \ .\label{7.34}
\end{eqnarray}
Then $L$ is monodromy invariant and can be shifted to $\hiii$.
Define $S:H^\infty\times H^\infty\to \C$ by
\begin{eqnarray} \label{7.35}
S(a,b) &:=& (-1)(2\pi i)^wL(a,\frac{1}{M_{mon}-\id}b) \mbox{ for }
a\in \hiii,b\in \hiin\ ,\\
S(a,b) &:=& (-1)(2\pi i)^wL(a,\frac{-N}{M_{mon}-\id}b) \mbox{ for }
a\in \hiii,b\in \hiie\ ,\label{7.36} 
\end{eqnarray}
where $\frac{-N}{M_{mon}-\id}:= -(\sum_{k\geq 1}\frac{1}{k!}N^{k-1})^{-1}$
on $\hiie$ is an automorphism of $\hiie$.

\begin{lemma}\label{t7.6}
The pairing $S$ is monodromy invariant and nondegenerate.
The restriction to $\hiin$ is $(-1)^{w-1}$-symmetric, the restriction
to $\hiie$ is $(-1)^{w}$-symmetric. It takes real values on $\hiii_\R$.
\end{lemma}

{\it Proof.}
Only the symmetry properties are not completely trivial.
For this proof denote $c:=(-1)(2\pi i)^w$.
For $a,b\in H^\infty_{\neq 1}$ and $\www a = (M_{mon}-\id)a$, $
\www b=(M_{mon}-\id)b$
\begin{eqnarray} \nonumber
S(\www b,\www a) &=& c \cdot L((M_{mon}-\id)b,a) \\
&=& c \cdot (-1)^w L(a,b)-c \cdot (-1)^w L(M_{mon}a,b)\nonumber\\
&=& c \cdot (-1)^{w-1}L((M_{mon}-\id)a,b) \nonumber\\
&=& (-1)^{w-1}S(\www a,\www b)\ .\label{7.37}
\end{eqnarray}
\eqref{7.34} implies $L(Na,b)=L(a,-Nb)$. This is used in the next calculation.
For $a,b\in \hiie$
\begin{eqnarray} \nonumber
S(b,a) &=& c \cdot L(b,\frac{-N}{e^N-\id}a) 
= c \cdot L(\frac{N}{e^{-N}-\id}b,a) \\
&=& c \cdot L(M_{mon}\frac{-N}{e^N-\id}b,a) 
=c \cdot (-1)^wL(a,\frac{-N}{e^N-\id}b)\nonumber\\
&=& (-1)^w S(a,b) \ .\label{7.38}
\end{eqnarray}
\hfill $\qed$

\bigskip

In section \ref{s7.3} filtrations in $\hiii$ will be related to additional
structure on $H\to \C^*$. For this we need the following 
sections $es(A,\alpha)$ and spaces $C^\alpha$ and $V^\alpha$.

Fix $A\in \hiil$ and $\alpha\in \R$ with $e^{-2\pi i\alpha}=\lambda$.
Define a global holomorphic section $es(A,\alpha)$ of $H\to \C^*$ by
\begin{eqnarray} 
es(A,\alpha) (z) &:=& e^{2\pi i\alpha\zeta}\exp (-\zeta N)A(\zeta)
\mbox{ for } \zeta \mbox{ with }e^{2\pi i\zeta}=z\label{7.39}\\
&=& \mbox{``}z^\alpha\exp(\log z\cdot \frac{-N}{2\pi i})A 
    \mbox{''} \ .\label{7.40}
\end{eqnarray}
One sees that the right hand side of \eqref{7.39} is independent of the 
choice of $\zeta$ with $e^{2\pi i\zeta}=z$. The section $es(A,\alpha)$
is called {\it an elementary section of order }$\alpha$. 
Let $C^\alpha$ be the space of all elementary sections of order $\alpha$.
The map
\begin{eqnarray} \label{7.41}
\psi_\alpha : \hiil \to C^\alpha\ ,\ A\mapsto es(A,\alpha)
\end{eqnarray}
is an isomorphism of vector spaces. By definition one has
\begin{eqnarray}\nonumber
z\cdot es(A,\alpha) &=& es(A,\alpha + 1)\ , 
\ z\circ \psi_\alpha = \psi_{\alpha +1}\ ,\\
\ndz es(A,\alpha) &=& \alpha\cdot es(A,\alpha -1) 
     + es(\frac{-N}{ 2\pi i}A, \alpha -1)\ , \nonumber \\
(\nzdz -\alpha) es(A,\alpha) &=& es (\frac{-N}{ 2\pi i}A,\alpha) 
  = \frac{-N}{2\pi i} es(A,\alpha)\ ,\label{7.42} \\
z&:& C^\alpha \to C^{\alpha +1} \qquad \mbox{ bijective, }\nonumber \\
\ndz &:& C^\alpha \to C^{\alpha- 1} \qquad \mbox{bijective iff }
     \alpha\neq 0\ , \nonumber \\
z\ndz -\alpha = \frac{-N}{2\pi i} &:& C^\alpha \to C^\alpha 
     \qquad \mbox{ nilpotent.} \label{7.43}
\end{eqnarray}
Let $i:\C^*\to \C$ be the inclusion. The space
\begin{eqnarray} \label{7.44} 
V^{>-\infty} := \sum_\alpha \C\{t\}C^\alpha \subset
(i_*\OO(H))_0
\end{eqnarray}
is the space of all germs at 0 of sections in $\OO(H)$ of moderate growth.
It is a $\Cz[z^{-1}]$-vector space with dimension equal to $\rk H$.
The Kashiwara--Malgrange filtration on $V^{>-\infty}$ is a decreasing
filtration, indexed by the set $\{\alpha\in \R\ |\ e^{-2\pi i \alpha}$
is an eigenvalue of $M_{mon}\}$, with subspaces
\begin{eqnarray} \label{7.45}
V^\alpha := \sum_{\beta\geq \alpha} \Cz C^\beta 
=\bigoplus_{\alpha\leq \beta<\alpha+1}\Cz C^\beta
\end{eqnarray}
and additionally
\begin{eqnarray} \label{7.46}
V^{>\alpha} = \sum_{\beta >\alpha}\Cz C^\beta 
=\bigoplus_{\alpha < \beta \leq \alpha+1}\Cz C^\beta\ .
\end{eqnarray}
Then $\Gr_V^\alpha := V^\alpha/V^{>\alpha} \cong C^\alpha$.

\bigskip
Finally, we will need a Fourier--Laplace transformation on elementary sections.
The following automorphisms $G^{(\alpha)}$ of $\hiia$ for $\alpha>0$ are
central in the next proposition,
\begin{eqnarray} \label{7.47}
G^{(\alpha)} &:=& \sum_{k\geq 0}\frac{1}{k!}\Gamma^{(k)}(\alpha) 
\left( \frac{-N}{2\pi i}\right)^k : \hiia\to \hiia \\
&=& \mbox{``} \Gamma (\alpha\cdot \id + \frac{-N}{2\pi i})
\mbox{''}\ .\label{7.48}
\end{eqnarray}
Here $\Gamma^{(k)}$ is the $k$-th derivative of the gamma function.

\begin{proposition}\label{t7.7}
(a) Let $\tau$ and $z$ both be coordinates on $\C$. Fiz $\alpha>0$ and
$A\in \hiia$. The Fourier--Laplace transformation $FL$ with
\begin{eqnarray} \label{7.49}
FL(es(A,\alpha-1)(\tau))(z) := \int_0^{\infty\cdot z}
e^{-\tau/z}\cdot es(A,\alpha-1)(\tau)\ddd\tau
\end{eqnarray}
is well defined and gives the elementary section
\begin{eqnarray} \label{7.50}
FL(es(A,\alpha-1)(\tau))(z) = es(G^{(\alpha)}A,\alpha)(z)\ .
\end{eqnarray}

(b) For $\alpha,\beta$ with $0<\alpha,\beta <1$ and for 
$A\in \hiia, B\in \hiib$
\begin{eqnarray} \label{7.51}
&& P\left( es(G^{(\alpha)}A,\alpha)(z),\, es(G^{(\beta)}B,\beta)(-z)\right)
=z\frac{1}{(2\pi i)^{w-1}}\cdot S(A,B)\ .
\end{eqnarray}
For $\alpha=\beta=1$ and for 
$A,B\in \hiie$
\begin{eqnarray} \label{7.52}
&& P\left( es(G^{(1)}A,1)(z),\, es(G^{(1)}B,1)(-z)\right)
=z^2\frac{-1}{(2\pi i)^{w}}\cdot S(A,B)\ .
\end{eqnarray}
\end{proposition}

{\it Proof.}
We omit the details. They are straightforward computations.
Part (a) uses the identity
\begin{eqnarray} \label{7.53}
\left(\frac{\ddd }{\ddd\alpha}\right)^k (\Gamma(\alpha+1)z^{\alpha+1})
= \int_0^\infty e^{-\tau/z} \tau^\alpha(\log \tau)^k\ddd \tau\ ,
\end{eqnarray}
part (b) uses
\begin{eqnarray} \label{7.54}
\Gamma(\frac{1}{2}+x)\Gamma(\frac{1}{2}-x)=\frac{\pi}{\cos \pi x}
\end{eqnarray}
and (which is a consequence of \eqref{7.54})
\begin{eqnarray} 
\Gamma(1+x)\Gamma(1-x)=\frac{\pi x}{\sin \pi x}\ .
\end{eqnarray}
\hfill $\qed$

\begin{remark}\label{t7.8}
The formulas \eqref{7.51} and \eqref{7.52} show that the definition
above of the polarizing form $S$ is a good one and not as artificial
as it might have looked.
\end{remark}

\subsection{Filtrations from regular singular $(TERP)$-structures}\label{s7.3}

We stick to the situation considered in section \ref{s7.2}, with all the
objects $H, \mu, \nnn, H_\R, P, w, \hiii, \hiii_\R, M_s, N, S, 
es(A,\alpha), C^\alpha, V^\alpha, \vgii, G^{(\alpha)}$.
The purpose of this section is to describe how an extension of the
bundle $H\to \C^*$ over $\{0\}$ may give rise to a filtration $F^\bullet$
on $\hiii$ which together with $(\hiii,\hiii_\R,M_s,-N,S)$ might give a 
sum of two PMHS's of weight $w$ and $w-1$ with an automorphism $M_s$.

In the following, $\www H$ will always mean a vector bundle on $\C$ which 
extends $H$, that means, $\www H|_{\C^*}=H$ and which satisfies the 
following two assumptions.
\begin{list}{}{}
\item[$(\alpha)$] $(\www H,\nnn, H_\R,P)$ is a $(TERP(w))$-structure.
\item[$(\beta)$]    The pair $(\www H,\nnn)$ has a regular singularity at 0,
that means, $\OO(\www H)_0\subset \vgii$, in other words, the 
holomorphic sections of $\www H$ have moderate growth at 0.
\end{list}

A $(TERP(w))$-structure on $M=\{pt\}$ with $(\beta)$ is called
a {\it regular singular} $(TERP(w))$-structure. 
Assumption $(\beta)$ implies that the endomorphism 
$\UU:\www H_0 \to \www H_0$ is  nilpotent (see e.g. \cite[II.4.1]{Sab4}).

\begin{remark} \label{t7.9}
In the case of a $(TERP(w))$-structure $(\www H,\nnn,H_\R,P)$
without assumption $(\beta)$, Sabbah \cite{Sab2} proposed a way to define
a filtration on $\hiii$ from $\www H$. In the case of tame functions on 
affine manifolds this filtration yields a mixed Hodge structure 
\cite{Sab1}\cite{Sab2}.
The way is similar, but also dual to the way described below. We refer to
\cite{Sab2} for it.
\end{remark}

The following way to define a filtration from a pair $(\www H,\nnn)$ with
regular singularity is essentially already used in \cite{Sch} and \cite{Va1}.
Consider the isomorphisms
$\psi_\alpha:\hiil\to C^\alpha$ from \eqref{7.41} 
and $C^\alpha\cong \Gr_V^\alpha$. Define two $M_s$-invariant decreasing
filtrations $F^\bullet$ and $\www F^\bullet$ on $\hiii$ by
\begin{eqnarray}  \label{7.56a}
F^p\hiil &:=& \psi_\alpha^{-1} z^{-w+1+p}\Gr_V^{\alpha+w-1-p}
\OO(H)_0\ ,\\
\www F^p \hiil &:=& (G^{(\alpha)})^{-1} F^p \hiil \label{7.56b}
\end{eqnarray}
for $\alpha\in (0,1]$ with $e^{-2\pi i\alpha}=\lambda$.

\begin{lemma}\label{t7.10}
The filtrations satisfy
\begin{eqnarray}  \label{7.57}
\dim F^p\hiil = \dim \www F^p\hiil\ ,\\
\dim \www F^p\hiil + \dim \www F^{w-p}\hiii_{\oooo \lambda} =
\dim \hiil \mbox{ \ for \ } \lambda \neq 1\ ,\label{7.58}\\
\dim \www F^p\hiie + \dim \www F^{w+1-p}\hiie =
\dim \hiie \ ,\label{7.59}\\
N(F^p)\subset F^{p-1}\ ,\ N(\www F^p)\subset \www F^{p-1}\ ,\label{7.60}\\
S(\www F^p\hiil,\, \www F^{w-p}\hiii_{\oooo\lambda})=0 
\mbox{ \ for \ }\lambda\neq 1\ ,\label{7.61}\\
S(\www F^p\hiie,\, \www F^{w+1-p}\hiie)=0 \ . \label{7.62}
\end{eqnarray}
\end{lemma}

{\it Proof.}
The definition of $F^\bullet$ can be rewritten as follows. Consider
any $\beta\in \R$.
\begin{eqnarray} \label{7.63}
A\in F^{w+[-\beta]}\hiib \iff es(A,\beta)\in \Gr_V^\beta \OO(\www H)_0\ .
\end{eqnarray}
Now $N(F^p)\subset F^{p-1}$ follows from the pole of Poincar\'e rank 1,
\begin{eqnarray} 
es(\frac{-N}{2\pi i}A,\beta +1) &=& 
(\nnn_\zdz -(\beta+1))es(A,\beta+1) \nonumber \\
&\in & (\nnn_\zdz -(\beta+1))z\Gr_V^\beta \OO(\www H)_0 \nonumber \\
&=& \Gr_V^{\beta+1} \left( (\nnn_\zdz -(\beta+1)) 
z\OO(\www H)_0\right)\nonumber\\
&\subset & \Gr_V^{\beta+1}\OO(\www H)_0\ .
\end{eqnarray}
The maps $G^{(\alpha)}$ are automorphisms with $NG^{(\alpha)}=G^{(\alpha)}N$.
This shows \eqref{7.60} and \eqref{7.57}. The pairing $P$ satisfies by 
definition
\begin{eqnarray} \label{7.65}
&& P: C^\alpha \times C^\beta \to 0 \mbox{ \ if } \alpha+\beta\notin \Z\ ,\\
&& P: C^\alpha \times C^\beta \to z^{\alpha+\beta}\C 
\mbox{ \ nondegenerate if }
\alpha+\beta\in \Z\ .\label{7.66}
\end{eqnarray}
Because of the $(TERP(w))$-structure one has for $b\in \vgii$
\begin{eqnarray} \label{7.67}
P(\OO(\www H)_0,b)\subset z^w\OO_{\C,0} \iff b\in \OO(\www H)_0\ .
\end{eqnarray}
This implies for $\beta \in\R$ and $b\in C^\beta$
\begin{eqnarray} \label{7.68}
P(\Gr_V^{w-1-\beta} \OO(\www H)_0,b)=0
\iff b\in \Gr_V^\beta \OO(\www H)_0\ .
\end{eqnarray}
With \eqref{7.63} and \eqref{7.51} this shows for $\lambda\neq 1$
and $B\in \hiii_{\oooo \lambda}$
\begin{eqnarray} \label{7.69}
S(\www F^p\hiil,B)=0\iff B\in \www F^{w-p}\hiii_{\oooo \lambda}\ .
\end{eqnarray}
This gives \eqref{7.61} and \eqref{7.58}, because $S$ is nondegenerate.
\eqref{7.62} and \eqref{7.59} follow analogously with \eqref{7.52}.
\hfill $\qed$

\bigskip

The lemma shows that the filtrations $F^\bullet$ and $\www F^\bullet$ have some
properties of Hodge filtrations of (mixed) Hodge structures of weight
$w$ on $\hiie$ and weight $w-1$ on $\hiin$. In fact, in the singularity case
$\www F^\bullet$ is the Hodge filtration of such PMHS's
(\cite[Theorem 10.29]{He2}, \cite[Theorem 3.5]{He1}, and section \ref{s8.1}
step 4 $(\gamma)$). 

In general the given germ $\OO(\www H)_0$ is richer than the filtered parts
$\Gr_V^\alpha\OO(\www H)_0$, which give the filtrations.
The proof of the lemma shows that the pole of Poincar\'e rank 1 and the
relation \eqref{7.67} are natural generalizations of the properties 
\eqref{7.58}--\eqref{7.62}. So one can consider
regular singular $(TERP(w))$-structures $(\www H\to \C,\nnn,H_\R,P)$ 
such that $\www F^\bullet$ is a Hodge filtration of PMHS's as a natural 
generalization of PMHS's.

This point of view was taken in \cite{He1}; but there the Fourier--Laplace
duals of $(TERP(w))$-structures with $\OO(\www H)_0\subset V^{>0}$ were 
considered and called Brieskorn lattices.

A regular singular $(TERP(w))$-structure has important discrete invariants,
its exponents. Their definition is modelled after the singularity case
(\cite{Va1}, cf. also \cite{He1}\cite{He2}). The {\it exponents}
are real numbers $\alpha_1,...,\alpha_\mu$ such that
\begin{eqnarray} \label{7.70}
\sharp (i\ |\ \alpha=\alpha_i) =\dim \Gr_V^\alpha \OO(\www H)_0 - 
\dim \Gr_V^{\alpha-1} \OO(\www H)_0\ .
\end{eqnarray}
They are ordered by $\alpha_1\leq ... \leq \alpha_\mu$.
From \eqref{7.68} one derives the symmetry
\begin{eqnarray} \label{7.71}
\alpha_i+\alpha_{\mu+1-i} =w\ .
\end{eqnarray}
One also sees easily
\begin{eqnarray} \label{7.72}
&& V^{\alpha_1}\supset \OO(\www H)_0 \supset V^{>\alpha_\mu -1}\ ,\\
&& \OO(\www H)_0 = \left( \OO(\www H)_0 \cap 
\sum_{\alpha_1\leq \beta \leq \alpha_\mu -1}C^\beta \right) \oplus 
V^{>\alpha_\mu-1}\ .\label{7.73}
\end{eqnarray}

\subsection{$(TERP)$-structures generated by elementary sections}\label{s7.4}

We stick to the situation considered in sections \ref{s7.2} and \ref{s7.3}.
Especially we fix $(H,\nnn, H_\R,P)$ and $(\hiii,M_s,N,S)$.
If $(\www H,\nnn, H_\R,P)$ is a $(TERP(w))$-structure with regular 
singularity at 0 then in general one cannot recover from the filtration
$F^\bullet$ the $\Cz$-lattice $\OO(\www H)_0$ and the $(TERP(w))$-structure.
But one can recover the quotients $\Gr_V^\alpha \OO(\www H)_0$ and define
with them a slightly simpler $(TERP(w))$-structure. One obtains a 
correspondence (corollary \ref{t7.11}) and can describe precisely when such
a $(TERP(w))$-structure is a $(trTERP(w))$-structure (lemma \ref{t7.12}).

Let $F^\bullet$ and $\www F^\bullet$ be decreasing $M_s$-invariant filtrations
on $\hiii$ with \eqref{7.56b} and with the properties in lemma \eqref{7.10}.
Let $\www H\to \C$ be the extension of $H\to \C^*$ such that the sheaf 
$\OO(\www H)$ is generated by the elementary sections in 
\begin{eqnarray} \label{7.74a}
\{es(A,\beta)\ |\ \beta\in \R, \ A\in F^{w+[-\beta]}\cap \hiib\}\ .
\end{eqnarray}
Then the proof of lemma \ref{7.10} and especially formula \eqref{7.63}
show that $(\www H,\nnn, H_\R,P)$ is a $(TERP(w))$-structure with a
regular singularity at 0 and that $F^\bullet$ and $\www F^\bullet$ are the
corresponding filtrations, i.e. \eqref{7.56a} and \eqref{7.56b} hold.
If $\alpha_1,...,\alpha_\mu$ are the exponents of this $(TERP(w))$-structure
then
\begin{eqnarray} \label{7.74b}
\OO(\www H)_0 = \bigoplus_{\alpha_1\leq\beta \leq\alpha_\mu -1}
\Gr_V^\beta \OO(\www H)_0 \oplus V^{>\alpha_\mu -1}\ .
\end{eqnarray}
The proof of lemma \ref{7.10} yields the following correspondence.

\begin{corollary}\label{t7.11}
The formulas \eqref{7.56a} and \eqref{7.74a} give a one-to-one correspondence
between the following data,
\begin{list}{}{}
\item[(a)] decreasing $M_s$-invariant filtrations $F^\bullet$ and 
$\www F^\bullet$ with \eqref{7.56b} and the properties in lemma 
\ref{t7.10}.
\item[(b)] $(TERP(w))$-structures $(\www H\to \C,\nnn, H_\R,P)$ such that
$\OO(\www H)_0$ is generated by elementary sections.
\end{list}
\end{corollary}

For such a $(TERP(w))$-structure the conditions to have a 
$(trTERP(w))$-structure and thus a CV-structure or even a CV$\oplus$-structure
can be expressed in properties of the filtration $F^\bullet$.
This generalizes in the case $M=\{pt\}$ the results of chapter \ref{s3}
(because now $\UU$ is nilpotent, but not necessarily zero)
and contains the viewpoint of $(TERP(w))$-structures, which was not
considered in chapter \ref{s3}.

\begin{lemma}\label{t7.12}
Let $(\www H\to \C,\nnn, H_\R,P)$ be a $(TERP(w))$-structure such that
$\OO(\www H)$ is generated by elementary sections and let $F^\bullet$
be the corresponding filtration defined in \eqref{7.56a}.

(a) The $(TERP(w))$-structure is a $(trTERP(w))$-structure if and only if
$F^\bullet\hiin$ gives a Hodge structure of weight $w-1$ and 
$F^\bullet\hiie$ gives a Hodge structure of weight $w$.

(b) In that case let $(K=\www H_0,\kappa,h,g,\UU,\QQ)$ be the corresponding
CV-structure (on $M=\{pt\}$; therefore $D,C,\www C$ vanish)
and let $\hat H\to \P^1$ be the extension of $\www H$ constructed in lemma
\ref{t2.14} (e). 
The Hodge decomposition can be written as
\begin{eqnarray} \label{7.75}
\hiii &=& \bigoplus_{\beta\in \R}\hiii(\beta)\ ,\\
\hiii (\beta) &:=& F^{w+[-\beta]}\hiib \cap \
\oooo{F^{[\beta]}\hiii_{e^{2\pi i\beta}}}\ .\label{7.76}
\end{eqnarray}
Then the space of global holomorphic sections of $\hat H$ is
\begin{eqnarray} \label{7.77}
\pi_* \OO(\hat H) = \bigoplus_{\beta\in \R} \{ es(A,\beta)\ |\ 
A\in \hiii (\beta)\}\ .
\end{eqnarray}
This yields a canonical isomorphism $\Psi:\hiii \to K$, which is the 
composition of 
\begin{eqnarray} \label{7.78}
\hiii =\bigoplus_{\beta\in \R}\hiii (\beta)
\to \pi_*\OO(\hat H)\mbox{ \ with \ }
\psi_\beta: \hiii(\beta) \to \pi_*\OO(\hat H)
\end{eqnarray}
and the isomorphism $\pi_*\OO(\hat H)\to K$.

The isomorphism $\Psi$ maps the complex conjugation on $\hiii$ to $\kappa$.
The hermitian form $h^\infty:= h\circ (\Psi\times \Psi)$ on $\hiii$ is given
as follows. Define a bilinear form $S^\infty$ on $\hiii$ by 
\begin{eqnarray} \label{7.79}
S^\infty (A,B) := S((G^{(\alpha)})^{-1}A,(G^{(\beta)})^{-1}B) 
\end{eqnarray}
for $A\in \hiia, B\in \hiib$, $\alpha,\beta\in (0,1]$.
The forms $S^\infty$ and $h^\infty$ are monodromy invariant, and
$h^\infty$ satisfies
\begin{eqnarray} \label{7.80}
&& h^\infty : \hiii (\alpha)\times \hiii (\beta)\to 0 \mbox{ \ for \ }
\alpha\neq \beta\ ,\\
&& h^\infty: (\hiii)^{p,q}\times (\hiii)^{r,s}\to 0 \mbox{ \ for \ }
(p,q)\neq (r,s)\ ,\label{7.81}\\
&& h^\infty (A,B) = \frac{(-1)^p}{(2\pi i)^{w-1}}S^\infty (A,\oooo B) 
\mbox{ \ for \ } A,B\in (\hiin)^{p,w-1-p}\ ,\label{7.82}\\
&& h^\infty (A,B) = \frac{(-1)^p}{(2\pi i)^{w}}S^\infty (A,\oooo B) 
\mbox{ \ for \ } A,B\in (\hiie)^{p,w-p}\ .\label{7.83}
\end{eqnarray}
\end{lemma}

{\it Proof.}
(a) The map $\tau$ from lemma \ref{t2.14} (d) maps elementary sections
to elementary sections, in the following way,
\begin{eqnarray} \nonumber
\tau(es(A,\alpha)) 
&=& \tau\left(z\mapsto z^\alpha \exp(\log z \frac{-N}{2\pi i})A\right)\\
&=& \left(\frac{1}{\oooo z}\mapsto \oooo{z^{-w}z^\alpha \exp(\log z 
\frac{-N}{2\pi i})A } \right) \nonumber \\
&=& \left( z\mapsto z^{-\alpha+w} \exp(\log z \frac{-N}{2\pi i})
\oooo A\right) \nonumber\\
&=& es(\oooo A,-\alpha+w)\ .\label{7.83b}
\end{eqnarray}
Let $\hat H\to \P^1$ be the extension of $\www H\to \C$ constructed in 
lemma \ref{t2.14} (e). Because $\OO(\www H)$ is generated by elementary 
sections, the same holds for $\OO(\hat H|_{\P^1-\{0\}})$. Therefore also the 
space $\pi_*\OO(\hat H)$ is generated by elementary sections.
The formulas \eqref{7.63} and \eqref{7.83b} show \eqref{7.77},
with $\hiii (\beta)$ defined by \eqref{7.76}. 
Therefore $\hat H$ is a trivial bundle if and only if \eqref{7.75} holds.

(b) It rests to show $\Psi\circ (\oooo \cdot )=\kappa \circ \Psi$ and 
the formulas for $h^\infty$. The compatibility 
$\Psi\circ(\oooo \cdot )=\kappa \circ \Psi$ follows from the definition
of $\kappa$ in theorem \ref{t2.19} (b) and from \eqref{7.83b} and 
\eqref{7.75}--\eqref{7.77}.
Consider $A\in \hiii(\alpha),B\in \hiii (\beta)$; then
\begin{eqnarray} \label{7.84}
h^\infty (A,B) &=& h(\Psi(A),\Psi(B)) = g(\Psi(A),\kappa\Psi(B))\\
&=& z^{-w}\cdot P(es(A,\alpha)(z),\, es(\oooo B,-\beta+w)(-z))\ ,
\nonumber
\end{eqnarray}
and this vanishes for $\alpha\neq \beta$ because of \eqref{7.65}+\eqref{7.66}.
Suppose $\alpha=\beta$; then \eqref{7.51}+\eqref{7.52} and the 
$z$-sesquilinearity of $P$ show for $\alpha\notin \Z$
\begin{eqnarray} \label{7.85}
h^\infty (A,B) = \frac{(-1)^{w+[-\alpha]}}{(2\pi i)^{w-1}}
S\left( (G^{(\alpha-[\alpha])})^{-1}A, 
\, (G^{(1+[\alpha]-\alpha)})^{-1}\oooo B\right)
\ ,
\end{eqnarray}
and for $\alpha\in \Z$
\begin{eqnarray} \label{7.86}
h^\infty (A,B) = \frac{(-1)^{w+[-\alpha]}}{(2\pi i)^{w}}
S\left( (G^{(1)})^{-1}A, \, (G^{(1)})^{-1}\oooo B\right)
\ .
\end{eqnarray}
This shows \eqref{7.80}--\eqref{7.83}. The forms $S^\infty$ and 
$h^\infty$ are monodromy invariant because $S$ is monodromy invariant.
\hfill $\qed$

\begin{remarks}\label{t7.13}
(i) The form $S^\infty$ shares many properties with $S$; it is monodromy
invariant, nondegenerate, $(-1)^{w-1}$-symmetric on $\hiin$ and
$(-1)^w$-symmetric on $\hiie$. But for $N\neq 0$ it does not take real
values on $H_\R$, because then $G^{(\alpha)}$ does not respect $H_\R$.

(ii) The CV-structure in lemma \ref{t7.12} is a CV$\oplus$-structure 
if and only if $h^\infty$ is positive definite. The formulas
\eqref{7.82} and \eqref{7.83} have the same form as \eqref{3.9},
but with $S^\infty$ instead of $S$. Therefore for $N\neq 0$ the condition
to have a CV$\oplus$-structure does not mean that $F^\bullet$ gives
polarized Hodge structures on $\hiin$ and $\hiie$.
For $N\neq 0$ it does also not mean that $\www F^\bullet$ gives 
polarized Hodge structures on $\hiin$ and $\hiie$, because 
$G^{(\alpha)}\neq const.\cdot \oooo{G^{(1-\alpha)}}$ for 
$0<\alpha<1$ if $N\neq 0$. But the condition that $h^\infty$ is positive
definite is close enough to a polarizing condition. This will be used in
the proof of theorem \ref{t7.20}.
\end{remarks}

Finally, consider the case $N=0$. Then $G^{(\alpha)}=\Gamma(\alpha)\cdot \id$,
and it respects the real structure. The pairings $S^\infty$ and $S$ coincide
on $\hiil\times\hiii_{\oooo\lambda}$ up to the factor
$\frac{\pi}{\cos \pi(\alpha-\frac{1}{2})}$ where $\alpha\in (0,1)$ and
$e^{-2\pi i\alpha}=\lambda\neq 1$. They coincide on 
$\hiie\times \hiie$. Two filtrations $F^\bullet$ and $\www F^\bullet$
on $\hiii$ with \eqref{7.56b} are equal. The results from sections
\ref{s3.2} for $M=\{pt\}$ are recovered in the following corollary,
now including the point of view of $(TERP(w))$-structures.

\begin{corollary}\label{t7.14}
Let $(H\to \C^*,\nnn,H_\R,P)$ and $(\hiii,M_s,N,S)$ be as at the beginning
of section \ref{s7.2}. Suppose that $N=0$.

Let $(\www H\to \C,\nnn,H_\R,P)$ be a $(TERP(w))$-structure with 
$\www H|_{\C^*}=H$.

(a) The following conditions are equivalent.
\begin{list}{}{}
\item[$(\alpha)$] The sheaf $\OO(\www H)$ is generated by elementary sections.
\item[$(\beta)$] The pair $(\www H,\nnn)$ has a logarithmic pole at 0.
\item[$(\gamma)$] The endomorphism $\UU:\www H_0\to \www H_0$ is $\UU=0$.
\end{list}

(b) Corollary \ref{t7.11} gives a one-to-one correspondence between 
$(TERP(w))$-structures with $(\alpha)$--$(\gamma)$ and decreasing 
monodromy invariant filtrations
$\www F^\bullet$ on $\hiii$ as in lemma \ref{t7.10}.

(c) The correspondence restricts to a one-to-one correspondence between
$(trTERP(w))$-structures with $(\alpha)$--$(\gamma)$ and filtrations
$\www F^\bullet$ on $\hiii$ which give Hodge structures of weight
$w-1$ on $\hiin$ and of weight $w$ on $\hiie$ and which satisfy 
\eqref{7.61} and \eqref{7.62}.

(d) The CV-structures of $(trTERP(w))$-structures with $(\alpha)$--$(\gamma)$
are CV$\oplus$-structures if and only if the corresponding Hodge structures
on $\hiin$ and $\hiie$ are polarized.

(e) In the situation in (d) theorem \ref{t3.1} and lemma \ref{t3.4} give
polarized Hodge structures on $\www H_0=K$ of weight $w-1$ and $w$.
The isomorphism $\Psi:\hiii\to K$ maps the filtration $\www F^\bullet$ to the
filtration in lemma \ref{t3.4}, the pairing $S^\infty$ to the pairing in 
lemma \ref{t3.4} and the automorphism $M_s$ to 
$(-1)^we^{2\pi i\QQ}$.
\end{corollary}

{\it Proof.} Everything follows easily from the discussion above.
\hfill $\qed$

\begin{remarks}\label{t7.15}
In the case of a variation of polarized Hodge structures on a bundle on $M$
as in section \ref{s3.2} the corresponding $(trTERP(w))$-structure
has a pole of Poincar\'e rank 1 along $\nmmm$, but because of $\UU=0$
each pair $(\www H|_{\C\times \{t\}},\nnn)$ for $t\in M$ has a logarithmic pole
and determines a Hodge filtration $\www F^\bullet (t)$ on $\hiii$ as above.
The pole of Poincar\'e rank 1 for derivations $\nnn_X$ with $X\in \tm$
just rephrases the Griffiths transversality.
\end{remarks}

\subsection{Filtrations and logarithmic poles}\label{s7.5}

In this section the situation considered in sections \ref{s7.2}, \ref{s7.3}
and \ref{s7.4} is extended along a manifold $M$.
Now $H\to \csmmm$ is a flat complex vector bundle with connection $\nnn$
and monodromy with eigenvalues in $S^1$ and 
with $\nnn$-flat real subspace $H_\R$ and a $\nnn$-flat 
pairing $P$ as that in section \ref{s7.2} (c) for some integer $w$.
The vector space $\hiii$ and the objects
$\hiii_\R, M_s,N,S,G^{(\alpha)}$ are defined
as in section \ref{s7.2}. Also elementary sections $es(A,\alpha)$
on $\csmmm$ are defined as in section \ref{s7.2}. Such a section
satisfies $\nnn_X es(A,\alpha)=0$ for $X\in \tm$ by definition.
Here $X\in \tm$ is identified with the canonical lift to $\csmmm$.

Extensions of $H$ to vector bundles $\www H\to \cmmm$ with logarithmic pole 
along $\nmmm$ are simple objects. The sheaves $\OO(\www H)$ are generated
by elementary sections, and therefore they can be encoded in filtrations
on $\hiii$. More precisely, the following holds.

\begin{proposition}\label{t7.16}
There is a one-to-one correspondence between $(LEP(w))$-structures 
$(\www H\to \cmmm,\nnn,P)$ (definition \ref{t5.1} (c)) with 
$\www H|_\csmmm =H$ and decreasing monodromy invariant filtrations
$F^\bullet =\www F^\bullet$ with the properties in lemma \ref{t7.10}.

For any $t\in M$ one recovers $F^\bullet$ from $(H|_{\C\times \{t\}},\nnn)$
with \eqref{7.56a}. The sheaf $\OO(\www H)$ is generated by the elementary
sections in \eqref{7.74a}.
\end{proposition}

{\it Proof.}
One can also formulate a correspondence without the pairings $P$ and $S$.
This is well known. Precise statements and proofs can for example
be found in \cite[III.1.4]{Sab4}\cite[8.2]{He2}.
The extension to a correspondence with $P$ and $S$ follows from the
proof of lemma \ref{t7.10}.
\hfill $\qed$

\bigskip

An important step in the construction of a Frobenius manifold is to
extend a given $(TEP(w))$-structure, or better a $(TERP(w))$-structure
(definition \ref{t2.12}), to a $(trTLEP(w))$-structure (definition \ref{t5.5}).
The construction in singularity theory by K. Saito and M. Saito takes the
point of view of filtrations.

\begin{theorem}\label{t7.17}
Let $(\www H\to \cmmm,\nnn,H_\R,P)$ be a $(TERP(w))$-structure with 
$\www H|_{\csmmm}=H$ such that for some $t_0\in M$ the restricted
$(TERP(w))$-structure $(\www H,\nnn,H_\R,P)|_{\C\times \{t_0\}}$ is 
regular singular and the filtration $\www F^\bullet$ defined in 
\eqref{7.56b} for this $(TERP(w))$-structure
gives PMHS's on $\hiin$ and $\hiie$ of weights $w-1$ and $w$.

Then any opposite filtration (definition \ref{t7.18} below)
gives via proposition \ref{t7.16} (applying it to $\immm$ instead of 
$\nmmm$) rise to an extension of 
$\www H$ to a bundle $\hat H\to \pmmm$ such that 
$(\hat H,\nnn,P)$ is a $(trTLEP(w))$-structure in a neighborhood of $t_0\in M$.

The endomorphism $\Nu:K\to K$ of the corresponding Frobenius type structure
(theorem \ref{t5.7}) is semisimple and has eigenvalues 
$-\alpha_1+\frac{w}{2},...,-\alpha_\mu+\frac{w}{2}$, where $\alpha_1,...,
\alpha_\mu$ are the exponents of the restricted 
$(TERP(w))$-structure on $\C\times \{t_0\}$. One has for $\alpha\in \R$
\begin{eqnarray} \label{7.87}
\sum_{\beta\geq\alpha}
    \ker (\Nu-(-\beta+\frac{w}{2}\id):K_{t_0}\to K_{t_0})\\
=  \{ [a]\in K_{t_0}\ |\ a\in 
    \OO(\www H|_{\C\times \{t_0\}})_0\cap V^\alpha\}\ .\nonumber
\end{eqnarray}
\end{theorem}

{\it Proof.}
One has to adopt the proof of \cite[Theorem 7.16]{He2} and include the 
pairings $P$ and $S$, using arguments in the proof of lemma \ref{t7.10}.
\hfill $\qed$

\bigskip
Opposite filtrations exist because of Deligne's $I^{p,q}$ (lemma \ref{t7.3});
the space of all opposite filtrations is isomorphic to 
$\C^{N_{opp}}$ as algebraic manifold for some $N_{opp}\in \Z_{\geq 0}$
\cite[Lemma 10.21]{He2}.
The following definition is adapted to the situation here.

\begin{definition}\label{t7.18}
Let $(\hiii,\hiii_\R,M_s,N,S)$ be as above and let $\www F^\bullet$ be an
$M_s$-invariant filtration on $\hiii$ such that the restrictions to 
$\hiin$ and $\hiie$ give PMHS's of weight $w-1$ and $w$.
A filtration $U^\bullet$ on $\hiii$ is called {\it opposite to }
$\www F^\bullet$ if it is decreasing and monodromy invariant with
\begin{eqnarray} \label{7.88}
N(U^p)\subset U^{p+1}\ ,\\
S(U^p\hiin,\, U^{-w-p}\hiin)=0\ ,\label{7.89}\\
S(U^p\hiie,\, U^{-w+1-p}\hiie)=0\ ,\label{7.90}\\
\hiin = \bigoplus_p \www F^p\cap U^{-1-p}\cap\hiin\ ,\label{7.91}\\
\hiie = \bigoplus_p \www F^p\cap U^{-p}\cap \hiie\ .
\end{eqnarray}
\end{definition}

\subsection{$E$-orbits of regular singular $(TERP)$-structures}\label{s7.6}

We consider $(H\to \csmmm,\nnn,H_\R,P)$ and $(\hiii,\hiii_\R, M_s,N,S)$
and elementary sections $es(A,\alpha)$ on $\csmmm$ as in section 
\ref{s7.5}. 

Suppose for a moment $M=\{pt\}$. Because of lemma \ref{7.10} the filtration
$\www F^\bullet$ in section \ref{s7.3} of a regular singular 
$(TERP(w))$-structure $(\www H\to \C,\nnn,H_\R,P)$ has a good chance to be
the Hodge filtration of PMHS's on $\hiin$ and $\hiie$.

Theorem \ref{t7.5} says that PMHS's correspond to nilpotent orbits of
Hodge like filtrations. In theorem \ref{t7.20} a part of this correspondence
is generalized to regular singular $(TERP(w))$-structures.
Nilpotent orbits generalize to Euler field orbits.
Theorem \ref{t7.20} is the main result of chapter \ref{s7}.
Together with the conjectural inverse statement (remark \ref{t7.21})
it shows the close relationship between PMHS's and CV-structures.

Suppose now that $M\subset \C$ is a ball around $0\in \C$ with
coordinate $\rho$ on $\C$. 
{\it A $(TERP(w))$-structure on $M$ with Euler field}
is a $(TERP(w))$-structure $(\www H\to \cmmm,\nnn,H_\R,P)$ 
with $\www H|_\csmmm=H$ such that $E:=\frac{\paa}{\paa \rho }$ satisfies
$\UU=-C_E$. The vector field $E$ is called {\it Euler field}.

\begin{lemma}\label{t7.19}
(a) A $(TERP(w))$-structure on $M$ with Euler field 
$E=\frac{\paa}{\paa \rho}$ is determined by its restriction to $0\in M$
in the following way. Let
\begin{eqnarray} \label{7.93}
\gamma_\rho:H_{(e^{-\rho}z,0)} \to H_{(z,\rho)} \mbox{ \ for }
(z,\rho)\in \csmmm 
\end{eqnarray}
be the isomorphism by parallel transport along the orbit of $\zdz+E$.
Let $(z\mapsto \sigma (z,0)\in H_{(z,0)})$ be a section in 
$\OO(\www H|_{\C\times \{0\}})$ in a neighborhood of 0. Then
\begin{eqnarray} \label{7.94}
(z,\rho)\mapsto \sigma (z,\rho) := \gamma_\rho \sigma (e^{-\rho}z,0)\in 
H_{(z,\rho)}
\end{eqnarray}
is a section in $\OO(\www H)$ in a neighborhood of $\nmmm$.
Such sections generate $\OO(\www H)$.

(b) The $(TERP(w))$-structure in (a) extends canonically to $\C\supset M$.

(c) If the restriction of the $(TERP(w))$-structure in (a) to one parameter
$\rho\in M$ is regular singular then the restriction to any parameter
is regular singular. If the restriction to one parameter is generated
by elementary sections then the restriction to any parameter is 
generated by elementary sections.
\end{lemma}

{\it Proof.} Let $\www {\sigma_1},...,\www{\sigma_\mu}$ be an 
$\OO_{\cmmm,0}$-basis of sections of the germ $\OO(\www H)_0$ and let
$\sigma_1,...,\sigma_\mu$ be the sections on $\www H|_\csmmm$
in a neighborhood of $0\in\cmmm$ with \eqref{7.94} and with 
$\sigma_i(z,0)=\www{\sigma_i}(z,0)$ for $z\in\C^*$ close to 0.
Write $\www \sigma :=(\www{\sigma_1},...,\www{\sigma_\mu})$ and
$\sigma:=(\sigma_1,...,\sigma_\mu)$. Then
\begin{eqnarray} \label{7.95}
\sigma = \www \sigma \cdot A \mbox{ \ for some }A\in 
\Gl(\mu,\OO_\csmmm)
\end{eqnarray}
with $A|_{\rho=0}=\id$. The orbits of $\zdz+E$ are the sets 
$\{(z,\rho)\ |\ e^{-\rho}z=\mbox{const.}\}$.
By definition $\nnn_{\zdz+E}\sigma =0$, and because of $\UU=-C_E$
\begin{eqnarray} \label{7.96}
\nnn_{\zdz+E}\www \sigma = \www \sigma \cdot B \mbox{ \ for some }
B\in M(\mu\times \mu,\OO_{\cmmm,0})\ .
\end{eqnarray}
Therefore 
\begin{eqnarray} \label{7.97}
(\zdz+E)A=-B\cdot A \ .
\end{eqnarray}
This shows $A\in \Gl(\mu,\OO_{\cmmm,0})$ and part (a).
The parts (b) and (c) are now trivial.
\hfill $\qed$

\begin{theorem}\label{t7.20}
Let $(\www H\to \C\times \C,\nnn,H_\R,P)$ be a $(TERP(w))$-structure
on $M=\C$ (with coordinate $\rho$) with Euler field 
$E=\frac{\paa}{\paa \rho}$.
Suppose that the restrictions to fixed parameters $\rho\in M$ are
regular singular $(TERP(w))$-structures.

(a) Let $F^\bullet (\rho)$ and $\www F^\bullet (\rho)$ be the 
corresponding filtrations from \eqref{7.56a} and \eqref{7.56b}.
Then 
\begin{eqnarray} \label{7.98}
F^\bullet (\rho) &=& \exp(\frac{\rho}{2\pi i}N)F^\bullet (0)\ ,\\
\www F^\bullet (\rho) &=& \exp(\frac{\rho}{2\pi i}N)\www F^\bullet (0)\ .
\label{7.99}
\end{eqnarray}

(b) Suppose that $(\hiin,S,-N,\www F^\bullet(0))$ and 
$(\hiie,S,-N,\www F^\bullet(0))$ are PMHS's of weights $w-1$ and $w$.
Then a bound $b\in \R$ exists such that the restriction of the 
$(TERP(w))$-structure to $\{\rho\ |\ \Re\rho>b\}$ is a 
$(trTERP(w))$-structure with positive definite metric $h$.

(c) Consider the situation in (b). Let $\QQ$ be the endomorphism of the
corresponding CV$\oplus$-structure on $\{\rho\ |\ \Re\rho>b\}$.
Let $\alpha_1,...,\alpha_\mu$ be the exponents of the $(TERP(w))$-structure
at $\rho=0$. Then the eigenvalues of $\QQ$ tend to 
$-\alpha_1+\frac{w}{2},...,-\alpha_\mu+\frac{w}{2}$ for $\Re\rho\to \infty$.
\end{theorem}

\begin{remark}\label{t7.21}
Part (b) generalizes the direction `$\Longrightarrow$' in the correspondence
between PMHS's and nilpotent orbits in theorem \ref{t7.5} (a).
Probably also the inverse statement to part (b) is true.
It would generalize the direction `$\Longleftarrow$' in theorem \ref{t7.5} (a),
which is the most important consequence of Schmid's $SL_2$-orbit theorem.
We do not use it and do not prove it here.
\end{remark}

{\it Proof of theorem \ref{7.20}.}
(a) An elementary section $es(A,\alpha)$ on $H\to \csmmm$ with 
$A\in \hiia$  transforms as follows under the shift $\gamma_\rho$ from
\eqref{7.93},
\begin{eqnarray} \label{7.99b}
\gamma_\rho es(A,\alpha)(e^{-\rho}z,0) = 
e^{-\rho\alpha}\exp\left(\frac{\rho}{2\pi i}N\right) es(A,\alpha)(z,\rho)\ .
\end{eqnarray}
Therefore the filtrations $F^\bullet (\rho)$ and $\www F^\bullet (\rho)$
from \eqref{7.56a} and \eqref{7.56b} satisfy \eqref{7.98} and
\eqref{7.99}.

(b) and (c) The proof is divided into five parts. In the parts 
(I)--(III) parts (b) and (c) are proved for the case when the 
$(TERP(w))$-structures for fixed $\rho\in M$ are generated by elementary 
sections, in (IV) and (V) they are proved in the general case.

(I) By theorem \ref{t7.5} (a) the tuples $(\hiin,S,\www F^\bullet(\rho))$
and $(\hiie,S,\www F^\bullet(\rho))$ are pure polarized Hodge structures of 
weight $w-1$ and $w$, respectively, if $\Re \rho$ is sufficiently large.
But in order to apply lemma \ref{t7.12} (b) we need statements about 
$F^\bullet (\rho)$. We have to connect $\www F^\bullet (\rho)$ and 
$F^\bullet (\rho)$ by a one parameter family of filtrations.

With a new parameter $r\in [0,1]$ the automorphism
$G^{(\alpha)}:\hiia\to\hiia$ from \eqref{7.47} is generalized to
\begin{eqnarray} \label{7.100}
G^{(\alpha,r)}:=\sum_{k\geq 0}\frac{1}{k!}\Gamma^{(k)}(\alpha) 
\left(r\frac{-N}{2\pi i}\right)^k : \hiia\to \hiia\ ;
\end{eqnarray}
then $G^{(\alpha,1)}=G^{(\alpha)}$ and $G^{(\alpha,0)}=\Gamma(\alpha)\cdot\id$.
The filtrations $F^\bullet (\rho)$ are generalized to
\begin{eqnarray} \label{7.101}
F^\bullet_{\rho,r}\hiia := G^{(\alpha,r)} \www F^\bullet(\rho) \hiia
\end{eqnarray}
for $r\in [0,1]$ and $\alpha\in (0,1]$; then
$F^\bullet_{\rho,0}=\www F^\bullet(\rho)$ and 
$F^\bullet_{\rho,1}=F^\bullet(\rho)$.

All these filtrations define mixed Hodge structures on $\hiin$ and $\hiie$
with the weight filtration from $-N$, 
because they coincide on the quotients of the weight filtration.
The proof of theorem \ref{7.5} shows that a bound $b_1\in \R$ exists such that
for $\Re \rho>b_1$ they define pure Hodge structures of weight $w-1$
on $\hiin$ and of weight $w$ on $\hiie$.
(Here it is important that the first part of the proof of theorem
\ref{t7.5} does not use any pairings.)

With the parameter $r\in [0,1]$ the pairing $S^\infty$ on $\hiii$ which 
was defined in lemma \ref{t7.12} is generalized to pairings 
$S^{\infty,r}$ by formula \eqref{7.79} with $G^{(\alpha)}$ and
$G^{(\beta)}$ replaced by $G^{(\alpha,r)}$ and $G^{(\beta,r)}$.

For $r\in [0,1]$ and $\rho\in \C$ with $\Re \rho>b_1$ the pairing
$h^\infty$ in lemma \ref{t7.12} is generalized to pairings
$h^{\infty,r,\rho}$ by the formulas \eqref{7.80}--\eqref{7.83}
with $S^{\infty,r}$ instead of $S^\infty$ and 
with the Hodge decomposition of the filtration 
$F^\bullet_{\rho,r}$. Then $h^{\infty,0,\rho}$ is a positive definite
hermitian form for $\www F^\bullet (\rho)$.

One sees easily that all the forms $h^{\infty,r,\rho}$ for $\rho$ with
$\Re \rho>b_1$ are hermitian and nondegenerate. Because 
$h^{\infty,0,\rho}$ is positive definite, they all are positive definite.

As in section \ref{s7.4}, the filtration $F^\bullet(\rho)$ gives rise to a
$(TERP(w))$-structure which is generated by elementary sections.
Now lemma \ref{t7.12} shows that for $\rho$ with $\Re\rho>b_1$ this is a
$(trTERP(w))$-structure with positive definite metric.

This shows (b) if the $(TERP(w))$-structures $(\www H,\nnn,H_\R,P)|_\rho $
for $\rho\in \C$ are all (one is sufficient by lemma \ref{7.19} (c))
generated by elementary sections.

\bigskip

(II) Suppose that the $(TERP(w))$-structure at $\rho=0$ is generated by
elementary sections and that Deligne's $I^{p,q}(0)$ for the 
mixed Hodge structure $(\hiii,-N,F^\bullet(0))$ satisfy
$I^{q,p}(0)=\oooo{I^{p,q}(0)}$.

We will see that then the bound $b$ in (b) can be chosen as $b=0$ and that
the endomorphisms $\QQ,\UU$ and $\kappa\UU\kappa$ of the 
CV$\oplus$-structure can be calculated explicitely. As in the proof of
theorem \ref{t7.5} we choose a basis $\{v_a\}_{a\in A}$ of
$\bigoplus_{p,q}I^{p,q}_0(0)$ with \eqref{7.17}, $\oooo{v_a}=v_{k(a)}$,
\eqref{7.19}. We can even choose 
$v_a\in I^{p(a),q(a)}_0(0)\cap\hiii_{\lambda(a)}$. Then
$\lambda(k(a))=\oooo{\lambda(a)}$. Define
\begin{eqnarray} \label{7.102}
b_{j,k,n}:= (-1)^j \left({k\atop j}\right)\frac{1}{n(n-1)...(n-j+1)}
\end{eqnarray}
for $n\in \Z_{\geq 0}$, $0\leq k\leq n$, $0\leq j\leq k$. Define
for $a\in A$ and $p\in \Z$ with $0\leq p(a)-p\leq n(a)$
\begin{eqnarray} \label{7.103}
\varphi_{a,p}(\rho) := \exp\left( \frac{\rho}{2\pi i}N\right)
\sum_{j=0}^{p(a)-p} b_{j,p(a)-p,n(a)} \cdot
\left( \frac{\Re\rho}{\pi i}N\right)^j v_a\ .
\end{eqnarray}
For $a\in A$ and $p\in \Z$ without $0\leq p(a)-p\leq n(a)$ define
$\varphi_{a,p}(\rho):=0$.

{\bf Claim 1:} {\it For $\Re \rho>0$ the pairs $(\hiin,F^\bullet(\rho))$ and 
$(\hiie,F^\bullet(\rho))$ are pure Hodge structures of weight $w-1$ and $w$,
and }
\begin{eqnarray} \label{7.014}
\qquad F^p(\rho)\cap \oooo{F^{w-1-p}(\rho)}\cap \hiin = 
\bigoplus_{a:\lambda(a)\neq 1,\, 0\leq p(a)-p\leq n(a)} 
\C\cdot \varphi_{a,p}(\rho)\ ,\\
F^p(\rho)\cap \oooo{F^{w-p}(\rho)}\cap \hiie = 
\bigoplus_{a:\lambda(a)=1,\, 0\leq p(a)-p\leq n(a)} 
\C\cdot \varphi_{a,p}(\rho)\ ,
\label{7.105}
\end{eqnarray}

This is an easy consequence of \eqref{7.23}, \eqref{7.24} and the following
elementary claim. We leave the proof of this claim to the reader.

{\bf Claim 2:} {\it For $n\in \Z_{\geq 0}$, $0\leq k\leq n$,
\begin{eqnarray} \nonumber
&& e^x(1+b_1x+b_2x^2+...+b_kx^k)\in \C\oplus \C x\oplus ... \oplus \C x^{n-k}
\mod x^{n+1}\\
\qquad && \iff b_j =b_{j,k,n}\mbox{ \ for \ } 1\leq j\leq k\ .\label{7.106}
\end{eqnarray}  }

Also the proof of the following claim is left to the reader.

{\bf Claim 3:}
\begin{eqnarray} \nonumber
\left(\frac{\Re \rho }{\pi i}N\right)\varphi_{a,p} (\rho)
& = & (p(a)-p-n(a))\varphi_{a,p-1}(\rho) \\
&& + (n(a)-2p(a)+2p) \varphi_{a,p}(\rho)\nonumber \\
&& + (p(a)-p)\varphi_{a,p+1}(\rho) \ .\label{7.107}
\end{eqnarray}

For $a\in A$ and $p\in\Z$ define $\alpha(a,p)\in \R$ as the number such that
\begin{eqnarray} \label{7.108}
\lambda (a) = e^{-2\pi i\alpha(a,p)}\mbox{ \ and \ }w+[-\alpha(a,p)]=p\ .
\end{eqnarray}
The numbers $\alpha(a,p)$ for $a\in A$ and $p\in \Z$ with
$0\leq p(a)-p\leq n(a)$ form the exponents of the
$(TERP(w))$-structure at $\rho =0$.

Because of claim 1 and lemma \ref{t7.12}, the $(TERP(w))$-structure
$(\www H,\nnn,H_\R,P)$ restricts on $\{\rho\ |\ \Re \rho >0\}$ to a 
$(trTERP(w))$-structure. 
Let $\hat H\to \P^1\times \{\rho\ |\ \Re\rho>0\}$ be the extension of
$\www H|_{\C\times \{\rho\ |\ \Re\rho>0\}}$ constructed in lemma
\ref{2.14} (e). Because of \eqref{7.77}
\begin{eqnarray} \label{7.109}
\pi_*\OO(\hat H) = \bigoplus_{a\in A, p:\, 0\leq p(a)-p\leq n(a)} 
\OO_M\cdot es(\varphi_{a,p},\alpha(a,p))\ .
\end{eqnarray}
The endomorphisms $\QQ,\UU$ and $\kappa\UU\kappa$ on 
$K:= \hat H|_{\{0\}\times \{\rho\ |\ \Re\rho>0\}}$ lift canonically to $\hat H$
and $\pi_*\OO(\hat H)$; for $\sigma \in \pi_*\OO(\hat H)$ 
\begin{eqnarray} \label{7.110}
\nnn_\zdz \sigma = \eezz \UU\sigma - (\QQ-\frac{w}{2}\id)\sigma -
z\kappa\UU\kappa\sigma \ .
\end{eqnarray}
For $\sigma = es(\varphi_{a,p},\alpha(a,p))$ we obtain with claim 3
\begin{eqnarray} \nonumber
&&\nnn_\zdz es(\varphi_{a,p},\alpha(a,p)) = 
es((\alpha(a,p)\id + \frac{-N}{2\pi i}) \varphi_{a,p},\alpha(a,p))\\
&=& \frac{p(a)-p-n(a) }{-2\Re\rho}\cdot \eezz
      \cdot es(\varphi_{a,p-1},\alpha(a,p-1))  \nonumber\\
&+& \left( \alpha(a,p) + \frac{n(a)-2p(a)+2p}{-2\Re\rho}\right)
      \cdot es(\varphi_{a,p},\alpha(a,p))  \nonumber\\
&+& \frac{p(a)-p}{-2\Re\rho}\cdot 
      z\cdot es(\varphi_{a,p+1},\alpha(a,p+1))\ .  \label{7.111}
\end{eqnarray}
For $\Re \rho \to \infty$, the matrices of $\UU$ and $\kappa\UU\kappa$ with
respect to the basis of $\pi_*\OO(\hat H)$ in \eqref{7.109} tend to 0,
the matrix of $\QQ$ tends to the diagonal matrix with eigenvalues 
$-\alpha(a,p)+\frac{w}{2}$.
This shows part (c) of the theorem under the assumptions above.

\bigskip

(III) Suppose that the $(TERP(w))$-structure at $\rho =0$ is generated by
elementary sections. Because of (I) a bound $b\in \R$ exists such that 
for $\rho$ with $\Re\rho>b$ the filtrations $F^\bullet (\rho)$ and 
$\www F^\bullet(\rho)$ give pure Hodge structures on $\hiin$ and $\hiie$
and the $(TERP(w))$-structure at $\rho$ is a $(trTERP(w))$-structure.
It rests to estimate $\QQ$ for $\Re\rho \to \infty$.

As in the proof of theorem \ref{t7.5}, we choose a basis $\{v_a\}_{a\in A}$
of $\bigoplus_{p,q}I^{p,q}_0(0)$ with \eqref{7.17}--\eqref{7.19}
and with $v_a\in I^{p(a),q(a)}_0\cap \hiii_{\lambda(a)}$.
We define $b_{j,k,n}$, $\varphi_{a,p}(\rho)$ and $\alpha(a,p)$
by \eqref{7.102}, \eqref{7.103} and \eqref{7.108}.
We define as in \eqref{7.76} for $\rho$ with $\Re\rho>b$
\begin{eqnarray} \label{7.112}
\hiii(\beta,\rho) := F^{w+[-\beta]}(\rho) \hiib \cap 
\oooo{F^{[\beta]}(\rho)\hiii_{e^{2\pi i\beta}}}\ .
\end{eqnarray}
Then $\hiii = \bigoplus_\beta \hiii(\beta,\rho)$, and the numbers $\beta$
with multiplicities $\dim \hiii(\beta,\rho)$
are the exponents $\alpha_1,...,\alpha_\mu$.

{\bf Claim 4:}
{\it A bound $b_2\geq b$ exists such that for $\rho$ with $\Re\rho >b_2$
and any $a\in A,p\in \Z$ with $ 0\leq p(a)-p\leq n(a) $
the set $\hiii(\alpha(a,p),\rho)$ contains
a unique element $\psi_{a,p}(\rho)$ with 
\begin{eqnarray} \label{7.113}
\psi_{a,p}(\rho) = \varphi_{a,p}(\rho) + 
\sum_{b,q:\, \alpha(b,q)\neq \alpha(a,p)} 
c_{a,p,b,q}(\rho) \varphi_{b,q}(\rho)\ .
\end{eqnarray}
The coefficients $c_{a,p,b,q}(\rho)$ are of order 
$O((\Re \rho)^{-\varepsilon})$ 
for some $\varepsilon>0$. The elements $\psi_{a,p}(\rho)$ for fixed $\rho$ 
form a basis of $\hiii$.    }

It can be proved be refining the arguments in the proof of theorem \ref{t7.5}.
One uses that the coefficients of the matrix $B$ in the proof of theorem
\ref{t7.5} which give the leading contribution to its determinant
are the same in the special case $I^{q,p}=\oooo{I^{p,q}}$ and in the 
general case $I^{q,p}\neq \oooo{I^{p,q}}$. The details are left to the 
reader.

Because of lemma \ref{t7.10} $N(F^p(\rho))\subset F^{p-1}(\rho)$ and
$N(\hiii(\beta,\rho)) \subset \bigoplus_{\gamma-\beta\in \{-1,0,1\}}
\hiii(\gamma,\rho)$. 
Claim 3 and claim 4 show
\begin{eqnarray} \label{7.114}
\left(\frac{\Re\rho}{\pi i}N\right) \psi_{a,p}(\rho) &=& 
(p(a)-p-n(a))\psi_{a,p-1}(\rho) \\
&& + (n(a)-2p(a)+2p)\psi_{a,p}(\rho) \nonumber\\
&& + (p(a)-p) \psi_{a,p+1}(\rho) \nonumber\\
&& + \sum_{b,q:\, \alpha(b,q)-\alpha(a,p)\in \{-1,0,1\}}
\www c_{a,p,b,q}(\rho)\psi_{b,q}(\rho)\ ,\nonumber
\end{eqnarray}
and the coefficients $\www c_{a,p,b,q}(\rho)$ are of order
$O((\Re \rho)^{-\varepsilon})$.

One can now repeat the arguments in (II) after claim 3 with
$\psi_{a,p}(\rho)$ instead of $\varphi_{a,p}(\rho)$.
One obtains formula \eqref{7.111} with an additional sum on the
right hand side with coefficients of order 
$O((\Re \rho)^{-1-\varepsilon})$.
This shows part (c) if the $(TERP(w))$-structure at $\rho=0$ is 
generated by elementary sections.

\bigskip
(IV) 
Now part (b) will be proved in the general case.
We choose a basis $\{v_a\}_{a\in A}$
of $\bigoplus_{p,q}I^{p,q}_0(0)$ as in step (III)
and define $\varphi_{a,p}(\rho)$, $\alpha(a,p)$, $\hiii(\beta,\rho)$,
$\psi_{a,p}(\rho)$ and the bound $b_2\geq b$ as in step (III).

By step (II) and corollary \ref{t7.11}, for $\rho$ with $\Re\rho >b_2$
the filtration $F^\bullet (\rho)$ corresponds to a $(trTERP(w))$-structure
which is generated by elementary sections. Its sheaf of fiberwise
global sections on $\P^1\times \{\rho\ |\ \Re\rho>b_2\}$ is
\begin{eqnarray} \label{7.115}
\bigoplus_{a\in A,p:\, 0\leq p(a)-p\leq n(a)}
\OO_M\cdot es(\psi_{a,p},\alpha(a,p))\ .
\end{eqnarray}
We want to show that for $\rho$ with large $\Re \rho$ the 
$(TERP(w))$-structure $(\www H,\nnn,H_\R,P)|_\rho$ is so close to it
that it is also a $(trTERP(w))$-structure.

{\bf Claim 5:} {\it 
For $\rho$ with $\Re\rho>b_2$, $a\in A$, $p\in \Z$ with 
$0\leq p(a)-p\leq n(a)$, 
the germ $\OO(\www H|_{\C\times\{\rho\}})_0$
contains a unique element $\chi_{a,p}(\rho)$ with
\begin{eqnarray} \label{7.116}
\chi_{a,p}(\rho) &=& es(\psi_{a,p}(\rho),\alpha(a,p)) \\
&& +  \sum_{(b,q,k)\in B_1(a,p)} c_{a,p,b,q}^{(k)} (\rho)\cdot 
z^{-k}es(\psi_{b,q}(\rho),\alpha(b,q))\ ,\nonumber 
\end{eqnarray}
where
\begin{eqnarray} \label{7.117}
 B_1(a,p) = \{ (b,q,k)\ |\ b\in A,q\in \Z,k\in \Z_{\geq 1},
\alpha(b,q)-k>\alpha(a,p)\} .
\end{eqnarray}
The sets $B_1(a,p)$ are finite. 
The sections $\chi_{a,p}(\rho)$ form a $\Cz$-basis of the germ
$\OO(\www H|_{\C\times \{\rho\}})_0$. The coefficients 
$c_{a,p,b,q}^{(k)}(\rho)$ are of an exponential order 
$O(e^{-\varepsilon_2\cdot \Re\rho})$ for some $\varepsilon_2>0$.  }

\bigskip
{\it Proof of claim 5.}
An elementary section is the leading elementary section of a section in
$\OO(\www H|_{\C\times \{\rho\}})_0$ if and only if it is a linear
combination of sections $z^{k}es(\psi_{b,q}(\rho),\alpha(b,q))$
with $k\geq 0$. This follows from the definition of the filtration
$F^\bullet(\rho)$. Because of \eqref{7.73} the germ 
$\OO(\www H|_{\C\times \{\rho\}})_0$ is generated by sections which
are linear combinations of finitely many elementary sections.

Now one can inductively find a section as in \eqref{7.116}. 
The uniqueness is also clear; just consider the leading elementary
section of the difference of two such sections.
The set $B_1(a,p)$ is obviously finite.

Fix $\rho_0\in \C$ with $\Re\rho_0>b_2$. Because of lemma \ref{t7.19}
and formula \eqref{7.99b}, one $\Cz$-basis of 
$\OO(\www H|_{\C\times \{\rho\}})_0$ for $\rho$ with $\Re\rho>b_2$
is given by the sections
\begin{eqnarray} \label{7.118}
&&es(\exp\left(\frac{\rho-\rho_0}{2\pi i}N\right)\psi_{a,p}(\rho_0),
\alpha(a,p)) \\
&+& \sum_{(b,q,k)\in B_1(a,p)} 
e^{-(\rho-\rho_0)(\alpha(b,q)-k-\alpha(a,p))} \nonumber\\
&& \cdot c_{a,p,b,q}^{(k)}(\rho_0) \cdot z^{-k}
es(\exp\left(\frac{\rho-\rho_0}{2\pi i}N\right)
\psi_{b,q}(\rho_0),\alpha(b,q))\ .\nonumber
\end{eqnarray}
If one expresses the sections 
$\exp\left(\frac{\rho-\rho_0}{2\pi i}N\right)\psi_{a,p}(\rho_0)$ in terms of 
the sections $\psi_{a,p}(\rho)$ or vice versa, the coefficients which turn up
have at most polynomial order in $\Re\rho$. This follows from claim 4
and \eqref{7.103}. If one constructs the $\Cz$-basis $\chi_{a,p}(\rho)$
starting from the basis in \eqref{7.118}, the exponential terms
$e^{-(\rho-\rho_0)(\alpha(b,q)-k-\alpha(a,p))}$ in \eqref{7.118}
also dominate the coefficients of the $\chi_{a,p}(\rho)$.
This shows claim 5.
\hfill $\qed${}

\bigskip
Let $\hat H\to \P^1\times \{\rho\ |\ \Re\rho>b_2\}$ be the extension of the
$(TERP(w))$-structure $(\www H,\nnn,H_\R,P)$ on 
$\C\times \{\rho\ |\ \Re\rho> b_2\}$ which was constructed in lemma \ref{t2.14}
(e). For $\rho$ with $\Re\rho\to \infty$ the coefficients 
$c_{a,p,b,q}^{(k)}(\rho)$ tend exponentially fast to 0.
The bundle $\hat H|_{\P^1\times \{\rho\}}$ tends to the trivial bundle
on $\P^1\times \{\rho\}$  with global sections 
in \eqref{7.115}. Because triviality is an
open condition, a bound $b_3\geq b_2$ exists such that 
$\hat H|_{\P^1\times \{\rho\ |\ \Re\rho> b_3\}}$ is a trivial bundle.
This shows part (b) in the general case.
It rests to prove part (c).

\bigskip

(V) One can construct fiberwise global sections on 
$\hat H|_{\P^1\times \{\rho\ |\ \Re\rho> b_3\}}$  
for some $b_3\geq b_2$ starting with the sections
$z^k\chi_{a,p}$ for $k\geq 0$ and their images under the map $\tau$.
By an infinite iteration one easily sees the following.

{\bf Claim 6:}{\it 
For $\rho$ with $\Re \rho>b_3$ for some $b_3\geq b_2$ the space 
$\pi_*\OO(\hat H|_{\P^1\times \{\rho\}})$ is generated by unique sections
\begin{eqnarray} \label{7.119}
&& es(\psi_{a,p}(\rho),\alpha(a,p)) \\
&+& \sum_{(b,q,k)\in B_2(a,p)} \www c_{a,p,b,q}^{(k)}(\rho) \cdot z^{-k}
es(\psi_{b,q}(\rho),\alpha(b,q))\ ,\nonumber 
\end{eqnarray}
where 
\begin{eqnarray} \label{7.120}
  B_2(a,p)= \{(b,q,k)|\, b\in A,q\in \Z,k\in \Z-\{0\},
\alpha_1\leq \alpha(b,q)-k\leq \alpha_\mu\}.&&
\end{eqnarray}
The sets $B_2(a,p)$ are finite. The coefficients 
$\www c_{a,p,b,q}^{(k)}(\rho)$ are of exponential order
$O(e^{-\varepsilon_2\cdot \Re\rho})$.       }

As in step (II) and step (III) one applies $\nnn_\zdz$ to these
sections in order to determine $\QQ,\UU$ and $\kappa\UU\kappa$.
The leading coefficients are the same as in formula \eqref{7.111}.
There are additional terms with coefficients of order
$O((\Re\rho)^{-1-\varepsilon})$ from step (III) and additional terms with
coefficients of order $O(e^{-\varepsilon_2\cdot\Re\rho})$ from claim 6.
This shows part (c) in the general case and finishes the proof of 
theorem \ref{t7.20}.
\hfill $\qed$

\begin{lemma}\label{t7.22}
Let $(\www H\to \C\times \C,\nnn,H_\R,P)$ be a $(TERP(w))$-structure
on $M=\C$ (with coordinate $\rho$) with Euler field $E=\frac{\paa}{\paa \rho}$.
Let $\hat H\to \P^1\times \C$ be the extension of $\www H$ constructed
in lemma \ref{t2.14} (e).

Two bundles $\hat H|_{\P^1\times \{\rho_1\}}$ and 
$\hat H|_{\P^1\times \{\rho_2\}}$ are canonically isomorphic if 
$\Re\rho_1=\Re\rho_2$. 
The set $\{\rho\ |\ \hat H|_{\P^1\times \{\rho\}} \mbox{ is not trivial}\}$
is either $\C$ or empty or a discrete union of real lines
$\{\rho \ |\ \Re\rho = const.\}$.
\end{lemma}

{\it Proof.}
Consider $\rho_1=0$ and $\rho_2=ir$ with $r\in \R$.

Recall the maps $\gamma:\P^1\to \P^1$, $z\mapsto\frac{1}{\oooo z}$, 
$\gamma_\nnn$, $\tau_{real}$ and $\tau$ from lemma \ref{t2.14}.
The flow of $\Im (\zdz)$ on $\C^*$ commutes with $\gamma$. The
$\nnn$-parallel transport in $H\to \csmmm$ along the orbits of 
$\Im(\zdz)+\Im (E)$ commutes with $\gamma_\nnn$ and $\tau_{real}$.
It nearly commutes with $\tau$; more precisely, for $\rho=ir$
with $r\in \R$,
\begin{eqnarray} \label{7.121}
\tau \circ \gamma_{ir}=e^{irw}\cdot \gamma_{ir}\circ \tau\ .
\end{eqnarray}
Because of this and \eqref{7.94}, the maps $\gamma_\rho$ in \eqref{7.93}
extend for $\rho=ir$ with $r\in \R$ to isomorphisms
\begin{eqnarray} \label{7.122}
\hat H_{(e^{-ir}z,0)} \to \hat H_{(z,ir)} \mbox{ \ for }
(z,ir)\in \P^1\times i\R\ .
\end{eqnarray}
This shows the first part of the lemma. The second part follows from the first
part and from the fact that
the set $\{\rho\ |\ \hat H|_{\P^1\times \{\rho\}}$ is not trivial$\}$ is 
real analytic.
\hfill $\qed$

\begin{remarks}
(i) Let $(\www H\to \cmmm,\nnn,H_\R,P)$ be a $(TERP(w))$-structure
on $M=\C$ (with coordinate $\rho$) with Euler field $E=\frac{\paa}{\paa \rho}$.
It induces a $(TERP(w))$-structure on the quotient $M_{quot}:=M/2\pi i\Z$;
\begin{eqnarray} \label{7.123}
\C=M \to M_{quot}\cong \C^*,\ 
\rho \mapsto e^{-\rho}=: \zeta\ .
\end{eqnarray}
The deck transformation $\rho\mapsto \rho+2\pi i$ on $M$ 
is lifted with the maps $\gamma_{2\pi i+\rho}\circ \gamma_{\rho}^{-1}$ to
an automorphism of the $(TERP(w))$-structure on $M$.
This gives the $(TERP(w))$-structure on $M_{quot}$.
The coordinate $\zeta$ on $M_{quot}\cong \C^*$ is chosen such that 
$E$ is mapped to $-\zeta\paa_\zeta$. The flow of $\Re(-\zeta\paa_\zeta)$
tends to 0.
Also with this choice of the coordinate on $M_{quot}$, 
the monodromy on $\{z\}\times \C^*$ is equal to
the usual monodromy on $\C^*\times \{\zeta\}$.

(ii) Suppose that for any (one is sufficient) $\rho$ and $\zeta$ the
$(TERP(w))$-structures in (i) are regular singular.
Then $\UU=-C_E$ is nilpotent. Then the $(TERP(w))$-structure on 
$M_{quot}$ is {\it tame} in the sense of \cite[ch. 2]{Si2}.
The estimates for $\UU$ and $\kappa \UU\kappa$ in the proof of 
theorem \ref{t7.20} can also be obtained
from \cite[Theorem 1]{Si2}, but not the estimate for $\QQ$.
\end{remarks}

\clearpage

\section{The case of hypersurface singularities}\label{s8}
\setcounter{equation}{0}

\noindent
In section \ref{s8.1} a version of the construction of Frobenius manifolds
for hypersurface singularities is presented. 
A central step is the construction of a $(TERP)$-structure, 
using essentially oscillating integrals.
In section \ref{s8.2} this is extended to a construction of $tt^*$ geometry
and of CDV-structures. With chapter \ref{s6} and theorem \ref{t7.20} a part of 
a conjecture about the behaviour of the $tt^*$ geometry along the
real vector field $E+\oooo E$ is proved.
Some examples in section \ref{s8.3} illustrate this behaviour.

\subsection{Construction of Frobenius manifolds with oscillating 
integrals}\label{s8.1}

In the beginning of the 80ies K. Saito \cite{SK3}\cite{SK4} studied
the Gau{\ss}--Manin connection of hypersurface singularities and 
developed the notion of primitive forms. His work was completed by
M. Saito \cite{SM2} and resulted in a construction of Frobenius
manifolds. A detailed and simplified account is given in \cite{He2}.

In \cite{Sab3}\cite{Sab4} a modified version is outlined which uses
a Fourier--Laplace dual of the Gau{\ss}--Manin connection,
a $(TERP)$-structure in our notations. Two center pieces of this
construction are recalled in the theorems \ref{t5.12} and \ref{t7.17}.
But Sabbah considered the case of global functions on affine manifolds.
In this section we will sketch his version for the classical case of
function germs. The following data will be constructed in four steps.

Step 1: A function $f:(\C^{n+1},0)\to (\C,0)$, a semiuniversal unfolding
$F$, the base space $M$ as an F-manifold (definition \ref{t4.2}),
the $\mu$-constant stratum $S_\mu\subset M$; the integer $w$
of the earlier chapters will be $w=n+1$.

Step 2: The topological part $(H\to \csmmm,\nnn,H_\R,P)$ of a 
$(TERP(n+1))$-structure and the data $(\hiii,M_s,N,S)$ (section \ref{s7.2}).

Step 3: A $(TERP(n+1))$-structure $(\www H\to \cmmm,\nnn, H_\R,P)$
(definition \ref{2.12}) with $\www H|_\csmmm = H$, the bundle
$K:=\www H|_\nmmm$ with Higgs field $C^K$, endomorphism $\UU^K$ and
holomorphic metric $g^K$ (lemma \ref{t2.14}).

Step 4: A Hodge filtration $\www F^\bullet$ for PMHS's on $\hiin$ and
$\hiie$ as in section \ref{s7.3}, an extension of the $(TERP(n+1))$-structure
to a $(trTLEP(n+1))$-structure $(\www H^{(U^\bullet)}\to \pmmm,\nnn,P)$
with theorem \ref{t7.17}, a Frobenius manifold structure on $M$ with
theorem \ref{t5.12}.

\bigskip
{\bf Step 1:}
Let $f:(\C^{n+1},0)\to (\C,0)$ be a holomorphic function germ with an
isolated singularity at 0. Its Milnor number $\mu$ is the dimension of the
Jacobi algebra
$\OO_{\C^{n+1},0}/\left(\frac{\paa f}{\paa x_0},...,
\frac{\paa f}{\paa x_n}\right)$.

Let $F:(\C^{n+1}\times \C^\mu,0)\to (\C,0)$ be a semiuniversal unfolding.
That means that $F|_{(\C^{n+1}\times \{0\},0)} = f$ and that the derivations
$\frac{\paa F}{\paa t_i}|_{(\C^{n+1}\times \{0\},0)}$, $i=1,...,\mu$,
represent a basis of the Jacobi algebra of $f$;
here $(x,t)=(x_0,...,x_n,t_1,...,t_\mu)\in \C^{n+1}\times \C^\mu$.

One can choose a representative $F:\XX\to \Delta $ with
$\Delta = B^1_\eta = \{\tau\ |\ |\tau|<\eta\}\subset \C$,
$M=B^\mu_\theta\subset \C^\mu$ and $\XX=F^{-1}(\Delta)\cap 
(B^{n+1}_\varepsilon\times M)\subset \C^{n+1}\times \C^\mu$
for suitable small $\varepsilon, \eta, \theta>0$.
This $F$ is a good representative in the sense of \cite[ch. 2D]{Lo},
and the map
\begin{eqnarray}\label{8.1}
\varphi:\XX\to \Delta\times M, \ (x,t)\mapsto (F(x,t),t)\ ,
\end{eqnarray}
is a $C^\infty$-fibration of Milnor fibers outside of a discriminant
$\check \DD \subset \dmmm$. The critical space $C\subset \XX$ of
$\varphi$ is defined by the ideal 
$\left(\frac{\paa F}{\paa x_0},...,
\frac{\paa F}{\paa x_n}\right)$.
The discriminant is $\check \DD= \varphi(C)$. 
The projection $pr_{C,M}:C\to M$ is finite and flat of degree $\mu$.
The Kodaira--Spencer map
\begin{eqnarray}\label{8.2}
\aaa_1:\tm \to (pr_{C,M})_*\OO_C,\ \frac{\paa }{\paa t_i}
\mapsto {\frac{\paa F}{\paa t_i}}|_C
\end{eqnarray}
is an isomorphism of free $\OO_M$-modules of rank $\mu$, because the 
unfolding $F$ is semiuniversal.
Here $M=B^\mu_\theta$ has to be chosen sufficiently small. 
The multiplication on the right hand side of \eqref{8.2} induces a 
multiplication $\circ$ on $\tm$. 
With $e:=\aaa_1^{-1}(1|_C)$ and $E:=\aaa_1^{-1}(F|_C)$, the tuple
$(M,\circ,e,E)$ is an F-manifold with Euler field $E$ \cite[theorem 5.3]{He2}.
This will also follow from the $(TERP(n+1))$-structure in step 3 and 
from lemma \ref{t4.3}.
One should consider $F$ as a family of functions 
\begin{eqnarray}\label{8.3}
F_t:\XX_t\to \Delta , \mbox{ where } \XX_t:=\XX\cap \C^{n+1}\times \{t\}\ ,
\end{eqnarray}
parameterized by $t\in M$, with $F_0=f$.
The Kodaira--Spencer map gives for each tangent space $T_tM$ an isomorphism
\begin{eqnarray}\label{8.4}
(T_tM,\circ,E|_t)  \cong (\bigoplus_{x\in Sing(F_t)} \mbox{ Jacobi algebra of }
(F_t,x), \mbox{mult.},[F_t])\ .
\end{eqnarray}
The F-manifold is generically semisimple, because for generic $t$ the function
$F_t$ has $\mu$ $A_1$-singularities.

The $\mu$-constant stratum $S_\mu\subset M$ is the subvariety
\begin{eqnarray}\label{8.5}
S_\mu = \{ t\in M\ |\ Sing(F_t)=\{x\} \mbox{ for some point }
x\in F_t^{-1}(0)\}\ .
\end{eqnarray}
For singularities with $n\geq 2$ it is in general not smooth. Define
\begin{eqnarray}\label{8.6}
\UU:= E\circ :\tm \to \tm\ .
\end{eqnarray}
Then
\begin{eqnarray}\label{8.7}
S_\mu = \{t\in M\ |\ \UU:T_tM\to T_tM \mbox{ is nilpotent}\}\ .
\end{eqnarray}
The inclusion `$\subset$' is obvious because of \eqref{8.4}.
The inclusion `$\supset$' follows from the result of Gabrielov, Lazzeri 
and L\^e that the singularities of $F_t$ are not all in the same fiber
if $F_t$ has several singularities 
(cf. \cite[Theorem 12.1]{He2} and references there).

\bigskip
{\bf Step 2:}
The bundle $H\to \csmmm$ will be the dual bundle to a bundle of 
homology classes of Lefschetz thimbles.
Fix for a moment $(z,t)\in \csmmm$.
A Lefschetz thimble in $\XX_t$ with boundary in 
$F_t^{-1}\left(\eta\cdot\frac{z}{|z|}\right) \subset \XX_t$ is an 
$(n+1)$-cycle $\Gamma = \bigcup_{\tau\in \gamma([0,1])}\delta(\tau)$,
where $\gamma:[0,1]\to \oooo\Delta$ 
is a path with $\gamma(0)\in F_t(Sing(F_t))$,
$\gamma((0,1))\cap F_t(Sing(F_t))=\emptyset$, 
$\gamma(1)=\eta\cdot\frac{z}{|z|}$ and $\delta(\tau)$ for 
$\tau\in \gamma([0,1])$
is a continuous family of $n$-cycles in $F_t^{-1}(\tau)$ which vanish
for $\tau\to \gamma(0)$. For more details see \cite{Ph3}\cite{Ph4}.

If $t$ is generic then $F_t$ has $\mu$ $A_1$-singularities with different
critical values. In that case standard arguments 
(e.g. \cite[ch. 2]{AGV2}, \cite[(5.11)]{Lo}) show that one obtains
 a space homotopy equivalent to $\XX_t$ if one glues $\mu$ Lefschetz thimbles
to the fiber $F_t^{-1}(\eta\cdot\frac{z}{|z|})$.
Therefore for generic $t$
\begin{eqnarray}\label{8.8}
H^{Lef}_{(z,t),\Z} := H_{n+1}(\XX_t,F_t^{-1}(\eta\cdot\frac{z}{|z|}),\Z)
\cong \Z^\mu\ .
\end{eqnarray}
Because $\XX_t$ is contractible, the homology sequence of the pair 
$(\XX_t,F_t^{-1}(\eta\cdot\frac{z}{|z|}))$ shows that the map 
\begin{eqnarray}\label{8.9}
H^{Lef}_{(z,t),\Z} \to H_{n}(F_t^{-1}(\eta\cdot\frac{z}{|z|}),\Z),
\ [\Gamma] \to [\delta(\gamma(1))] 
\end{eqnarray}
is an isomorphism.

Because the pairs $(\XX_t,F_t^{-1}(\eta\cdot\frac{z}{|z|}))$ are 
diffeomorphic for all $(z,t)\in\csmmm$, \eqref{8.8}\and \eqref{8.9}
hold for all $(z,t)\in \csmmm$, and 
$\bigcup_{(z,t)\in \csmmm} H^{Lef}_{(z,t),\Z} $ is a $\Z$-bundle
of rank $\mu$. Therefore
\begin{eqnarray}\label{8.10}
H:= \bigcup_{(z,t)\in \csmmm} \Hom (H^{Lef}_{(z,t),\C},\C)
\end{eqnarray}
is a flat vector bundle with  a connection $\nnn$ and a real flat
subbundle $H_\R$. The monodromy $M_{mon}$ is the same as the monodromy of the
Milnor fibration of $f$, because of the isomorphism in \eqref{8.9}.
It is quasiunipotent. Semisimple and unipotent part are called $M_s$ and 
$M_u$, the nilpotent part is $N:=\log M_u$.

The intersection form for Lefschetz thimbles with boundaries in 
$F_t^{-1}(\eta\cdot\frac{z}{|z|})$ and in $F_t^{-1}(\eta\cdot\frac{-z}{|z|})$
gives a well defined nondegenerate pairing \cite{Ph4}
\begin{eqnarray}\label{8.11}
<\cdot,\cdot> : H^{Lef}_{(z,t),\Z}\times H^{Lef}_{(-z,t),\Z}\to \Z\ .
\end{eqnarray}
It induces an isomorphism $H^{Lef}_{(z,t),\Z}\to
\Hom(H^{Lef}_{(-z,t),\Z},\Z)$. With this isomorphism for $z$ and for $-z$
one obtains a pairing $<\cdot ,\cdot >^*$ on the dual spaces.
Define the pairing
\begin{eqnarray}\nonumber
P &:& H_{(z,t)}\times H_{(-z,t)}\to \C\ ,\\
&& (a,b)\mapsto (-1)^{\frac{n(n+1)}{2}}\frac{1}{(2\pi i)^{n+1}}
<a,b>^*\ .\label{8.12}
\end{eqnarray}
It is $\nnn$-flat, nondegenerate, $(-1)^{n+1}$-symmetric and takes values
in $i^{n+1}\R$ on $H_\R$.

Now the data $(\hiii,\hiii_\R,M_s,N,S)$ are defined as in section \ref{s7.2}.
The space $\hiii$ is the space of manyvalued global flat sections of the
bundle $H$, see \eqref{7.29}. The space $\hiii_\R\subset \hiii$ is the
real subspace of sections in $H_\R$. The induced monodromy operator
on $\hiii$ is also called $M_{mon}$ with semisimple part $M_s$
and nilpotent part $N$. The eigenspaces of $M_s$ are the 
$\hiil$; and $\hiin:=\bigoplus_{\lambda\neq 1}\hiil$.

A pairing $S$ on $\hiii$ is defined as in section \ref{s7.2} $(f)$.
It satisfies the properties in lemma \ref{t7.6} for $w=n+1$.
In \cite[ch. 10]{He2} another space $\hiii$ with a pairing $S$ was defined.
The spaces $\hiii$ are canonically isomorphic by the isomorphism in
\eqref{8.9}. Using topological arguments from \cite{Ph4} or \cite{AGV2}
one can show that the pairings $S$ coincide.

\bigskip
{\bf Step 3:}
One has to start with the Gau{\ss}--Manin connection of the map
$\varphi:\XX\to\dmmm$ in \eqref{8.1}. 
Detailed discussions from different points of view are given in
\cite{Mal1}, \cite{Gr}, \cite{Ph1}, \cite{SK3}\cite{SK4}, \cite{Lo},
\cite{Od}, \cite{SM2}, \cite{AGV2}, \cite{He2}.
The sketchy version here is close to \cite{He2}.

Work of Brieskorn on the Gau{\ss}--Manin connection for $f$ motivated the
definition of the following sheaf by Greuel and K. Saito.
The sheaf 
\begin{eqnarray}\label{8.13}
\HH^{(GM)} := \varphi_*\oomm^{n+1}_{\XX/M} / \ddd F\land \ddd \varphi_*
\oomm^{n-1}_{\XX/M}
\end{eqnarray}
is a coherent and even a free $\OO_\dmmm$-module of rank $\mu$
\cite{Gr}.
Let $H^{(GM)}$ be the corresponding vector bundle. 
The restriction of $H^{(GM)}$ to $\dmmm - \check \DD$ is canonically 
isomorphic to the cohomology bundle
\begin{eqnarray}\label{8.14}
\bigcup_{(\tau,t)\in \dmmm-\check\DD} H^n(F_t^{-1}(\tau),\C)\ .
\end{eqnarray}
The flat Gau{\ss}--Manin connection $\nnn^{(GM)}$ on this cohomology
bundle has a logarithmic pole along the discriminant 
$\check \DD$ with respect to the bundle $H^{(GM)}$ \cite{SK3}\cite{SK4}.

In order to carry out a Fourier--Laplace transformation in the sense of
\cite[V.2.c]{Sab4} one has to make things algebraic in the
$\tau$-direction, where now $\tau$ is the coordinate on $\Delta\subset\C$
used in \eqref{8.14}.
Because $\check \DD$ does not meet $(\paa\Delta)\times M$, one can extend
the bundle $H^{(GM)}\to\dmmm$ to a bundle
$\www H^{(GM)}\to \cmmm$ with flat connection outside of $\check \DD$.

The sheaf $\OO(\www H^{(GM)})$ is a free $\OO_\cmmm$-module of rank
$\mu$ and extends the sheaf $\HH^{(GM)}$.
It contains the sheaf $\OO^{mod}(\www H^{(GM)})$ of sections of 
$\www H^{(GM)}$ which have moderate growth along $\immm$.
Let $\pi_1:\cmmm\to M$ be the projection.
Then the sheaf $\HH^{alg,GM} := (\pi_1)_*\OO^{mod}(\www H^{(GM)})$
is a free $\OO_M[\tau]$-module of rank $\mu$. It has a logarithmic pole
along $\check \DD$ and a regular singularity along $\immm$.

Let $\pi_2:\dmmm\to M$ be the projection. Any element of 
$(\pi_2)_*\HH^{(GM)}$ has a unique preimage in $(\pi_2)_*\HH^{(GM)}$
under $\nnn^{(GM)}_{\paa_\tau}$ (\cite{SK3}\cite{SK4}\cite{Od}
\cite[Theorem 10.7]{He2}). Therefore $(\nnn^{(GM)}_{\paa_\tau})^{-1}$
acts on $(\pi_2)_*\HH^{(GM)}$ and also on $\HH^{alg,GM} $.

With the following Fourier--Laplace transformation one obtains an 
$\OO_M[z]$-module $\HH^{alg,osc} $ with a meromorphic connection
$\nnn^{(osc)}$ \cite[V.2]{Sab4}:
As sheaves $\HH^{alg,osc}=\HH^{alg,GM}$. Functions in $\OO_M$ and
vector fields $X\in \tm$ are lifted canonically to $\cmmm$.
Their actions on $\HH^{alg,osc} $ and $\HH^{alg,GM} $ coincide,
\begin{eqnarray}\label{8.15}
\nnn^{(GM)}_X = \nnn^{(osc)}_X\mbox{ \ for }X\in\tm\ .
\end{eqnarray}
The actions of $z$ and of $\nnn^{(osc)}_{\paa_z}$ on $\HH^{alg,osc}$
are given by
\begin{eqnarray}\label{8.16}
z &=& (\nnn^{(GM)}_{\paa_\tau})^{-1}\ ,\\
z^2 \nnn^{(osc)}_{\paa_z} &=& \tau \ .\label{8.17}
\end{eqnarray}
Then $\HH^{alg,osc}$ is a free $\OO_M[z]$-module of rank $\mu$ and 
$\nnn^{(osc)}$ is a flat connection on it with regular singularity
along $\immm$ and pole of Poincar\'e rank 1 along $\nmmm$, see
\cite[V.2.7]{Sab4} for the case without parameters
$t_1,...,t_\mu$. The order 1 pole for derivations 
$\nnn^{(osc)}_{\paa/\paa t_i}$ follows from formula \eqref{8.23} below.
The sheaf $\HH^{alg,osc}$ is the sheaf
$(\pi_1)_*\OO^{(mod)}(\www H)$ of fiberwise global sections with 
moderate growth along $\immm$ of a vector bundle $\www H\to\cmmm$
with flat connection $\nnn^{(osc)}$ on $\www H|_\csmmm$. 
We claim that there is a canonical isomorphism
\begin{eqnarray}\label{8.18}
(\www H|_\csmmm,\nnn^{(osc)})\cong (H,\nnn)
\end{eqnarray}
with $(H,\nnn)$ from step 2.
Formula \eqref{8.19} will make this isomorphism explicit.
Consider a section $\omega_{GM}\in \HH^{alg,GM}$, 
the corresponding section $\omega_{osc}\in \HH^{alg,osc}$, a point 
$(z,t)\in\csmmm$ and a Lefschetz thimble 
$\Gamma=\bigcup_{\tau \in \gamma([0,1])} \delta(\tau)$ with $\delta(\tau)$
and $\gamma:[0,1]\to \oooo\Delta$ as in step 2.
Then $[\Gamma]\in H_{n+1}(\XX_t,F_t^{-1}(\eta\cdot \frac{z}{|z|}),\Z)$
and we want to define a number $\omega_{osc}([\Gamma])$.

Extend $\gamma$ to a path $\gamma:[0,\infty)\to \C$ with $\gamma([1,\infty))$
the half line from $\gamma(1)$ to $\infty\cdot z$. Extend the homology bundle
which is dual to the cohomology bundle in \eqref{8.14} to a flat bundle
on $\cmmm-\check\DD$. Extend the family $[\delta(\tau)]$ of 
homology classes of cycles to a flat section in this bundle over
$\gamma([0,\infty))$. 
Define
\begin{eqnarray}\label{8.19}
\omega_{osc}([\Gamma]) :=\int_{\gamma([0,\infty))} e^{-\tau/z}
\omega_{GM}([\delta(\tau)])\ddd \tau\ .
\end{eqnarray}
This imitates {\it oscillating integrals}.
One has to show that it is well defined and that $\omega_{osc}$ gives
in this way a holomorphic section in $H$
which has moderate growth along $\immm$. Then it is easy to see
that this makes the Fourier--Laplace transformation above explicit.
We will not carry out the details here. Closely related considerations
(but without the extension to $\immm$) are presented in \cite{Ph3}\cite{Ph4}.

Also the pairing $P$ is treated there. In \cite[2i\`eme partie 4]{Ph4}
it is shown that the pairing $P$ on $\HH^{alg,osc}$ which is defined in 
\eqref{8.12} is the image under
the Fourier--Laplace transformation of K. Saito's higher residue pairings
(\cite{SK1}, cf. also \cite[10.4]{He2}). This implies that
$(\www H,\nnn,H_\R,P)$ is a $(TERP(n+1))$-structure.

Let us describe the bundle $K:= \www H|_\nmmm$, its Higgs field $C^K$,
the endomorphism $\UU^K$ and the pairing $g^K$ (lemma \ref{t2.14}).
The sheaf $\OO(K)$ is 
\begin{eqnarray}\nonumber
\OO(K) &\cong & \HH^{alg,osc}/z\HH^{alg,osc} \cong 
\HH^{alg,GM}/(\nnn^{(GM)}_{\paa_\tau})^{-1} \HH^{alg,GM} \\
&\cong & (\pi_2)_* \oomm^{n+1}_{\XX/\dmmm }=: \oomm_F\ ,\label{8.20}
\end{eqnarray}
where $\pi_2:\dmmm\to M$ is the projection 
\cite{SK3}\cite{SK4}\cite{He2}. Let us denote this isomorphism by
$\aaa_2:\OO(K)\to \oomm_F$.
The sheaf $\oomm_F$ is an $\OO_M$-module of rank $\mu$ and a free
$(pr_{C,M})_*\OO_C$-module of rank 1. With $\aaa_1$ as in \eqref{8.2},
the Higgs bundle structure on $K$ is given by
\begin{eqnarray}\label{8.21}
\aaa_2 (C^K_{\paa/\paa t_i} b) = -\aaa_1 (\frac{\paa}{\paa t_i})\cdot
\aaa_2(b)
\end{eqnarray}
for $b\in \OO(K)$. Therefore $\OO(K)$ is a free $\tm$-module of rank 1
and 
\begin{eqnarray}\label{8.22}
C^K_XC^K_Y = -C^K_{X\circ Y}\mbox{ and } C^K_e= -\id\ .
\end{eqnarray}
Formula \eqref{8.21} follows from the formula for the Gau{\ss}--Manin 
connection (e.g. \cite[Theorem 10.5]{He2}),
\begin{eqnarray}\label{8.23}
\nnn^{(GM)}_{\paa/\paa t_i} [\omega] = [\Lie_{\paa/\paa t_i}\omega]
-\nnn^{(GM)}_{\paa_\tau} [\frac{\paa F}{\paa t_i}\omega]
\end{eqnarray}
for $\omega\in \oomm^{n+1}_\XX$ with $[\omega]\in \HH^{(GM)}$.

The endomorphism $\UU^K$ is defined by the action of $z^2\nnn_\dz$
on $\OO(K)$ (lemma \ref{t2.14}). The definition of the Euler field,
\eqref{8.17} and \eqref{8.21} show
\begin{eqnarray}\label{8.24}
\UU^K = -C_E^K\ .
\end{eqnarray}
From Pham's result about the pairing $P$ and the higher residue pairings
and from their properties, it follows that the pairing $g^K$ on $K$
is mapped with $\aaa_2$ to Grothendieck's residue pairing 
on $\oomm_F$ (see e.g. \cite[10.4]{He2}).

\bigskip

{\bf Step 4:}
The pair $(\HH^{(GM)},\nnn^{(GM)})$ has a logarithmic pole along the 
discriminant $\check \DD$. Therefore the restriction to $t=0$ has a 
regular singularity at 0. The germ $\OO(H^{(GM)}|_{\Delta\times \{0\}})_0$
at 0 is called Brieskorn lattice and has been studied a lot.
Its properties are put together in \cite[10.6]{He2}.
They can be translated to properties of the restricted $(TERP(n+1))$-structure
$(\www H,\nnn,H_\R,P)|_{\C\times \{0\}}$. Some of them are the following:

{\bf $(\alpha)$} This $(TERP(n+1))$-structure is regular singular at 0;
therefore all the notions in sections \ref{s7.2}--\ref{s7.4}
can be used to discuss it.

{\bf $(\beta)$} Its exponents $\alpha_1,...,\alpha_\mu$ 
(defined in \eqref{7.70}) satisfy \cite{Mal1}
\begin{eqnarray}\label{8.25}
0<\alpha_1\leq ... \leq \alpha_\mu<n+1\ ,\ \alpha_i+\alpha_{\mu+1-i}=n+1\ .
\end{eqnarray}

{\bf $(\gamma)$} Define $F^\bullet$ and $\www F^\bullet$ for this
$(TERP(n+1))$-structure with \eqref{7.56a} and \eqref{7.56b}.
Then $(\hiin,M_s,-N,S,\www F^\bullet)$ is a PMHS of weight $n$ with 
automorphism $M_s$, and $(\hiie,-N,S,\www F^\bullet)$ is a PMHS of weight
$n+1$ \cite[Theorem 10.29]{He2} \cite[Theorem 3.5]{He1}
(the minus sign in $-N$ stems from the different definitions of the monodromy
here and in \cite{Sch}).

{\bf $(\delta)$} The smallest exponent has multiplicity one, i.e.
$\alpha_1<\alpha_2$, and for $\omega \in \OO(\www H|_{\C\times \{0\}})_0$ with
$[\omega]\in K_0$
\begin{eqnarray}\nonumber
&& \mbox{the map }C^K_\bullet [\omega] :T_0M\to K_0 \mbox{ is an isomorphism}\\
&&\iff \omega \in V^{\alpha_1}-V^{\alpha_2}\ .\label{8.26}
\end{eqnarray}
\cite[3.11]{SM3} (cf. \cite[Theorem 10.32]{He2})
($V^\alpha$ was defined in \eqref{7.45}).

Because of $(\gamma)$ one can apply theorem \ref{t7.17}.
The choice of any filtration $U^\bullet$ on $\hiii$ which is opposite
to $\www F^\bullet$ (definition \ref{7.18}) yields an extension
$\www H^{(U^\bullet)}\to \pmmm$ of $\www H\to \cmmm$ such that 
the tuple $(\www H^{(U^\bullet)},\nnn,H_\R,P)$ is a $(trTLEP(n+1))$-structure.

Let $(K,\nnn^r,C^K,g^K,\UU^K,\Nu^K)$ be the corresponding Frobenius type 
structure on $K$. By theorem \ref{t7.17} $\Nu^K$ is semisimple with
eigenvalues $-\alpha_1+\frac{n+1}{2},...,-\alpha_\mu+\frac{n+1}{2}$.
In order to obtain a Frobenius manifold one wants to apply theorem 
\ref{t5.12}. One needs a $\nnn^r$-flat section $\zeta\in \OO(K)$ such
that $C^K_\bullet \zeta:\tm \to \OO(K)$ is an isomorphism and 
$\zeta$ is an eigenvector of $\Nu^K$.
Because of \eqref{7.87} and $(\delta)$ above, the space of all such sections
is $\{\zeta\in \OO(K)\ |\ \nnn^r\zeta=0,\ \Nu^K\zeta = 
(-\alpha_1+\frac{n+1}{2})\zeta\}$, and it is one dimensional.
This recovers \cite[Theorem 11.1]{He2}.

\begin{theorem}\label{t8.1}
Each choice $(U^\bullet,\zeta)$ of an opposite filtration $U^\bullet$ and
a generator $\zeta$ of a one dimensional space gives a Frobenius manifold
structure $(M,\circ,e,E,g)$ on $M$. The flat endomorphism $\nnn^g_\bullet E$
of $TM$ is semisimple with eigenvalues $1+\alpha_1-\alpha_i$,
$i=1,...,\mu$; and $\Lie_E(g) = (2 -(\alpha_\mu-\alpha_1))g$.
The choice $(U^\bullet,c\zeta)$ for $c\in \C^*$ gives the Frobenius manifold
$(M,\circ,e,E,c^2g)$.
\end{theorem}

\subsection{$tt^*$ geometry in the singularity case}\label{s8.2}

In step 3 in section \ref{s8.1} a $(TERP(n+1))$-structure 
$(\www H\to \cmmm,\nnn,H_\R,P)$ was constructed on the base space $M$
of a semiuniversal unfolding $F$ of a singularity
$f:(\C^{n+1},0)\to (\C,0)$, using essentially oscillating integrals.
The results in chapter \ref{s6} show that it is generically a 
$(trTERP(n+1))$-structure.

\begin{theorem}\label{t8.2}
(a) The set $R$ of points in $M$ where the $(TERP(n+1))$-structure is not
a $(trTERP(n+1))$-structure is a real analytic subvariety $R\subset M$,
$R\neq M$. The set $R$ is invariant under the flow of $E-\oooo E$.
The restriction to $M-R$ of the bundle $K=\www H|_\nmmm$ carries a
canonical CV-structure.

(b) If $(M,\circ,e,E,g)$ is one of the Frobenius manifold structures
in theorem \ref{t8.1}, then on $M-R$ a real structure $\kappa$ exists
such that $(M-R,\circ,e,E,g,\kappa)$ gives a CDV-structure.
The hermitian metric $h$ and the endomorphism $\QQ$ of this 
CDV-structure are invariant under the flow of $E-\oooo E$.
\end{theorem}

{\it Proof.}
(a) The F-manifold $(M,\circ,e,E)$ is generically semisimple.
The Higgs bundle $(K=\www H|_\nmmm,C^K)$ is a free $TM$-module of rank 1.
Therefore the $(TERP(n+1))$-structure is at generic points locally isomorphic
to a universal semisimple $(TERP(n+1))$-structure in the sense of 
chapter \ref{s6}. Theorem \ref{t6.1} (a) and remark \ref{t6.2} (iii)
apply.
The set $R$ is invariant under the flow of $E-\oooo E$ because of
lemma \ref{t7.22}.

(b) This follows from theorem \ref{t5.15}, from the construction of the 
Frobenius manifold structure and from the definition of a CDV-structure.
\hfill $\qed$

\bigskip
In \cite{CV1}\cite{CV2} an important aspect of the study of the CV-structures
there is the behaviour with respect to the renormalization group flow.
Lemma \ref{t8.6} will show that it corresponds here to the real 
vector field $E+\oooo E$.

Especially interesting are the behaviour of the hermitian metric $h^K$
and the endomorphism $\QQ^K$ of the CV-structure on $K|_{M-R}$.
We formulate conjecture \ref{t8.3} and prove a part of it in 
theorem \ref{t8.4}. Part (b) of theorem \ref{t8.4} is the main result
of this chapter. It uses theorem \ref{t7.20}, which 
showed the relevance of PMHS's for CV-structures.

The restriction $F_t$ of the semiuniversal unfolding $F$ to a point
$t\in M$ has a finite set $Sing(F_t)$ of singularities.
For each singularity $x\in Sing(F_t)$ one has a tuple of exponents 
$Exp(F_t,x)\subset \Q\cap (0,n+1)$. They are defined precisely as the
exponents $\alpha_1,...,\alpha_\mu$ for $f$.
They are symmetric around $\frac{n+1}{2}$. Let
$Exp(F_t):=\bigcup_{x\in Sing(F_t)} Exp(F_t,x)$ be the union.
Along an orbit of $E$ this tuple is constant because there the singularities
do not split and because the exponents are constant along 
$\mu$-constant deformations \cite{Va2}.

\begin{conjecture}\label{t8.3}
Consider any $t\in M$. The set $R$ does not contain the $E+\oooo E$ orbit
of $t$. If one goes far enough along the flow of $E+\oooo E$ then one will
not meet anymore the set $R$, the hermitian metric will be positive
definite, and the eigenvalues of $\QQ^K$ will tend to 
$Exp(F_t)-\frac{n+1}{2}$.
\end{conjecture}

\begin{theorem}\label{t8.4}
(a) Conjecture \ref{t8.3} is true for points $t\in M$ such that $F_t$
has $\mu$ different critical values.

(b) Conjecture \ref{t8.3} is true for points $t$ in the $\mu$-constant
stratum $S_\mu$.
\end{theorem}

The points in (a) are the points $t$ such that $\UU|_t$ has $\mu$ different
eigenvalues, the points in (b) are the points $t$ such that $\UU|_t$ is 
nilpotent.

\bigskip
{\it Proof of theorem \ref{t8.4}.}
(a) The eigenvalues of $\UU$ are at generic points $t$ canonical coordinates
$u_1,...,u_\mu$ with $E=\sum_i u_i\frac{\paa }{\paa u_i}$. Along the
flow of $E+\oooo E$ their phases are constant and their absolute values
increase. Theorem \ref{t6.1} (b) applies.

(b) It is sufficient to prove the conjecture for the point 
$0\in S_\mu\subset M$.
At $t=0$ the $(TERP(n+1))$-structure from step 3 in section \ref{s8.1}
is regular singular and induces PMHS's
$(\hiin,M_s,-N,S,\www F^\bullet )$ and $(\hiie,-N,S,\www F^\bullet)$
of weight $n$ and $n+1$ (section \ref{s8.1}, step 4 $(\gamma)$).

The Euler field orbit of 0 is a point, i.e. $E|_0=0$, if and only if
the singularity is quasihomogeneous. Then $N=0$, and $\www F^\bullet$ gives
pure polarized Hodge structures. Corollary \ref{t7.14} (d) applies.

If $E|_0\neq 0$ then lemma \ref{t7.19} and theorem \ref{t7.20} apply.
\hfill $\qed$

\begin{remark}\label{t8.5}
The base space $M$ of a semiuniversal unfolding was constructed in 
section \ref{s8.1} as a small ball in $\C^\mu$.
Within $M$ one cannot go very far along the flow of $E$.

But if $E|_t\neq 0$ for all $t\in M$ then one can extend $M$ uniquely to 
a manifold $M_{orb}$ with Euler field $E$ such that each $E$-orbit
in $M_{orb}$ is isomorphic to $(\C,\frac{\paa}{\paa \rho})$ 
(with $\rho$ as coordinate on $\C$) and such that each $E$-orbit
in $M_{orb}$ meets $M$. Lemma \ref{t7.19} works also with additional
holomorphic parameters and shows that the $(TERP(n+1))$-structure
on $M$ extends to a $(TERP(n+1))$-structure on $M_{orb}$.

If $E|_t=0$ for some $t$ then 
\begin{eqnarray}\nonumber 
\{t\in M\ |\ E|_t=0\} = \{ t\in M\ |\ F_t \mbox{ is a quasihomogeneous
singularity}\}\ .
\end{eqnarray}
In that case one can choose new coordinates $\tau_i$ on $M$ such that
the Euler field takes the form
\begin{eqnarray}\label{8.27}
E=\sum_{i=1}^\mu d_i\tau_i\paa_{\tau_i}
\end{eqnarray}
where $d_i\in\Q_{\leq 1}$, and the numbers $1-d_i$ are the weighted
degrees of a monomial basis of the Jacobi algebra of a quasihomogeneous
singularity. Then the closure of $\tau(M)\subset \C^\mu$ under the flow
of $E$ is an open subset $M_{orb}\subset \C^\mu$ such that all $E$-orbits
with Euler fields are either points or are isomorphic to
$(\C^*,c\cdot \zeta\paa_\zeta)$ with $\zeta$ as coordinate on $\C$ and
$c\in \Q^*$ ($c$ may vary with the orbits). Again the 
$(TERP(n+1))$-structure on $\tau(M)$ extends to a $(TERP(n+1))$-structure
on $M_{orb}$. Therefore, in any case one can extend $M$ and the 
$(TERP(n+1))$-structure on $M$ arbitrarily along $E$-orbits.
\end{remark}

A semiuniversal unfolding $F$ of a singularity is a family of functions
$F_t$, $t\in M$. In \cite{CV1}\cite{CV2} an additional parameter
$\lambda$ and one parameter unfoldings
$\lambda\cdot F_t$, $\lambda\in \C^*$ (or $\R^*$), are considered
and related to the renormalization group flow.
Lemma \ref{t8.6} states the relation with the Euler field $E$.

\begin{lemma}\label{t8.6}
Choose in the base space $M$ of a semiuniversal unfolding $F$
a parameter $t$. Consider the one parameter unfolding 
$e^\rho\cdot F_t$, $\rho\in \C$, of $F_t$.

A neighborhood $U\subset \C$ of $\rho=0$ and a unique map $u:U\to M$
with $u(0)=t$ exist such that $u^*F$ and the unfolding
$e^\rho\cdot F_t$, $\rho \in  U$, are isomorphic in the following 
sense: a holomorphic family of (multi-germ) isomorphisms $\varphi_\rho$
of suitable neighborhoods of the critical points
of $u^*F|_\rho$ and $e^\rho\cdot F_t$ exists such that $u^*F|_\rho 
=e^\rho\cdot F_t \circ \varphi_\rho$ in these neighborhoods 
and such that $\varphi_0=\id$.

The image $u(U)\subset M$ is an $E$-orbit, and 
$\ddd u(\frac{\paa}{\paa \rho})=E$ on $u(U)$.
\end{lemma}

{\it Proof.}
The first part follows from the fact that the unfolding $F$ is for any 
$t\in M$ a semiuniversal unfolding of $F_t$ as a multi-germ
around its singular points. The set $u(U)\subset M$ is an $E$-orbit with
$\ddd u(\frac{\paa}{\paa \rho})=E$ because of the calculation
\begin{eqnarray}\label{8.28}
\frac{\paa}{\paa \rho} \left( e^\rho\cdot F_t\right) 
= e^\rho\cdot F_t
\end{eqnarray}
and because of the Kodaira--Spencer map in \eqref{8.2}.
\hfill $\qed$

\subsection{Examples}\label{s8.3}

In \cite{CV1}\cite{CV2} quasihomogeneous singularities are considered.
A semiuniversal unfolding of a quasihomogeneous singularity $f\in 
\C[x_0,...,x_n]$ can be chosen as
\begin{eqnarray}\label{8.29}
F=f+\sum_{i=1}^\mu t_im_i\ ,
\end{eqnarray}
where $m_1,...,m_\mu$ are weighted homogeneous polynomials which 
represent a basis of the Jacobi algebra of $f$.
The $m_i$ have weighted degrees $\deg_w m_i \in \Q_{\geq 0}$, where
$\deg_w f=1$, the $t_i$ inherit weighted degrees 
$d_i:=\deg_w t_i := 1-\deg_w m_i$.
The Euler field is 
\begin{eqnarray}\label{8.30}
E=\sum_{i=1}^\mu d_it_i\paa_{t_i}\ .
\end{eqnarray}
In \cite{CV1}\cite{CV2} a deformation $F_t$ of $f$ is called
{\it relevant} if $t_i=0$ for $d_i\leq 0$,
{\it marginal} if $t_i=0$ for $d_i\neq 0$,
{\it irrelevant} if $t_i=0$ for $d_i\geq 0$.

In \cite{CV1}\cite{CV2} mainly deformations $F_t$ with $t_i=0$ for $d_i<0$
are considered. One reason is that such deformations have as global
functions on $\C^{n+1}$ still only the singularities which come
from the singularity 0 of $f=F_0$, so that
$\sum_{x\in \C^{n+1}}\mu(F_t,x) =\mu$.
Deformations $F_t$ with $t_i\neq0$ for some $d_i<0$ have in general
additional singularities coming from infinity, so that
$\sum_{x\in \C^{n+1}}\mu(F_t,x) >\mu$.

The deformations $F_t$ with $t_i=0$ for $d_i>0$ are the $\mu$-constant
deformations. The results in chapter \ref{s7} and theorem 
\ref{t8.4} (b) apply precisely to these deformations and thus form
a complement to \cite{CV1}\cite{CV2}. Also, they apply to
nonquasihomogeneous singularities, as well.

Theorem \ref{t8.4} (b) will be illustrated by the simplest examples
of singularities with nontrivial $\mu$-constant stratum.
There are three classes of singularities with one dimensional
$\mu$-constant stratum:\\
(A) the simple elliptic singularities $\www {E_6},\www {E_7},\www {E_8}$;\\
(B) the hyperbolic singularities $T_{pqr}$ for 
$p,q,r\in \Z_{\geq 2}$ with $\frac{1}{p}+\frac{1}{q}+\frac{1}{r}<1$;\\
(C) the 14 exceptional unimodal singularities.\\
There are also three classes of singularities with two dimensional
$\mu$-constant stratum. One consists of \\
(D) the 14 exceptional bimodal singularities.

The singularities in (A) are all quasihomogeneous. We consider the 
singularities in (B), (C) and (D) and describe their
$\mu$-constant strata $S_\mu$, the subsets $R\cap S_\mu$ of points 
where the $(TERP(n+1))$-structure is not a $(trTERP(n+1))$-structure,\
the metric $h$ and the endomorphism $\QQ$ of the CV-structure on 
the complement $S_\mu-R$. In \cite{He-1}\cite{He0} 
(see also \cite{Ku}) the Gau{\ss}--Manin connections for these singularities
were studied carefully. We use the results.

In the following $S_\mu$ will be the parameter space of a standard
$\mu$-constant family $\www F_t$, $t\in S_\mu$, which contains representatives
of each right equivalence class of singularities of the given 
topological type. In (B) it is not simply connected.
We work with surface singularities ($n=2$), but that is not important.

\bigskip
(B) The hyperbolic surface singularities $T_{p,q,r}$
for $p,q,r\in \Z_{\geq 2}$ with
$\frac{1}{p}+\frac{1}{q}+\frac{1}{r}<1$ form a one parameter family
\begin{eqnarray}\label{8.31}
\www F_t = x_0^p+x_1^q+x_2^r+t\cdot x_0x_1x_2,\ t\in S_\mu=\C^*\ ,
\end{eqnarray}
with $\mu=p+q+r-1$ and exponents $\alpha_1,...,\alpha_\mu$ with 
$1=\alpha_1<\alpha_2\leq ...\leq \alpha_{\mu-1}<\alpha_\mu=2$.
The Euler field on $S_\mu$ is 
\begin{eqnarray}\label{8.32}
E= (1-\frac{1}{p}-\frac{1}{q}-\frac{1}{r})t\paa_t \ .
\end{eqnarray}
The nilpotent endomorphism $N$ of $\hiii$ vanishes on $\hiin$, but
it has a $2\times 2$-Jordan block on the two dimensional space
$\hiie$. For each $t\in S_\mu$ the $(TERP(n+1))$-structure at $t$
is generated by elementary sections.
It splits into one $(TERP(n+1))$-structure for $\hiin$ and one for $\hiie$.

The one for $\hiin$ is constant along $S_\mu$, it is a 
$(trTERP(n+1))$-structure,
its metric $h$ is positive definite, its endomorphism $\QQ$ has eigenvalues
$\alpha_2-\frac{3}{2},...,\alpha_{\mu-1}-\frac{3}{2}$ everywhere.

Using the calculations in \cite{He-1}\cite{He0} and doing some more,
one obtains the following.
The $(TERP(n+1))$-structure for $\hiie$ is not a $(trTERP(n+1))$-structure
on the set
\begin{eqnarray}\label{8.33}
S_\mu \cap R=\{t\in S_\mu\ |\ |t|=r_1\}
\end{eqnarray}
for some $r_1\in \R_{>0}$. The metric $h$ of its CV-structure on 
$S_\mu-R$ is positive definite on $\{t\in S_\mu\ |\ |t|>r_1\}$
and negative definite on $\{t\in S_\mu\ |\ |t|<r_1\}$.
Its endomorphism $\QQ$ is on $S_\mu-R$ semisimple with eigenvalues
\begin{eqnarray}\label{8.34}
\pm \left( \frac{1}{2} + \left(1-\frac{1}{p}-\frac{1}{q}-\frac{1}{r}\right)
\frac{1}{2} \cdot \frac{1}{\log |t|-\log r_1}\right) \ .
\end{eqnarray}

\bigskip
(C) The 14 exceptional unimodal surface singularities
form each a one parameter family
\begin{eqnarray}\label{8.35}
\www F_t = f+t\cdot m\ ,\ \ t\in S_\mu=\C\ ,
\end{eqnarray}
where $f$ is quasihomogeneous and $m$ is a monomial of weighted 
degree $\deg_w m>1$. The exponents satisfy
\begin{eqnarray}\label{8.36}
\alpha_1<1,\ 1+(1-\alpha_1)<\alpha_2,\ \alpha_{\mu+1-i}=3-\alpha_i  \ .
\end{eqnarray}
The Euler field $E$ on $S_\mu$ is 
\begin{eqnarray}\label{8.37}
E= -2(1-\alpha_1)t\paa_t\ .
\end{eqnarray}
The $(TERP(n+1))$-structure splits into two $(TERP(n+1))$-structures,
one for $\sum_{j=2}^{\mu-1} \hiii_{e^{-2\pi i\alpha_j}}$, the other
for the two dimensional space
$\bigoplus_{j\in \{1,\mu\}} \hiii_{e^{-2\pi i\alpha_j}}$.
The first one has the same properties as the first one in (B).

With the calculations in \cite{He-1}\cite{He0} and some more one finds
the following. The second one is not a $(trTERP(n+1))$-structure on 
\begin{eqnarray}\label{8.38}
S_\mu \cap R=\{t\in S_\mu\ |\ |t|=r_2\}
\end{eqnarray}
for some $r_2\in \R_{>0}$. The metric $h$ of its CV-structure on 
$S_\mu-R$ is positive definite on $\{t\in S_\mu\ |\ |t|<r_2\}$
and negative definite on $\{t\in S_\mu\ |\ |t|>r_2\}$.
Its endomorphism $\QQ$ is on $S_\mu-R$ semisimple with eigenvalues
\begin{eqnarray}\label{8.39}
\pm \frac{\alpha_1-\frac{3}{2}-|t|^2(\frac{1}{2}-\alpha_1) }
{1-|t|^2} \ .
\end{eqnarray}

\bigskip
(D) The 14 exceptional bimodal surface singularities form each a two 
parameter family 
\begin{eqnarray}\label{8.40}
\www F_t = f+t_1\cdot m_1+t_2\cdot m_2\ ,\ \ t\in S_\mu=C^2\ ,
\end{eqnarray}
where $f$ is quasihomogeneous and $m_1$ and $m_2$ are monomials
of weighted degrees $>1$. The exponents satisfy
\begin{eqnarray}\label{8.41}
\alpha_1<1,\ 1<\alpha_2< 1+(1-\alpha_1) <\alpha_3,
\ \alpha_{\mu+1-i}=3-\alpha_i \ .
\end{eqnarray}
The Euler field $E$ on $S_\mu$ is 
\begin{eqnarray}\label{8.42}
E= -2(1-\alpha_1)t_1\paa_{t_1}+ (-2+\alpha_1+\alpha_2) t_2\paa_{t_2}\ .
\end{eqnarray}
The $(TERP(n+1))$-structure splits into two $(TERP(n+1))$-structures,
one for $\sum_{j=3}^{\mu-2} \hiii_{e^{-2\pi i\alpha_j}}$, the other
for the four dimensional space 
$\bigoplus_{j\in \{1,2,\mu-1,\mu\}} \hiii_{e^{-2\pi i\alpha_j}}$.
The first one has the same properties as the first ones in (B) and (C).
Its endomorphism $\QQ$ has eigenvalues $\alpha_j-\frac{3}{2}$,
$j\in\{3,...,\mu-2\}$. 

With the calculations in \cite{He-1}\cite{He0} and some more one finds
the following. The second one is not a $(trTERP(n+1))$-structure on 
\begin{eqnarray}\label{8.43}
S_\mu \cap R=\{(t_1,t_2)\in S_\mu\ |\ 
  r_3|t_1+r_5t_2^{r_6}|^2-(1-r_4|t_2|^2)^2 =0\}
\end{eqnarray}  
for some $r_3,r_4\in \R_{>0}$, $r_5\in \R$, $r_6=\deg_wt_2/\deg_w t_1
\in \{4,3,\frac{5}{2}\}$ and $r_5=0$ for $r_6=\frac{5}{2}$. 
The set $S_\mu\cap R$ has singularities. It splits $S_\mu-R$ into 
three connectivity components.
I did not determine the metric $h$ and the endomorphism $\QQ$ precisely.
It is clear that the metric $h$ is positive definite on the component
which contains $t=0$, and that the eigenvalues of $\QQ$ tend to 
$\alpha_j-\frac{3}{2}$ for  $j\in\{1,2,\mu-1,\mu\}$  for $t\to 0$.

\clearpage


\bigskip
\centerline{---------------------}

\centerline{
Max-Planck-Institut f\"ur Mathematik,
Vivatsgasse 7, 53111 Bonn, Germany}

\centerline{hertling\char64 mpim-bonn.mpg.de}

\end{document}